\documentclass[preprint,12pt]{elsarticle}
\pdfoutput=1

\usepackage{hyperref}
\usepackage{graphicx}
\graphicspath{ {./graphics/} }
\usepackage{caption}
\usepackage{subcaption}
\usepackage{multirow}
\usepackage{amssymb}
\usepackage{amsthm}
\usepackage{amsmath}
\usepackage{algorithm} 
\usepackage{algpseudocode} 

\usepackage{lineno}

\let\oldequation\equation
\let\oldendequation\endequation
\renewenvironment{equation}
  {\linenomathNonumbers\oldequation}
  {\oldendequation\endlinenomath}

\def\mysep{\!\!\!&}
\def\Wi{W^{1,\infty}}
\def\bx{\boldsymbol{x}}
\DeclareMathOperator{\sgn}{sgn}

\newcommand{\dmat}[2]{\mathcal{#1}=(\boldsymbol{#2}_{ij})}
\newcommand{\dvec}[2]{\mathcal{#1}=(\boldsymbol{#2}_{i})}
\newcommand{\Rd}[2]{\mathbb{R}^{#1\times #2,d}}
\newcommand{\Rdvec}[1]{\mathbb{R}^{#1,d}}
\newcommand{\no}[1]{}

\newcommand{\boxedgraphics}[1]{\fbox{\includegraphics[width=\textwidth]{#1}}}
\newtheorem*{definition}{Definition}
\newcommand{\comment}[1]{}
\journal{Elsevier}

\begin{document}
\begin{frontmatter}

\title{Accurate Stabilization Techniques for RBF-FD Meshless Discretizations with Neumann Boundary Conditions}

\author[inst1]{Riccardo Zamolo}\corref{correspondingauthor}
\cortext[correspondingauthor]{Corresponding author}
\ead{rzamolo@units.it}

\affiliation[inst1]{
	organization={Department of Engineering and Architecture, University of Trieste, via Alfonso Valerio 10, 34127},
    city={Trieste},
    country={Italy}}

\author[inst1]{Davide Miotti}
\author[inst1]{Enrico Nobile}

\begin{abstract}
A major obstacle to the application of the standard Radial Basis Function-generated Finite Difference (RBF-FD) meshless method is constituted by its inability to accurately and consistently solve boundary value problems involving Neumann boundary conditions (BCs). This is also due to ill-conditioning issues affecting the interpolation matrix when boundary derivatives are imposed in strong form.
In this paper these ill-conditioning issues and subsequent instabilities affecting the application of the RBF-FD method in presence of Neumann BCs are analyzed both theoretically and numerically. The theoretical motivations for the onset of such issues are derived by highlighting the dependence of the determinant of the local interpolation matrix upon the boundary normals. Qualitative investigations are also carried out numerically by studying a reference stencil and looking for correlations between its geometry and the properties of the associated interpolation matrix. Based on the previous analyses, two approaches are derived to overcome the initial problem. The corresponding stabilization properties are finally assessed by succesfully applying such approaches to the stabilization of the Helmholtz-Hodge decomposition.
\end{abstract}



\end{frontmatter}


\section{Introduction}\label{s:Intro}
Some of the main issues with traditional CAE (Computer Aided Engineering) solvers, i.e. those based on Finite Volume (FVM) \cite{ver07intro} and Finite Elements (FEM) methods \cite{zien77finite}, are caused by the limited flexibility of the mesh data structure and the complexity of its generation \cite{liu05intro}.
In order to leverage the increased CPU capabilities and allow the user to focus on other aspects of the simulation process, meshless methods \cite{li13meshless,ngu08meshless,patel20meshless,liu09meshfree,duar95review,slak2021medusa} might be adopted as a replacement for the more traditional mesh based ones.

Before they can be considered as an improvement over the traditional ones, however, meshless methods must guarantee reliable results also in presence of generic geometries, e.g., arbitrarily complex or affected by deformations, and boundary conditions, e.g., Neumann boundary conditions in particular. This would allow the user to skip the mesh generation process altogether without sacrificing the accuracy or the stability of the solution procedure.

In the present paper the  Radial Basis Function - Finite Difference (RBF-FD) collocation method augmented with polynomials is employed, a brief discussion on the main features is presented in the following sections, however an interested reader is encouraged to also see \cite{mio21fully}.

The discretization of the domain is usually achieved by placing a number of nodes both on the boundary and on the interior, when the RBF-FD collocation method is applied there is no need to store connectivity information and partial differential equations can be solved in strong form directly \cite{liu05intro}.

In practical cases the basic RBF-FD method has to be enhanced in order to ensure good convergence close to the boundaries \cite{fornberg2002observations}. Among others, the most popular improvement is polynomial augmentation which allows to reach higher convergence rates and to improve interpolation near boundaries as long as Dirichlet boundary conditions (BCs) are considered \cite{bay19insight,flyer2016role,bayona2017role,bay2019role}. The accuracy of this approximation can therefore be controlled by properly increasing the polynomial order and the number of nodes.

Neumann BCs can be enforced both at the stencil level, i.e., when the local interpolants exactly satisfy the BCs at the boundary nodes \cite{Mavric_2017,mavrivc2017application,mavrivc2015local,mio21fully,ZAMOLO2020109730}, or at the assembly level, i.e., when the equations for the BCs at the boundary nodes appear explicitly in the final sparse matrix \cite{Slak_phd,Janmono}. In both cases the RBF-FD collocation method, however, may lead to large errors and/or ill-conditioning problems. These issues constitute a major obstacle for the robust and accurate solution of boundary value problems of practical interest. Some effective solutions for this last problem have been suggested, however no one has been universally adopted yet.

In \cite{liu2006stabilized,tom2021least,tom2021unfitted} the stabilization is achieved through different least squares procedures. The one in \cite{tom2021unfitted} seems especially capable of providing both stable and accurate results by using two different node sets.

Another approach to the stabilization is the one based on the so called ghost nodes. In \cite{fly2016enhancing} one layer of ghost nodes is placed outside the domain at a distance roughly continuing the pattern of the interior nodes in order to address the issues related to one-sided boundary stencils. In \cite{lin2020radial} the authors propose instead the use of fictitious nodes that extend far outside the domain, resembling the approach from the method of fundamental solutions \cite{fairweather1998method}.
In another implementation \cite{fed02improved} the authors add a set of nodes, which can lie inside or outside of the domain, adjacent to the boundary and, correspondingly, add an additional set of equations obtained via collocation of the PDE on the boundary.
The strategy of ghost nodes is successfully employed also in \cite{de2022parallel}, for the stable RBF-FD solution of elliptic PDEs problems on complex 3D geometries and even in higher dimension in \cite{Janmono}.

Alternatively to the ghost nodes, in \cite{shankar2017overlapped} additional nodes are placed close to the boundary, but within the domain, in order to reduce the error in those areas in the context of the overlapped RBF-FD method. The overlapped variant differs from the traditional RBF-FD in the selection of the interpolation stencils and can be intended as a generalization of the latter.

Yet another approach is adopted in \cite{chen2006meshless,Mavric_2017}, where every time a Neumann boundary node has to be included in the stencil, its normal is checked according to different geometric criteria.

The aim of the present paper is to present some remarks on the roots of those issues related to the enforcement of Neumann BCs at the stencil level and to suggest new improvements of the traditional interpolation scheme. It is important to remark that the final implementation must be capable of dealing with any possible shape of the boundary, as a consequence it is not possible to make any a priori assumption on the direction of the normals, i.e., the outward facing vectors defined at each boundary node and orthogonal to the local tangent plane.

Building on top of the simplest RBF-FD scheme, the main source of errors and stability issues is found in the poor conditioning of the local interpolation matrix in presence of Neumann BCs. We found that, given a certain node distribution, it is always possible to have normals that make the local interpolation matrix singular. Such singular directions can be calculated as the ones that send the determinant of the said matrix to zero, therefore a control can be implemented in order not to include the corresponding boundary nodes in the stencil. This approach, somewhat similar to the one adopted in \cite{Mavric_2017}, is more expensive but totally robust and it is tested with good results. At the best of the authors' knowledge it represents the only possible solution that carries on the traditional RBF-FD method totally unchanged and doesn't require any intervention on the node generation. The presented mathematical formulation leaves great freedom for further developments aimed at improving both accuracy and stability of the RBF-FD scheme furthermore.

Alternatively, we found that it is possible to properly move the nodes on the boundary in order to rule out the possibility of singular interpolation matrices. This operations appears to be the most computationally efficient and the simplest to implement on smooth and regular domains, although requiring non-trivial geometrical operations on generic domains. However, a robust implementation of this approach on complex-shaped domains is not simple and would require a proper modification of the node generation algorithm. 

The proposed strategies for the prevention of ill-conditioning problems due to Neumann BCs, which doesn't require human intervention, paves the way for the implementation of robust RBF-FD meshless methods in combination with fully automated node generation.
\section{RBF-FD overview}\label{s:RBFFD}
Follows a brief description of the RBF-FD method which constitutes the starting point from which all the variants discussed above are developed.

The actual solution of a boundary value problem can be split in $3$ steps:
\begin{enumerate}
    \item node generation;
    \item discretization through the RBF-FD method;
    \item solution procedure for the resulting sparse linear system.
\end{enumerate}

\subsection{Node generation}
The employed node generation procedure consists of two phases:
\begin{enumerate}
    \item Generation of a node distribution within the domain $\Omega$ that satisfies a defined spacing function $s(\bx)$ on average. The boundary $\partial \Omega$ is not taken into account yet.
    \item Iterative refinement of the initial node distribution. At this stage the already existing nodes can be moved and projected onto $\partial \Omega$ appropriately.
\end{enumerate}

Implementation details fall outside the topic of this paper, more details can be found in \cite{ZAMOLO20184305}. Throughout this paper we will assume that nodes are evenly distributed, this can be attained with an adequate number of refinement iterations. This allows us to assume a negligible difference between fill distance and separation distance \cite{de2010stability,de2005near}. As a consequence such a difference can reasonably be ignored as a source of interpolation instabilities.

\subsection{RBF-FD method}\label{ss:RBF_FD}
{
\renewcommand{\arraystretch}{1.5}
\begin{table} \footnotesize
    \centering
    \begin{tabular}{l l l l}
        Smooth & Type & Definition $\Phi(r)$ & Positive definite, order\\
        \hline
        Infinitely & Gaussian (GA) &  $e^{-(\varepsilon \, r)^2}$ & strictly\\
        & Multiquadric (MQ) & $\sqrt{1 + (\varepsilon \, r)^2}$  & conditionally, 1\\
        & Inverse Multiquadric (IMQ) & $1/\sqrt{1 + (\varepsilon \, r)^2}$ & strictly\\
        & Inverse Quadratic (IQ) & $1/(1 + (\varepsilon \, r)^2)$ & strictly\\
        \hline
        Piecewise  & Monomial PHS &  $r^{2k+1}$, $k \in \mathbb{N}$ & conditionally, $k+1$\\
         & Thin Plate Spline PHS & $r^{2k}\log r$, $k \in \mathbb{N}$ & conditionally, $k +1$\\
    \end{tabular}
    \caption{most common Radial Basis Functions}
    \label{tab:RBFs}
\end{table}
}

In the simplest version of this method, any real valued function can be interpolated at a generic point $\bx$ as a linear combination of $m$ Radial Basis Functions (RBFs) $\Phi(||\bx-\bx_i||_2)$, and the set $\{\bx_i\}_1^m$ is the interpolation stencil, composed by distinct nodes. The most commonly used RBFs are reported in Table \ref{tab:RBFs}, where they are defined in terms of $r = ||\bx-\bx_i||_2$, for $\bx_i$ fixed.

By using the notation $\varphi_i(\bx) = \Phi(||\bx-\bx_i||_2)$, the interpolation $u^h(\bx)$ of a function $u(\bx)$ thus takes the following form:
\begin{equation}
    u^h(\bx) = \sum_{i=1}^m \alpha_i \varphi_i(\bx)
    \label{eq:bareboneRBF}
\end{equation}

In practical cases the previous interpolation has to be enhanced in order to ensure appropriate convergence properties \cite{fornberg2002observations}. A very popular improvement is polynomial augmentation, where the bare RBF interpolation scheme in equation \eqref{eq:bareboneRBF} is modified\comment{as follows: $u^h(\bx)$ is expanded also with respect to} by including the polynomial basis $\{p_j\}_1^q$ of the multivariate polynomial space $\Pi_P^d$ of total degree $P$ in $d$ dimensions \cite{bay19insight} as follows:
\begin{equation}
    u^h(\bx) = \sum_{i=1}^m \alpha_i \varphi_i(\bx) + \sum_{j=1}^q \beta_j p_j(\bx)
    \label{eq:polyRBF}
\end{equation}

We assume the nodes in the interpolation stencil are listed in the following order: the first $m_I$ are internal nodes, i.e., nodes lying within $\Omega$, while the following $m_B$ are boundary nodes, i.e., those lying on $\partial\Omega$, and $m_I+m_B=m$.

The following interpolation conditions at the internal nodes hold:
\begin{equation}
    u^h(\bx_i) = u(\bx_i) \, , \quad i = 1,\dots, m_I
    \label{eq:interp_cond}
\end{equation}

Since we will focus only on Neumann BCs, i.e., $\partial u/ \partial\boldsymbol{n} = g$, we will use the following notation for the normal derivative of a generic function $F$ at the boundary node $\bx_k$:
\begin{equation}
    \partial F(\bx_k) := \frac{\partial F}{\partial \boldsymbol{n}_k}(\bx_k)
    \label{eq:normal_der_notation}
\end{equation}
where $\boldsymbol{n}_k$ is the unit normal at $\bx_k$.

The interpolation conditions \eqref{eq:interp_cond} are then supplemented with Neumann BCs at the boundary nodes as follows:
\begin{equation}
    \partial u^h(\bx_k) = \partial u(\bx_k) = g(\bx_k)\, , \quad k = m_J,\dots, m
    \label{eq:NeumannBC}
\end{equation}
where $m_J=m_I+1$. 

In the case of polynomial augmentation, the following orthogonality conditions are enforced \cite{fasshauer2007meshfree}:
\begin{equation}
    \sum_{i=1}^m \alpha_i p_j(\bx_i) = 0 \, , \quad j=1,\dots,q
    \label{eq:ort_cond}
\end{equation}

By collecting equations \eqref{eq:interp_cond} and \eqref{eq:NeumannBC}-\eqref{eq:ort_cond} in matrix form, the following linear system is attained:
\begin{equation}
    \underbrace{\begin{bmatrix}
        \boldsymbol{\varphi}_{BC} & \mathbf{P}_{BC}\\
        \mathbf{P}^T & \boldsymbol{0}
    \end{bmatrix}}_{\displaystyle\mathbf{M}}
    \begin{Bmatrix}
            \boldsymbol{\alpha}\\
            \boldsymbol{\beta}
    \end{Bmatrix}
    =
    \begin{Bmatrix}
            \boldsymbol{u}\\
            \boldsymbol{g}\\
            \boldsymbol{0}
    \end{Bmatrix}
    \label{eq:interpSystem}
\end{equation}
where $\boldsymbol{\alpha}=\{ \alpha_i \}_{1}^{m}$, $\boldsymbol{\beta}=\{ \beta_j \}_{1}^{q}$, $\boldsymbol{u} = \{ u(\bx_i) \}_{1}^{m_I}$ and $\boldsymbol{g} = \{ g(\bx_k)\}_{m_J}^{m}$ are column vectors.
The matrix $\mathbf{M}$ in equation \eqref{eq:interpSystem} will be called the interpolation matrix and its blocks are expanded as follows:
\begin{equation}
    \boldsymbol{\varphi}_{BC}=
    \begin{bmatrix}
\varphi_1(\bx_1)               \mysep \cdots \mysep \varphi_m(\bx_1)               \\
\vdots                                    \mysep \ddots \mysep \vdots                                    \\
\varphi_1(\bx_{m_I})           \mysep \cdots \mysep \varphi_m(\bx_{m_I})           \\
\partial\varphi_1(\bx_{m_J}) \mysep \cdots \mysep \partial\varphi_m(\bx_{m_J}) \\
\vdots                                    \mysep \ddots \mysep \vdots                                    \\
\partial\varphi_1(\bx_m)       \mysep \cdots \mysep \partial\varphi_m(\bx_m)
    \end{bmatrix}
    \label{eq:phimat}
\end{equation}

\begin{equation}
    \mathbf{P}_{BC}=
    \begin{bmatrix}
p_1(\bx_1)                 \mysep \cdots \mysep p_q(\bx_1)               \\
\vdots                                \mysep \ddots \mysep \vdots                              \\
p_1(\bx_{m_I})             \mysep \cdots \mysep p_q(\bx_{m_I})           \\
\partial p_1 (\bx_{m_J}) \mysep \cdots \mysep \partial p_q(\bx_{m_J})\\
\vdots                                \mysep \ddots \mysep \vdots                              \\
\partial p_1(\bx_m)        \mysep \cdots \mysep \partial p_q(\bx_m)
    \end{bmatrix}
    \label{eq:pmat}
\end{equation}

\begin{equation}
    \mathbf{P}^T=
    \begin{bmatrix}
p_1(\bx_1) \mysep \cdots \mysep p_1(\bx_m) \\
\vdots                \mysep \ddots \mysep \vdots                \\
p_q(\bx_1) \mysep \cdots \mysep p_q(\bx_m) 
    \end{bmatrix}
    \label{eq:ptmat}
\end{equation}

The normal derivatives of RBFs in equation \eqref{eq:phimat} can be expressed as follows:
\begin{equation}
    \partial \varphi_i (\bx_k) = \frac{\partial \varphi_i}{\partial \boldsymbol{n}_k}(\bx_k) = 
    \Phi'(r_{i,k})\,\boldsymbol{e}_{i,k} \cdot \boldsymbol{n}_k
    \label{eq:normal_phi}
\end{equation}
where $r_{i,k} = ||\bx_k-\bx_i||_2$ and $\boldsymbol{e}_{i,k} =(\bx_k-\bx_i)/r_{i,k} $ is the unit direction from node $i$ towards node $k$. For the problem of interest, i.e., solution of PDEs, we also require the RBFs to satisfy equation \eqref{eq:PhiCondition} for $\varphi_i$ to be differentiable also at the nodes:
\begin{equation}
    \Phi'(0)=0
    \label{eq:PhiCondition}
\end{equation}
therefore, the diagonal entries $\partial\varphi_i(\bx_i)$ of the last $m_B$ rows of $\boldsymbol{\varphi}_{BC}$ are $0$.

Given a linear partial differential equation $\mathfrak{L}(u) = f$ in the unknown field $u$, the Weighted Residual Method is applied in combination with the Collocation technique as explained in \cite{liu05intro}.
The linear operator $\mathfrak{L}$ is therefore applied to $u^h$ and the following equation is made valid at each internal node $\bx$:
\begin{equation}
    \mathfrak{L}u^h(\bx) = f(\bx)
    \label{eq:collocation}
\end{equation}

When we substitute the definition of $u^h$ given in equation \eqref{eq:polyRBF} into equation \eqref{eq:collocation}, the linear operator can be applied to the basis functions as follows:
\begin{equation}
    \sum_{i=1}^m \alpha_i \, \mathfrak{L}\varphi_i(\bx) + \sum_{j=1}^q \beta_j \, \mathfrak{L}p_j(\bx) = f(\bx)
    \label{eq:collocation2}
\end{equation}

Equation \eqref{eq:collocation2} can be recast in vector form as follows, where $\boldsymbol{\Psi}(\bx) = \{ \mathfrak{L}\varphi_i(\bx) \}_1^m$ and $\boldsymbol{\Pi}(\bx) = \{ \mathfrak{L}p_j(\bx) \}_1^q$ are column vectors:
\begin{equation}
    \begin{Bmatrix}
            \boldsymbol{\Psi}(\bx)\\
            \boldsymbol{\Pi}(\bx)
    \end{Bmatrix}^T
    \begin{Bmatrix}
            \boldsymbol{\alpha}\\
            \boldsymbol{\beta}
    \end{Bmatrix}
    = f(\bx)
    \label{eq:colloCompact}
\end{equation}

By substituting vectors $\boldsymbol{\alpha}$ and $\boldsymbol{\beta}$ from equation  \eqref{eq:interpSystem} into equation \eqref{eq:colloCompact} we obtain:
\begin{equation}
    \underbrace{\begin{Bmatrix}
            \boldsymbol{\Psi}(\bx)\\
            \boldsymbol{\Pi}(\bx)
    \end{Bmatrix}^T
    \mathbf{M}^{-1}
    }_{\displaystyle\boldsymbol{c}^T}
     \begin{Bmatrix}
            \boldsymbol{u}\\
            \boldsymbol{g}\\
            \boldsymbol{0}
    \end{Bmatrix}
    =
    f(\bx)
    \label{eq:finalStencil}
\end{equation}

We point out that, instead of calculating the inverse matrix $\mathbf{M}^{-1}$, we can solve the following adjoint system for the coefficients $\boldsymbol{c}$:
\begin{equation}
    \mathbf{M}^T
    \boldsymbol{c}
    =
    \begin{Bmatrix}
        \boldsymbol{\Psi}(\bx)\\
        \boldsymbol{\Pi}(\bx)
    \end{Bmatrix}
    \label{eq:adjoint}
\end{equation}

Given an interpolation stencil built around point $\bx$, the vector $\boldsymbol{c}$ represents therefore the stencil coefficients for the RBF-FD approximation of the linear operator $\mathfrak{L}$ at point $\bx$. Only the first $m$ elements of $\boldsymbol{c}$ are kept, i.e., those multiplying nonzero elements in equation \eqref{eq:finalStencil}. These elements of $\boldsymbol{c}$ are then plugged in the corresponding row of a large $N_I\times N_I$ sparse matrix $\mathbf{L}$, $N_I$ being the total number of internal nodes, leading to the following linear system:
\begin{equation}
    \mathbf{L} \boldsymbol{u}^h = \boldsymbol{q} + \boldsymbol{f}
    \label{eq:sparseFinal}
\end{equation}
where $\boldsymbol{q}$ comes from the boundary conditions, i.e., from those terms of $\boldsymbol{c}$ which are associated to boundary nodes in equation \eqref{eq:finalStencil}, while $\boldsymbol{u}^h$ and $\boldsymbol{f}$ are the solution vector and the vector of the known term $f$, respectively, at the internal nodes. Equation \eqref{eq:sparseFinal}, which is the discrete version of the PDE $\mathfrak{L}u = f$ we started with, can now be solved for $\boldsymbol{u}^h\in \mathbb{R}^{N_I}$.

We point out that, if $\mathbf{M}$ is ill conditioned in equation \eqref{eq:interpSystem}, so is $\mathbf{M}^T$ of equation \eqref{eq:adjoint}. We will show that the onset of instabilities at the boundary is due to ill conditioning of the interpolation matrix $\mathbf{M}$, and that such ill conditioning is essentially influenced by the direction of the boundary normals appearing in the corresponding rows, see equations \eqref{eq:phimat} and \eqref{eq:pmat}. 

The final sparse linear system of equation \eqref{eq:sparseFinal} can now be solved using a proper iterative method. Because of the lack of connectivity information, local interpolation systems $\eqref{eq:adjoint}$ can be solved independently from one another.

\subsection{Cardinal functions and Lebesgue constant} \label{ss:card_func}
The interpolant $u^h$ can also be obtained by considering the identity operator $\mathfrak{L}=\mathcal{I}$ in the left side of equation \eqref{eq:collocation}. Then, by considering the left side of equation \eqref{eq:finalStencil}, we obtain explicitly:
\begin{equation}
    \mathcal{I}u^h(\bx) = u^h(\bx) = \sum_{i=1}^{m_I} u(\bx_i) \psi_i(\bx) + \! \sum_{k=m_J}^m 
    \! \partial u(\bx_k) \psi_k(\bx)
    \label{eq:cardinalExpansion}
\end{equation}
where the cardinal functions $\{\psi_i(\bx)\}_1^m$ are therefore given by the first $m$ elements of $\boldsymbol{c}$ in equation \eqref{eq:adjoint} in the case $\mathfrak{L}=\mathcal{I}$, i.e., when the right side of the same equation is the column vector of the bare basis functions:
\begin{equation}
    \mathbf{M}^T
    \begin{Bmatrix}
        \{\psi_i(\bx)\}_1^m\\
        \boldsymbol{0}
    \end{Bmatrix}
    =
    \begin{Bmatrix}
        \{ \varphi_i(\bx) \}_1^m\\
        \{ p_j(\bx) \}_1^q\\
    \end{Bmatrix}
    \label{eq:adjoint_cardinal}
\end{equation}

These cardinal functions extend the definition given in \cite{de2010stability,FOR2007379} to the case with Neumann BCs. The cardinal functions associated to the internal nodes, i.e., $\{\psi_i(\bx)\}_1^{m_I}$, satisfy the following properties:
\begin{equation}
    \begin{cases}
        \begin{aligned}
                    \psi_i(\bx_j) &=\delta_{ij} && \text{if $j\le m_I$} \\
            \partial\psi_i(\bx_j) &=0           && \text{otherwise}
        \end{aligned}
    \end{cases}
\end{equation}
while the ones associated to the boundary nodes, i.e., $\{\psi_k(\bx)\}_{m_J}^m$, satisfy the following properties:
\begin{equation}
    \begin{cases}
        \begin{aligned}
                    \psi_k(\bx_j)&=0           && \text{if $j\le m_I$} \\
            \partial\psi_k(\bx_j)&=\delta_{kj} && \text{otherwise}
        \end{aligned}
    \end{cases}
\end{equation}
where $\delta_{ij}$ is the Kronecker delta.

In order to decouple the contribution of internal and boundary nodes explicitly, we define two Lebesgue functions:
\begin{equation}
\begin{split}
    \lambda_I(\bx) &= \sum_{i=1}^{m_I} |\psi_i(\bx)|\\
    \lambda_B(\bx) & = \sum_{k = m_J}^m |\psi_k(\bx)|
\end{split}
\label{eq:LebFunc}
\end{equation}
where $\lambda_I$ is associated to the internal nodes, whereas $\lambda_B$ is associated to the boundary nodes. This can be seen as an extension of the traditional definition of a single Lebesgue function \cite{de2010stability}.

When the interpolation stencil is restricted to lie in a compact set $K \subset \mathbb{R}^d$, then the Lebesgue constant is defined as the maximum of the Lebesgue function \cite{de2010stability,FOR2007379}.
In our case we will consider $K$ to be the minimal convex set containing all of the stencil nodes (including those lying on the boundary).
In order to keep the distinction between boundary nodes and inner nodes, we define 2 constants:
\begin{equation}
\begin{split}
    \Lambda_I & = \max_{\bx \in K} \lambda_I(\bx)\\
    \Lambda_B & = \max_{\bx \in K} \lambda_B(\bx)
    \label{eq:LebConst_IB}
\end{split}
\end{equation}

The Lebesgue constant is especially useful since it provides an upper bound for the interpolation error, which in turn is closely related to the conditioning of matrix $\mathbf{M}^T$. From equation \eqref{eq:cardinalExpansion} $u^h$ can be seen as the image of a function $u$ with respect to the linear interpolation operator $u \mapsto u^h$. We can find an upper bound for $||u^h||_{{\infty}}$ on $K$ as follows:
\begin{equation}
\begin{split}
||u^h||_{{\infty}} &= \max_{x \in K} \left| \sum_{i=1}^{m_I} u(\bx_i) \psi_i(\bx) + \sum_{k=m_J}^m \partial u(\bx_k) \psi_k(\bx) \right|\\
 & \leq \max_{\bx \in K} \left| \sum_{i=1}^{m_I} u(\bx_i) \psi_i(\bx) \right| + \max_{\bx \in K} \left| \sum_{k=m_J}^m \partial u(\bx_k) \psi_k(\bx) \right| \\
 & \leq \left( \max_{i \leq m_I} |u(\bx_i)| \right) \Lambda_I + 
 \left( \max_{m_J \leq k \leq m} | \partial u(\bx_k)| \right) \Lambda_B\\
 & \leq ||u||_{{\infty}} \Lambda_I + \underbrace{ \max_{\bx \in K} || \nabla u(\bx)||_2 }_{\displaystyle C_{\infty}(u)} \Lambda_B
\end{split}
\label{eq:maj1}
\end{equation}
where in the last inequality the Euclidean length of the gradient $\nabla u$ is employed as the upper bound for the absolute value of the normal derivative and $C_{\infty}(u)$ is the $L^\infty$-norm on $K$ of the Euclidean length of $\nabla u$. When the stencil does not include any boundary node, $\Lambda_B=0$ and equation \eqref{eq:maj1} takes the usual form \cite{de2010stability,FOR2007379}:
\begin{equation}
    ||u^h||_{{\infty}} \leq \Lambda_I ||u||_{{\infty}}
    \label{eq:internal_case}
\end{equation}
where $\Lambda_I$ can therefore be seen as an upper bound for the operator norm of the RBF interpolation operator.

In the general case, from equation \eqref{eq:maj1} we can derive an estimate for the interpolation error; assume that $u^*$ is the best RBF approximation possible of $u$, while $u^h$ is a generic approximation:
\begin{equation}
\begin{split}
    || u - u^h ||_{{\infty}} & \leq ||u - u^*||_{{\infty}} + ||u^* - u^h||_{{\infty}}\\
    & = ||u - u^*||_{{\infty}} + ||(u - u^*)^h||_{{\infty}}\\
    & \leq (1 + \Lambda_I) ||u - u^*||_{{\infty}} +  \Lambda_B\underbrace{ \max_{\bx \in K} || \nabla (u-u^*)(\bx)||_2 }_{\displaystyle C_{\infty}(u - u^*)}
\end{split}
\label{eq:error1}
\end{equation}
Once again, in the case without boundary nodes we have $\Lambda_B = 0$.

In the attempt of reaching a more compact form, similar to the one in \eqref{eq:internal_case}, we now assume that equation \eqref{eq:PhiCondition} holds. Using multi-index notation and the definition of the $||\cdot||_{W^{1,\infty}}$ norm given in \cite{brezis2011functional} we can write:
\begin{equation}
    \def\SA{\sum_{|\alpha| \leq 1}}
    \def\SI{\sum_{i=1}^{m_I}}
    \def\SB{\sum_{k=m_J}^m}
    \begin{split}
        ||u^h||_{\Wi} &=
\SA \max_{\bx \in K}
\left| \SI          u(\bx_i) \psi_i^{(\alpha)}(\bx) +
       \SB \partial u(\bx_k)\psi_k^{(\alpha)}(\bx) \right| \\
& \leq \SA \max_{\bx \in K}
\left( \left| \SI          u(\bx_i) \psi_i^{(\alpha)}(\bx) \right| +
    \left| \SB \partial u(\bx_k)\psi_k^{(\alpha)}(\bx) \right| \right)\\
& \leq
\SA \left( \max_{i \leq m_I}        |          u(\bx_i)| \right) \Lambda_I^{(\alpha)} +
\SA \left( \max_{m_J \leq k \leq m} | \partial u(\bx_k)| \right) \Lambda_B^{(\alpha)}\\
& = \left( \max_{i \leq m_I}        |          u(\bx_i)| \right) \SA \Lambda_I^{(\alpha)} +
    \left( \max_{m_J \leq k \leq m} | \partial u(\bx_k)| \right) \SA \Lambda_B^{(\alpha)}\\
& \leq ||u||_{\Wi} \underbrace{\SA (\Lambda_I^{(\alpha)} + \Lambda_B^{(\alpha)})}_{\displaystyle \Lambda_m} 
    \end{split}
    \label{eq:maj2}
\end{equation}
where $\Lambda_I^{(\alpha)}$ and $\Lambda_B^{(\alpha)}$ are now defined by equations \eqref{eq:LebFunc}-\eqref{eq:LebConst_IB} with $\psi_i^{(\alpha)}$ instead of $\psi_i$.
From equation \eqref{eq:maj2} we see that $\Lambda_m$ can be used as an upper bound for an operator norm of the RBF interpolation operator when it is defined on the Sobolev space $\Wi$ \cite{brezis2011functional}. 
If $u \in \Wi$ and equation \eqref{eq:PhiCondition} holds, we can derive another estimate for the interpolation error which follows from equation \eqref{eq:maj2}:
\begin{equation}
    \begin{split}
         || u - u^h ||_{\Wi} &\leq ||u - u^*||_{\Wi} + ||u^* - u^h||_{\Wi}\\
         & = ||u - u^*||_{\Wi} + ||(u - u^*)^h||_{\Wi}\\
         & \leq (1+\Lambda_m)||u - u^*||_{\Wi}
    \end{split}
    \label{eq:error2}
\end{equation}

In conclusion, we remark that any $\psi_i^{(\alpha)}(\bx)$ is obtained as the $i^{th}$ component of the solution of a linear system similar to the one in equation \eqref{eq:adjoint_cardinal}, where the the left hand side is kept the same and $\{ \varphi_i(\bx) \}_1^m$ is changed with $\{ \varphi_i^{(\alpha)}(\bx) \}_1^m$ and $\{ p_j(\bx) \}_1^q$ with $\{ p_j^{(\alpha)}(\bx) \}_1^q$. We now see that, whenever the matrix $\mathbf{M}^T$ is ill-conditioned, the terms $\Lambda_I^{(\alpha)}$ and $\Lambda_B^{(\alpha)}$ are unbounded.

\section{Neumann Boundary Conditions}\label{s:Neumann}
\subsection{Preliminary considerations}
\begin{figure}[t]
    \centering
    \includegraphics[width=.5\textwidth]{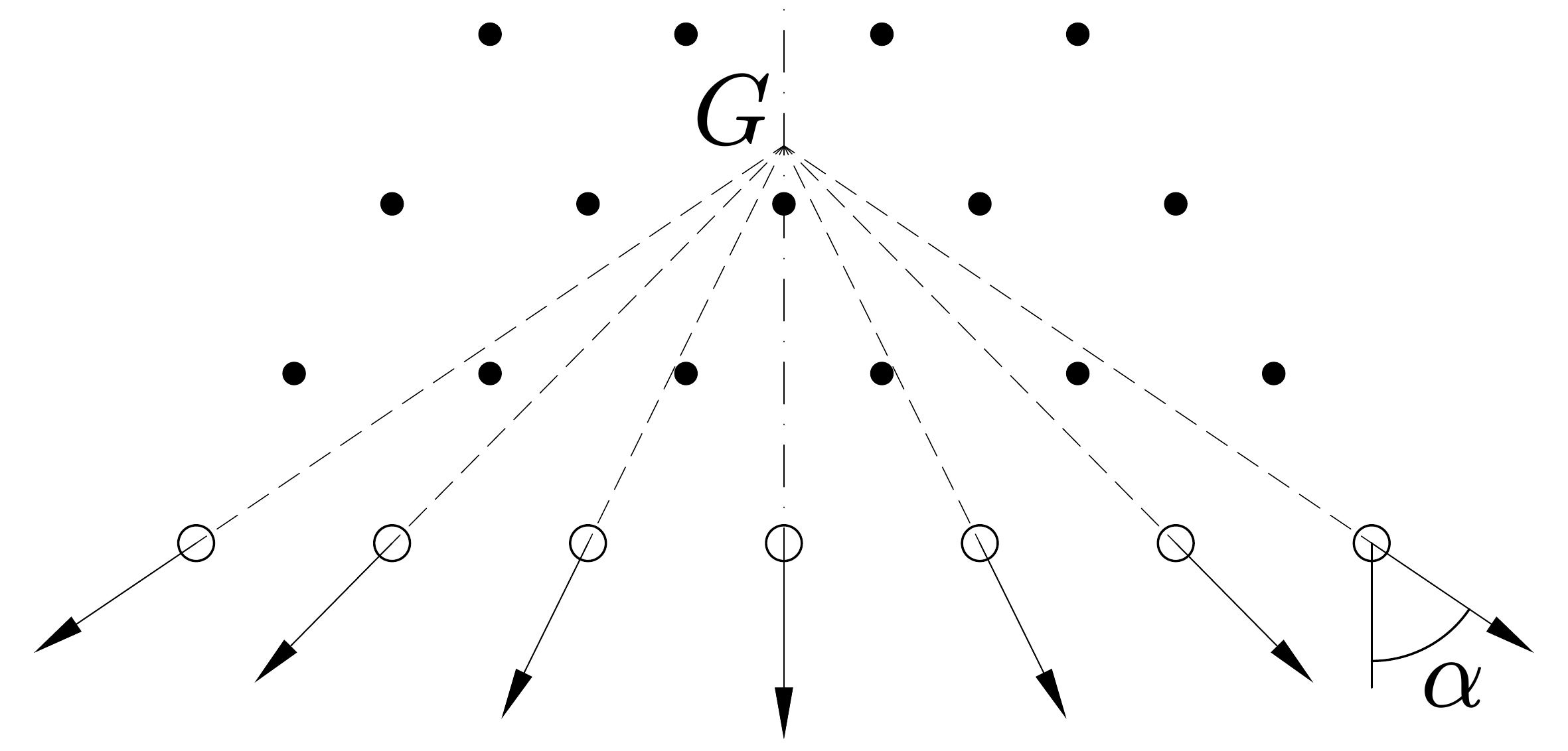}
    \caption{reference stencil with $m_I=15$ internal nodes (filled dots) and $m_B=7$ boundary nodes (empty circles).}
    \label{FIG:ref_stencil}
\end{figure}

\begin{figure}[t]
    \def\SpaceBelowText{.25em}
    \def\ImagesWidth{.49\textwidth}
    \centering
    \begin{subfigure}[b]{\ImagesWidth}
        \centering
        \includegraphics[width=\textwidth]{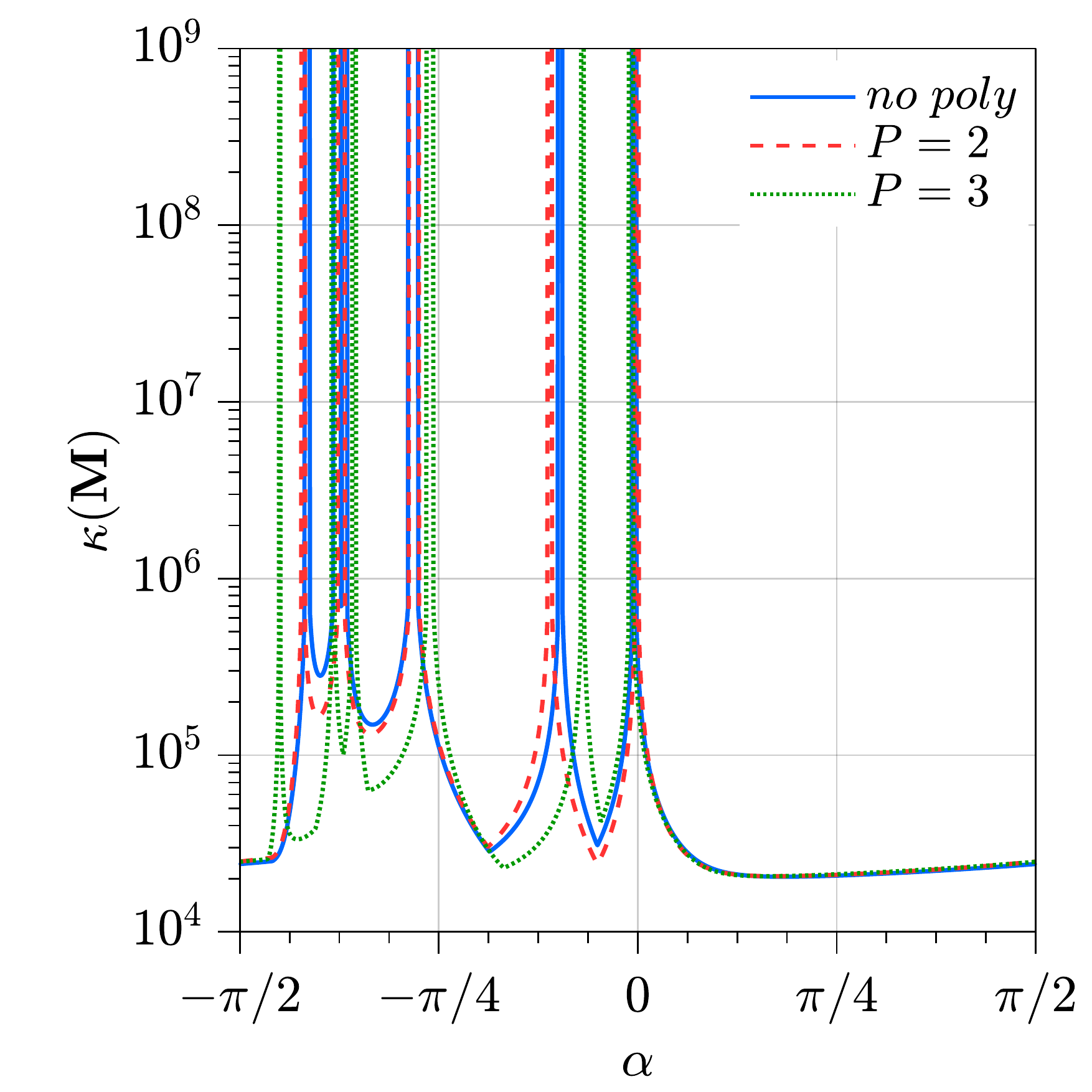}
    \end{subfigure}
    \hfill
    \begin{subfigure}[b]{\ImagesWidth}
        \centering
        \includegraphics[width=\textwidth]{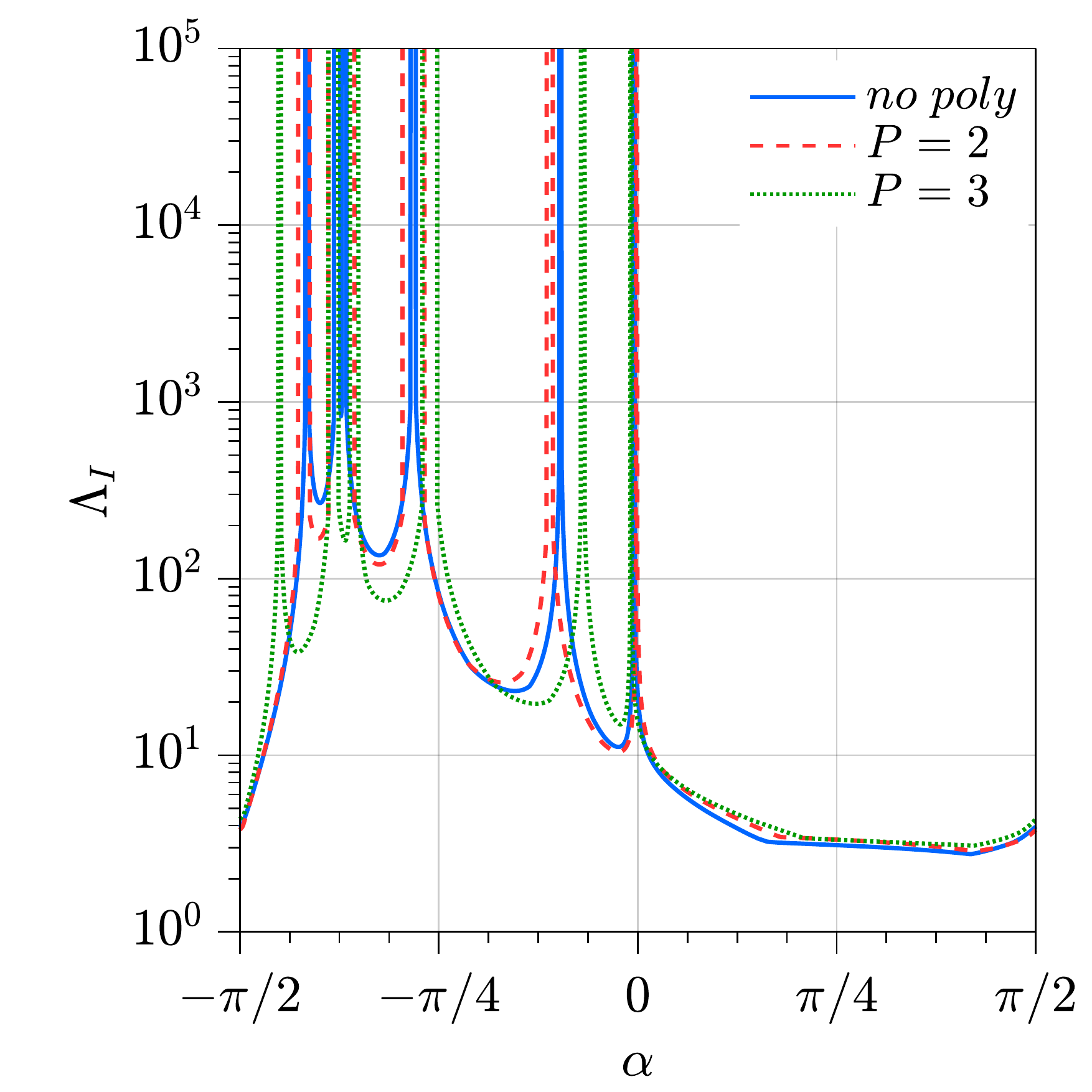}
    \end{subfigure}
    \caption{condition number of the interpolation matrix (left), and Lebesgue constant $\Lambda_I$ (right), for the RBF interpolation on the reference stencil of Figure \ref{FIG:ref_stencil}.}
    \label{FIG_singularities_alfa}
\end{figure}

\begin{figure}[t]
    \def\ImagesWidth{.49\textwidth}
    \def\SpaceBelowFirstRow{.015\textwidth}
    \centering
    \begin{subfigure}{\ImagesWidth}
        \centering
        \hspace{3em}\footnotesize{$\alpha = -\displaystyle\frac{\pi}{4}-0.02$}\\[1em]
        \includegraphics[width=\textwidth]{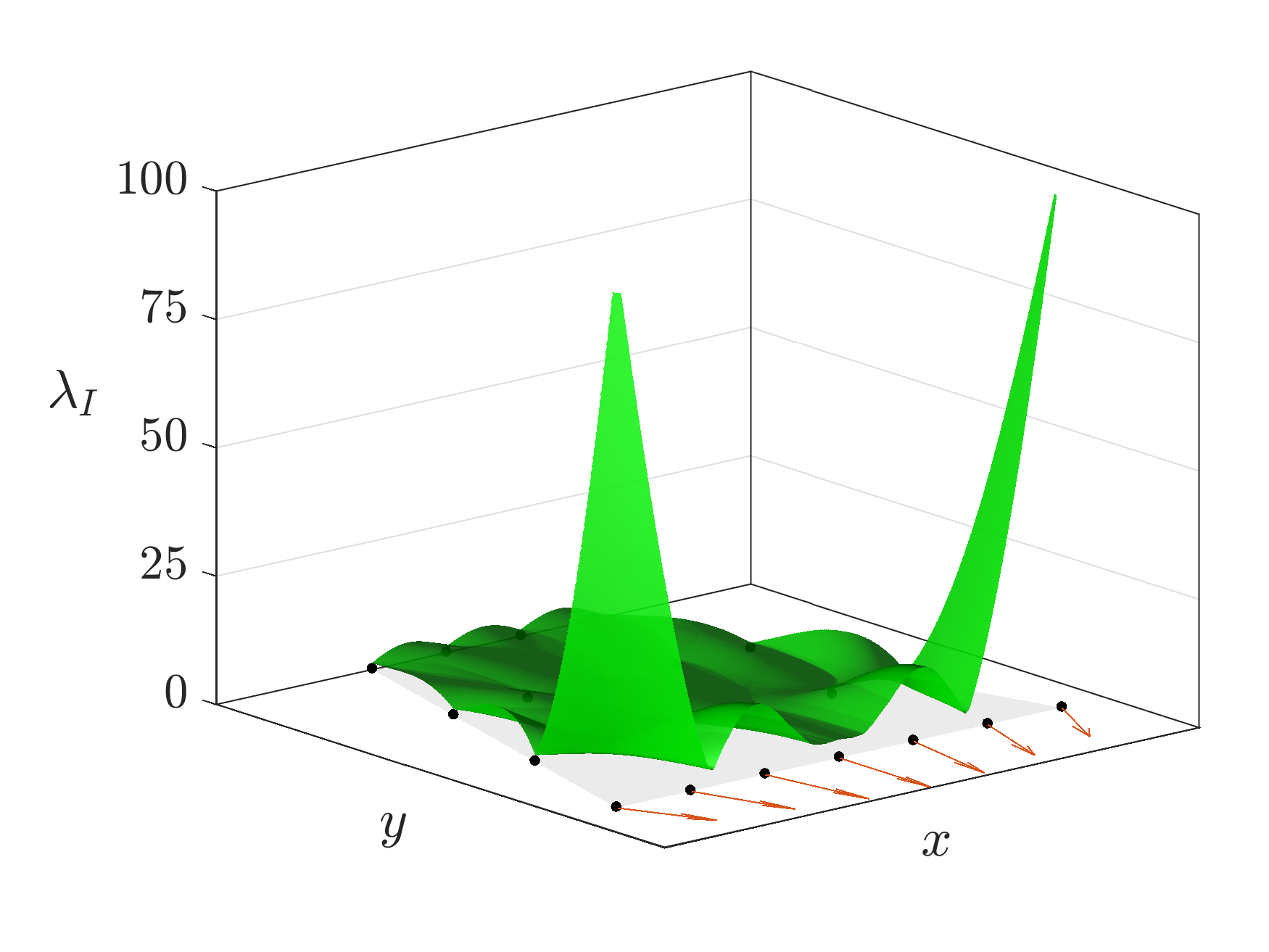}
    \end{subfigure}
    \hfill
    \begin{subfigure}{\ImagesWidth}
        \centering
        \hspace{3em}\footnotesize{$\alpha = \displaystyle\frac{\pi}{4}$}\\[1em]
        \includegraphics[width=\textwidth]{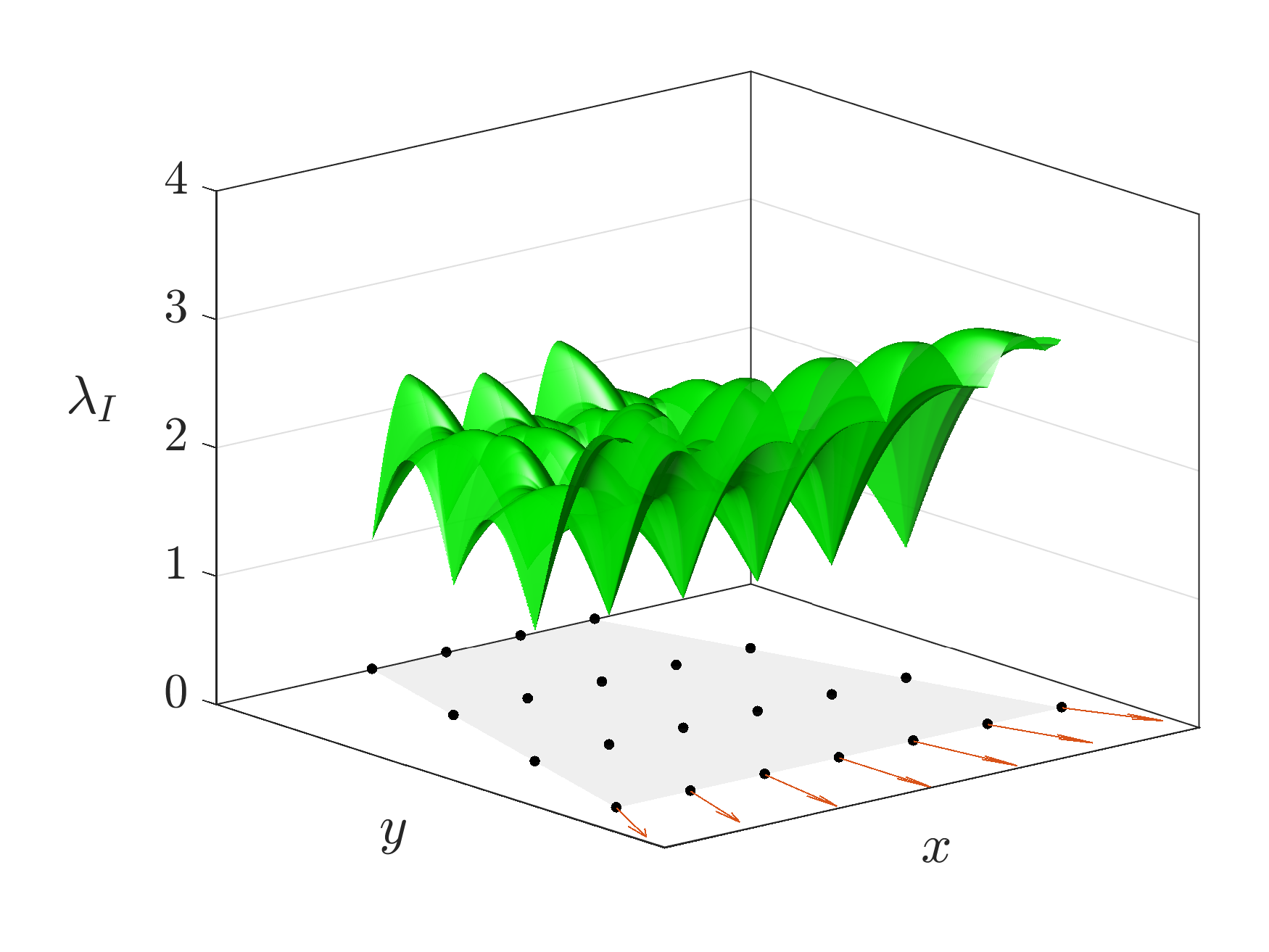}
    \end{subfigure}\\
     \begin{subfigure}{\ImagesWidth}
        \centering
        \includegraphics[width=\textwidth]{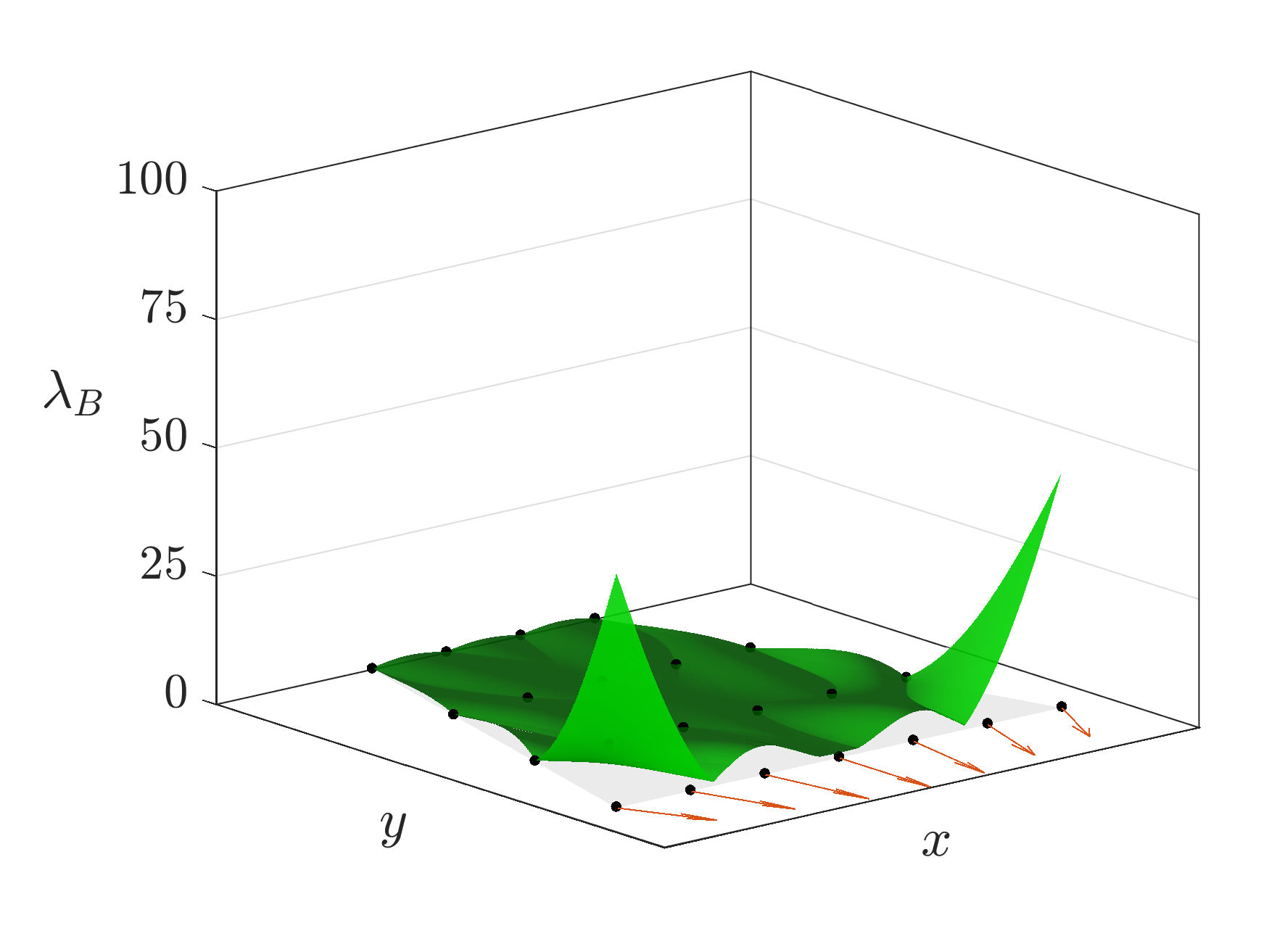}
    \end{subfigure}
    \hfill
    \begin{subfigure}{\ImagesWidth}
        \centering
        \includegraphics[width=\textwidth]{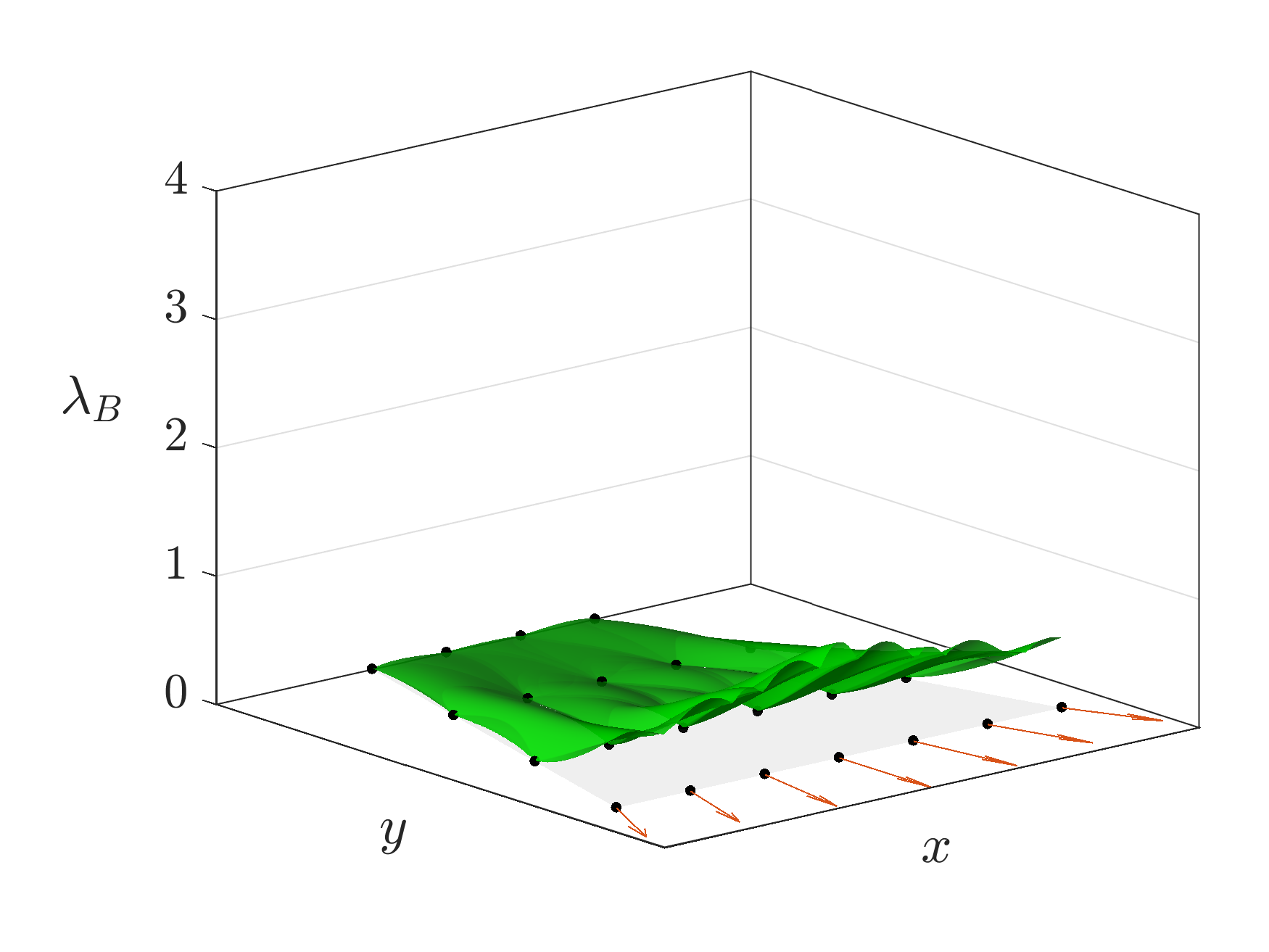}
    \end{subfigure}\\
    \caption{3D visualization of the Lebesgue functions $\lambda_I(\boldsymbol{x})$ and $\lambda_B(\boldsymbol{x})$ for the MQ RBF interpolation ($\varepsilon s=0.5$ and no polynomial augmentation) on the stencil of Figure \ref{FIG:ref_stencil}. Note the large difference in the $z$-axis scaling for the two values of $\alpha$.}
    \label{fig:cardinal_alpha}
\end{figure}

Figure \ref{FIG:ref_stencil} represents a stand-alone simplified model of the typical situation we encounter when the interpolation system of equation \eqref{eq:adjoint} is solved for a stencil close to the boundary. Nodes are arranged in a perfectly hexagonal layout and each unit normal is arranged along the line joining the corresponding boundary node with point $G$, and pointing outward.
The direction of the normals can thus be controlled by moving the point $G$ along the axis of symmetry, this way an angle $\alpha \in [ -\pi/2 \, , \, \pi/2 ]$ remains defined by the normal located further to the right. $\alpha>0$ if $G$ is above the line of the boundary nodes, $\alpha<0$ otherwise.

In Figure \ref{FIG_singularities_alfa} we report the condition number $\kappa(\mathbf{M})$ of the interpolation matrix and the constant $\Lambda_I$, defined in equation \eqref{eq:LebConst_IB}, against the angle $\alpha$ in the case of the reference stencil of Figure \ref{FIG:ref_stencil}. Whenever one of these quantities grows unbounded, we have a very ill-conditioned system in equation \eqref{eq:adjoint}, which translates into large errors and most likely stability issues. All plots of Figures \ref{FIG_singularities_alfa} and \ref{fig:cardinal_alpha} are attained using MQ RBF with shape parameter $\varepsilon$ satisfying $\varepsilon  s=0.5$ (see Table \ref{tab:RBFs}), where $s$ is the distance between the nodes.

It is possible to see that the negative values of $\alpha$ are the most dangerous, they are associated with locations of the point $G$ placed below the boundary in Figure \ref{FIG:ref_stencil}, i.e., with converging normals.
From Figure \ref{FIG_singularities_alfa} we see that multiple singularities appear for certain negative values of $\alpha$, near which both $\kappa(\mathbf{M})$ and $\Lambda_I$ assume very large values. Furthermore such singular configurations exist for modest values of $\alpha$ regardless of the degree of the polynomial augmentation and this remains true even with larger stencils and polynomials of higher degrees.
These remarks seem to tell us that no matter the node density or the polynomial degree, even a slight local concavity of the boundary might cause severe ill-conditioning issues for the considered RBF interpolations.

In Figure \ref{fig:cardinal_alpha} it is possible to qualitatively appreciate the functions $\lambda_I(\boldsymbol{x})$ and $\lambda_B(\boldsymbol{x})$ for two particular values of the angle $\alpha$, $\lambda_I$ and $\lambda_B$ being defined in equation \eqref{eq:LebFunc}. These plots are again attained using the same MQ RBF scheme with no polynomial augmentation, while the inclusion of the polynomial augmentation would have increased even more the values of $\lambda_I$ and $\lambda_B$ near the edges of the stencil. 

Once again, we can see that for certain negative values of $\alpha$, i.e., for critical directions of the normals, large interpolation errors can occur. We can therefore see that Neumann boundary conditions might compromise the accuracy and the well-posedness of the interpolation scheme of equation \eqref{eq:adjoint} and, through this, compromise the final solution of system \eqref{eq:sparseFinal}. In the following sections we aim at providing a more precise description for this phenomenon, along with some strategies for avoiding singular configurations of the normals.

\subsection{Bare RBF with one boundary node}
\subsubsection{Singular direction}\label{sss:singular_dir_1node}
Consider for simplicity a stencil with $m_B=1$ boundary node $\bx_m$, i.e., $m_I=m-1$, and an interpolant with no polynomial augmentation as in equation \eqref{eq:bareboneRBF}. The $m\times m$ interpolation matrix is therefore:
\begin{equation}
    \mathbf{M}=\boldsymbol{\varphi}_{BC}=
    \begin{bmatrix}
        \varphi_1(\bx_1)          \mysep \cdots \mysep \varphi_m(\bx_1)     \\
        \vdots                  \mysep \ddots \mysep \vdots             \\
        \varphi_1(\bx_{m_I})      \mysep \cdots \mysep \varphi_m(\bx_{m_I}) \\
        \partial\varphi_1(\bx_m) \mysep \cdots \mysep \partial\varphi_m(\bx_m)
    \end{bmatrix}
    \label{eq:InterpM_pureRBF_OneBoundaryNode}
\end{equation}

The determinant of $\mathbf{M}$ can be expressed through a cofactor expansion along the last row:
\begin{equation}
    \det(\mathbf{M})=\sum_{j=1}^{m_I}C_{m,j}\partial\varphi_j(\bx_m)
    \label{eq:OneBoundaryNode_DetCofactors}
\end{equation}

In equation \eqref{eq:OneBoundaryNode_DetCofactors}, $C_{m,j}=(-1)^{m+j}D_{m,j}$ is the cofactor of entry $M_{m,j}=\partial\varphi_j(\bx_m)$ in $\mathbf{M}$, where $D_{m,j}$ is the minor of the corresponding entry, i.e., the determinant of the submatrix obtained by deleting row $m$ and column $j$ from $\mathbf{M}$. We point out that the expansion in equation \eqref{eq:OneBoundaryNode_DetCofactors} has $m_I=m-1$ terms because the diagonal entry of the last row is $\partial\varphi_m(\bx_m)=0$, as outlined in section \ref{ss:RBF_FD}.

Recalling the form of the normal derivative and the notation of equation \eqref{eq:normal_phi}, the determinant can be expressed as:
\begin{equation}
    \det(\mathbf{M})=\sum_{j=1}^{m_I}C_{m,j}\Phi'(r_{j,m})\,\boldsymbol{e}_{j,m}\cdot\boldsymbol{n}=\boldsymbol{v}\cdot\boldsymbol{n}
    \label{eq:OneBoundaryNode_DetCofactors_v}
\end{equation}
where $\boldsymbol{n}=\boldsymbol{n}_m$ is the unit normal at the boundary node $\bx_m$. Thus $\boldsymbol{v}$ turns out to be a linear combination of $\boldsymbol{e}_{j,m}$, i.e., the unit directions pointing towards the boundary node $\bx_m$, as defined in section \ref{ss:RBF_FD}. Denoting with $w_{j,m}=C_{m,j}\Phi'(r_{j,m})$ the coefficient of such linear combination, it takes the form:
\begin{equation}
    \boldsymbol{v}=\sum_{j=1}^{m_I}w_{j,m}\boldsymbol{e}_{j,m}
    \label{eq:OneBoundaryNode_DetCofactors_optimal_v}
\end{equation}

When the location of the boundary node $\bx_m$ is fixed and the unit normal $\boldsymbol{n}$ can vary, the direction of $\boldsymbol{v}$ can be interpreted as the optimal direction for $\boldsymbol{n}$ which maximizes $\det(\mathbf{M})$ (or minimizes it in the opposite direction). On the other hand, equation \eqref{eq:OneBoundaryNode_DetCofactors_v} states that if $\boldsymbol{n}$ is orthogonal to $\boldsymbol{v}$, then $\mathbf{M}$ becomes singular.
\comment{
normals $\boldsymbol{n}$ orthogonal to $\boldsymbol{v}$ lead to a singular interpolation matrix $\mathbf{M}$.}

The existence of such singular directions can also be directly deduced from the following fact. If we denote the interpolation matrix in equation \eqref{eq:InterpM_pureRBF_OneBoundaryNode} by $\mathbf{M}(\boldsymbol{n})$, i.e., as a function of the unit normal $\boldsymbol{n}$, then  $\det(\mathbf{M}(\boldsymbol{n}))$ is a continuous function of $\boldsymbol{n}$ and we observe that $\det(\mathbf{M}(-\boldsymbol{n}))=-\det(\mathbf{M}(\boldsymbol{n}))$ for any $\boldsymbol{n}$, since the reversal of $\boldsymbol{n}$ results in the change of sign of the last row of $\mathbf{M}(\boldsymbol{n})$. Since $\boldsymbol{n}$ can be changed continuously to reach $-\boldsymbol{n}$, there exist at least one normal for which the determinant vanishes.

\subsubsection{Possible remedies}\label{sss:Possible remedies}
The aforementioned negative result can be addressed by using different strategies, as outlined in section \ref{s:Intro}. In this paper we will focus on the following two approaches:

{\flushleft
\begin{tabular}{l p{.75\textwidth}}
    Approach 1. & Discard a boundary node if the corresponding dot product $|\boldsymbol{v}\cdot\boldsymbol{n}|/\|\boldsymbol{v}\|_2$ is too small. \\
    Approach 2. & Change the boundary node location according to the corresponding normal.
\end{tabular}
}

Both approaches are inspired by the negative result on scattered data interpolation in multiple dimensions expressed by the Curtis-Mairhuber Theorem \cite{fasshauer2007meshfree}. In the case of interpolation with Neumann BCs this translates into the following: the basis functions should depend not only on the data locations, but also on the normals.
\comment{
This approach can be considered an extension of the principle induced by the negative result on scattered data interpolation in multiple dimensions expressed by the Curtis-Mairhuber Theorem \cite{fasshauer2007meshfree}. In the case of interpolation with conditions involving normal derivatives, the previous principle can be extended as follows: the basis functions should depend not only on the data locations, but also on the normals.}

\begin{figure}[t]
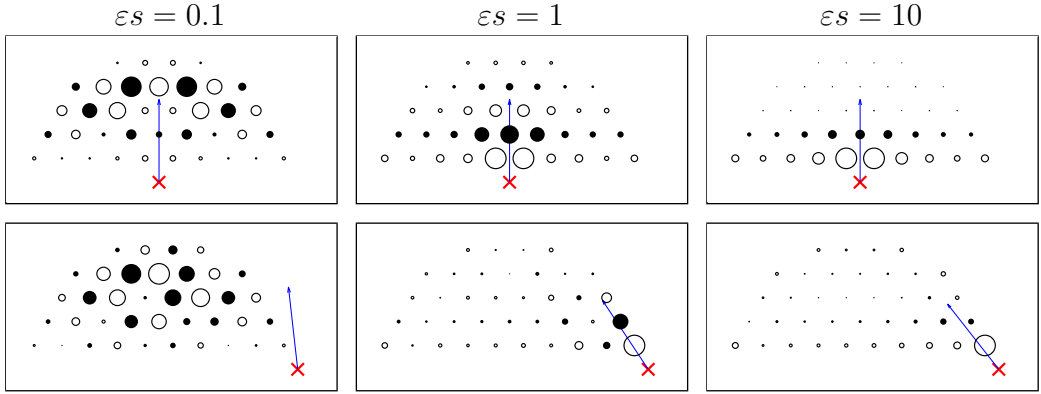

    \def\SpaceBelowText{.25em}
    \def\ImagesWidth{.32\textwidth}
    \def\SpaceBelowFirstRow{0.5em}
    \setlength\fboxsep{0pt}
    \setlength\fboxrule{0.5pt}
    \centering

    \begin{subfigure}[b]{\ImagesWidth}
        \centering
        $\varepsilon s=0.1$\\[\SpaceBelowText]
        \boxedgraphics{internal_weights_v/center_01.pdf}
    \end{subfigure}
    \hfill
    \begin{subfigure}[b]{\ImagesWidth}
        \centering
        $\varepsilon s=1$\\[\SpaceBelowText]
        \boxedgraphics{internal_weights_v/center_1.pdf}
    \end{subfigure}
    \hfill
    \begin{subfigure}[b]{\ImagesWidth}
        \centering
        $\varepsilon s=10$\\[\SpaceBelowText]
        \boxedgraphics{internal_weights_v/center_10.pdf}
    \end{subfigure}\\[\SpaceBelowFirstRow]
    \begin{subfigure}[b]{\ImagesWidth}
    \centering
        \boxedgraphics{internal_weights_v/side_01.pdf}
    \end{subfigure}
    \hfill
    \begin{subfigure}[b]{\ImagesWidth}
        \centering
        \boxedgraphics{internal_weights_v/side_1.pdf}
    \end{subfigure}
    \hfill
    \begin{subfigure}[b]{\ImagesWidth}
        \centering
        \boxedgraphics{internal_weights_v/side_10.pdf}
    \end{subfigure}
    \caption{coefficients $w_{j,m}$ in equation \eqref{eq:OneBoundaryNode_DetCofactors_v} for the optimal direction, shown as a blue arrow, associated to the boundary node, shown as a red cross. Positive and negative coefficients are shown as filled dots and empty circles, respectively, while the size of the markers is proportional to the magnitude of the coefficients. Top row: symmetric position of the boundary node; bottom row: asymmetric position of the boundary node.}
    \label{FIG_internal_weights_v}
\end{figure}

Approach 2 will be discussed in section \ref{ss:optimal_placement}.
Approach 1, on the other hand, requires either a calculation of $\boldsymbol{v}$ or a suitable approximation of it. Given the definition of $\boldsymbol{v}$ in equation \eqref{eq:OneBoundaryNode_DetCofactors_optimal_v}, the easiest approximation is to consider the direction pointing from the boundary node towards some reference point, e.g., the stencil centroid. The more constant the coefficients $w_{j,m}$ in equation \eqref{eq:OneBoundaryNode_DetCofactors_optimal_v}, the more accurate this approximation would be. Unfortunately, this condition never occurs, as shown in Figure \ref{FIG_internal_weights_v}. In this figure the magnitude and the sign of the coefficients $w_{j,m}$ are shown in the case of the MQ RBF for a boundary stencil, i.e., one-sided stencil close to a straight boundary, with hexagonal node arrangement and for two different positions of the boundary node, shown as a red cross. Three different values for the shape parameter $\varepsilon$ are considered, where the distances are scaled with respect to the constant nodal spacing $s$, i.e., $\varepsilon s=constant$ (the case $\varepsilon s=0.1$ has been computed with MATLAB variable precision arithmetic due to the extreme cancellations occuring for small $\varepsilon$ \cite{fornberg2004some}).

Figure \ref{FIG_internal_weights_v} shows that moving away from the boundary node, the coefficients $w_{j,m}$ have alternating signs and a decay rate growing with $\varepsilon$, which is in perfect accordance with theoretical observations for RBF interpolation on equispaced infinite lattices \cite{fornberg2008locality}. More specifically, by qualitative arguments based on the adaptation to a finite stencil, the cardinal expansion coefficients $\lambda_k$ defined in \cite{fornberg2008locality} correspond to the cofactors $C_{m,j}$ in equation \eqref{eq:OneBoundaryNode_DetCofactors} up to a multiplicative constant, where the subscript $k$ is the nondimensional distance from the boundary node.
\comment{
The qualitative behaviour of the magnitude of the coefficients $w_{j,m}=C_{m,j}\Phi'(r_{j,m})$ is therefore given by $|\lambda_k|\Phi'(ks)$, which in the case of MQ and with the exponential decay given by \cite{fornberg2008locality} becomes $\Tilde{w}(k)=ke^{-\mu k}/\sqrt{1+(\varepsilon sk)^2}$, after ruling out multiplicative constants.}

The qualitative behaviour of the magnitude of the coefficients $w_{j,m}=C_{m,j}\Phi'(r_{j,m})$ is therefore given by the following equation, where the second proportionality can be found in \cite{fornberg2008locality} for the exponential regime:
\begin{equation}
    |w_{j,m}| \propto |\lambda_k \Phi'(ks)| \propto e^{-\mu k}|\Phi'(ks)|
    \label{eq:qualitative}
\end{equation}

In the case of MQ RBF, equation \eqref{eq:qualitative} can be expressed in the form:
\begin{equation}
   |w_{j,m}| \propto \frac{ke^{-\mu k}}{\sqrt{1+(\varepsilon sk)^2}} = \Tilde{w}(k)
\end{equation}

When $\varepsilon s=0.1$, then $\mu=0.15$ \cite{fornberg2008locality} and $\Tilde{w}(k)$ has a maximum for $k\approx 5$, which in Figure \ref{FIG_internal_weights_v} corresponds to the large magnitude coefficients $w_{j,m}$ far from the boundary node. When $\varepsilon s=1$, then $\mu=1.04$ \cite{fornberg2008locality} and $\Tilde{w}(k)$ has a maximum for $k\approx 1$, which in Figure \ref{FIG_internal_weights_v} corresponds to the large magnitude coefficients $w_{j,m}$ immediately close to the boundary node, followed by a rapid decay. When $\varepsilon s=10$ the coefficients with largest magnitude are again associated to the nodes immediately close to the boundary node, followed by a stronger decay.

Figure \ref{FIG_internal_weights_v} also shows the optimal direction, i.e., the direction of $\boldsymbol{v}$ defined in equation \eqref{eq:OneBoundaryNode_DetCofactors_optimal_v}, as a blue vector. In the case of the symmetric position of the boundary node (top row in Figure \ref{FIG_internal_weights_v}), the optimal direction is vertical, regardless of $\varepsilon s$, as expected for symmetry. In the case of the asymmetric position of the boundary node (bottom row in Figure \ref{FIG_internal_weights_v}), the optimal direction changes from almost vertical for $\varepsilon s=0.1$ to a certain inclination for $\varepsilon s=10$, in accordance to the previous observations. Indeed, the larger the shape parameter $\varepsilon$, the larger the coefficients $w_{j,m}$ associated to the closest nodes and therefore the more the optimal direction points towards these closest nodes.

\begin{figure}[t]
    \def\SpaceBelowText{.25em}
    \centering
    \begin{subfigure}[b]{.31\textwidth}
        \centering
        \includegraphics[width=\textwidth]{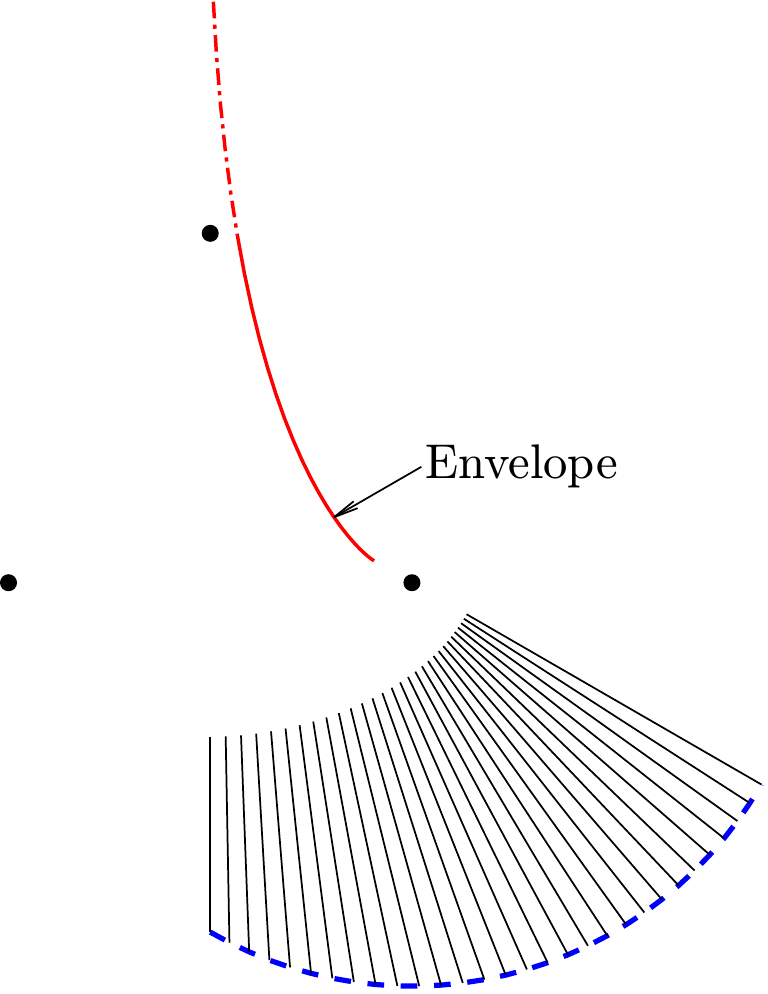}
    \end{subfigure}
    \hspace{2em}
    \begin{subfigure}[t]{.55\textwidth}
        \centering
        \includegraphics[width=\textwidth]{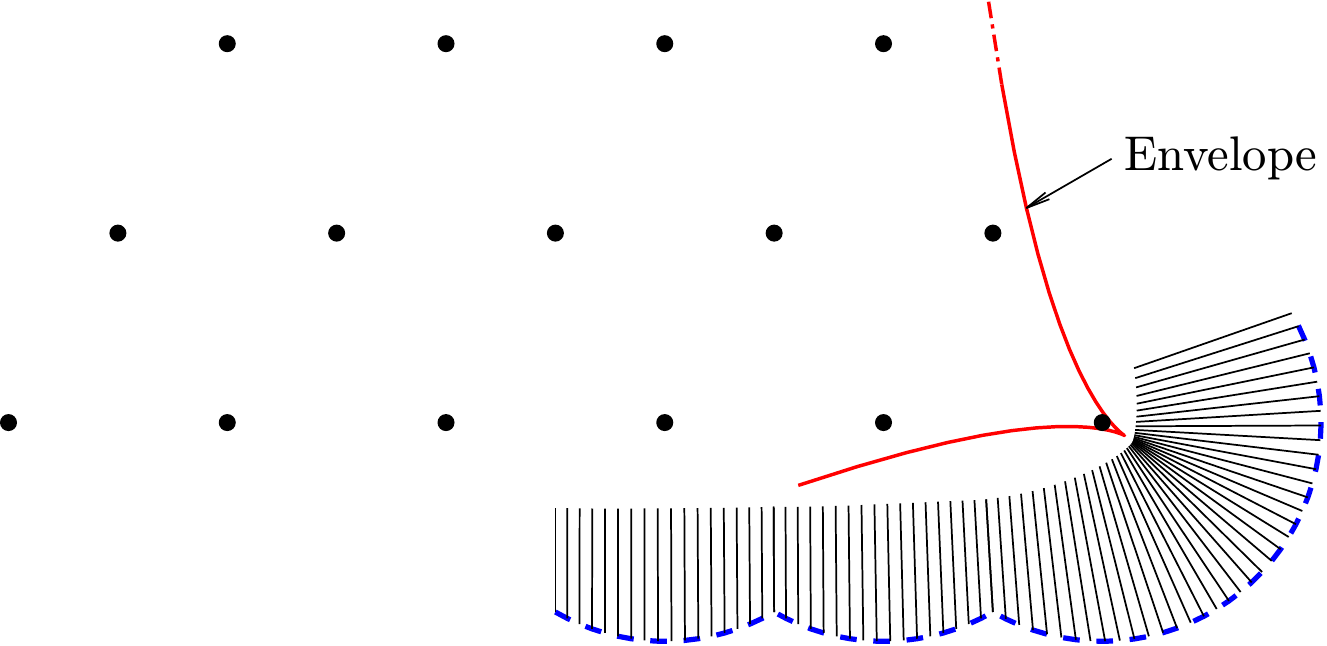}
    \end{subfigure}
    \caption{optimal vectors $\boldsymbol{v}$ (solid black lines) and their envelopes (solid red curves) for two stencils with $m_I=3$ (left) and $m_I=15$ (right) internal nodes as the boundary node moves along the blue dashed curve. This last curve represents the points having constant distance $s$, i.e., the nodal spacing, from the nearest internal node.}
    \label{FIG_envelopes}
\end{figure}

Figure \ref{FIG_envelopes} shows the direction and the magnitude of the optimal vector $\boldsymbol{v}$, depicted with solid black lines, for two hexagonal stencils, with $m_I=3$ (left) and $m_I=15$ (right) internal nodes, as the boundary node moves along the blue dashed curve. This blue dashed curve is formed by those points that have constant distance $s$ from the nearest internal node.\comment{and acts as a dummy boundary} MQ RBF with $\varepsilon s=0.5$ has been employed. The envelopes of the optimal directions of $\boldsymbol{v}$ are shown as solid red curves. In the case $m_I=3$ the optimal vectors seem to point qualitatively toward some reference mid point related to the internal nodes, although their envelope shows there is no such single point, or region, which can be heuristically assumed as a good reference point: the envelope indeed falls far out of the region delimited by the three internal nodes. We also note that the magnitude of $\boldsymbol{v}$, i.e., the length of the solid black lines, is slightly larger when the boundary node is located further to the right of the blue dashed curve. In such cases the optimal direction is pointing toward the closest internal node on the right and the envelope for these locations gets indeed closer to this node. This means that, given the role of $\boldsymbol{v}$ in equation \eqref{eq:OneBoundaryNode_DetCofactors_v}, these are the locations for the boundary nodes on the blue dashed curve that maximize $\det(\mathbf{M})$ when $\boldsymbol{n}=\boldsymbol{v}$.

In the case $m_I=3$, although a first look might suggests that all optimal directions converge at the same point, the fact that the envelope falls far out of the region delimited by the three internal nodes refutes this initial impression. Therefore, it is not possible to naively approximate the optimal directions as vectors radiating from a precise region.

The length of the solid black lines in Figure \ref{FIG_envelopes} are scaled according to the value of  $||\boldsymbol{v}||_2$. We see that whenever the boundary node is placed at the right margin, then the optimal direction, as well as $||\boldsymbol{v}||_2$, suggest that a better conditioning is achieved when the normal is parallel to the vector joining the boundary node with its closest internal neighbor.

The case $m_I=15$ is very similar to the previous one, with an envelope falling further away from the region delimited by the internal nodes, thus confirming that there is no such reasonable mid point for a simple geometric approximation of the optimal directions. This holds even in the case of very regular stencils with hexagonal arrangement like those employed here. The length of the solid black lines confirms that also in this case it is convenient to have a boundary normal pointing toward the nearest internal node at the right margin. By heuristic arguments, the enforcement of a Neumann BC at a boundary node should be ``supported" by a close internal node along the direction of the normal.

We conclude that an estimate for the optimal direction $\boldsymbol{v}$ can not be obtained through simple geometric heuristics and a more robust calculation is therefore mandatory in the general case. Such calculation follows here below.

\subsubsection{Practical calculation of \texorpdfstring{$\boldsymbol{v}$}{v}}\label{sss:practical_v}
The formulation given in equation \eqref{eq:OneBoundaryNode_DetCofactors_optimal_v} for $\boldsymbol{v}$ is helpful to highlight the previous observations but is not of practical interest. A more convenient approach can be derived by considering the following splitting of the interpolation matrix $\mathbf{M}$ given in equation \eqref{eq:InterpM_pureRBF_OneBoundaryNode}:
\begin{equation}
    \mathbf{M}=
    \begin{bmatrix}
    \begin{array}{ccc|c}
\varphi_1(\bx_1)         \mysep \cdots \mysep \varphi_{m_I}(\bx_1)     
& \varphi_m(\bx_1)            \\
\vdots                              \mysep \ddots \mysep \vdots                              
& \vdots                                 \\
\varphi_1(\bx_{m_I})     \mysep \cdots \mysep \varphi_{m_I}(\bx_{m_I}) 
& \varphi_m(\bx_{m_I})        \\ \hline
\partial\varphi_1(\bx_m) \mysep \cdots \mysep \partial\varphi_{m_I}(\bx_m)     
& 0
    \end{array}
    \end{bmatrix}
    =
    \begin{bmatrix}
    \begin{array}{c|c}
\boldsymbol{\varphi}_{II} & \boldsymbol{\varphi}_{IB} \\ \hline
\partial\boldsymbol{\varphi}_{BI} & 0
    \end{array}
    \end{bmatrix}
    \label{eq:InterpM_pureRBF_OneBoundaryNode_splitting}
\end{equation}
where the bottom-right diagonal entry is $\partial\varphi_m(\bx_m)=0$. $\boldsymbol{\varphi}_{II}$ is the symmetric $m_I\times m_I$ RBF interpolation matrix obtained when considering the $m_I$ internal nodes without the boundary node $\bx_m$, $\boldsymbol{\varphi}_{IB}$ is the column vector of the values of function $\varphi_m$, i.e., the RBF centered at $\bx_m$, evaluated at the $m_I$ internal nodes, while $\partial\boldsymbol{\varphi}_{BI}$ is the row vector of the normal derivatives of the RBFs centered at the $m_I$ internal nodes and evaluated at $\bx_m$.

By taking the Schur complement $S_{BB}$ of the block $\boldsymbol{\varphi}_{II}$ in equation \eqref{eq:InterpM_pureRBF_OneBoundaryNode_splitting}, we obtain the following scalar Schur complement:
\begin{equation}
    S_{BB} = -\partial\boldsymbol{\varphi}_{BI}\boldsymbol{\varphi}_{II}^{-1}\boldsymbol{\varphi}_{IB}=-\boldsymbol{\varphi}_{IB}^T\boldsymbol{\varphi}_{II}^{-1}\partial\boldsymbol{\varphi}_{BI}^T
    \label{eq:Schur_1node}
\end{equation}
where the second equality derives from a simple transposition since $\boldsymbol{\varphi}_{II}$ is symmetric. The following equality also holds:
\begin{equation}
    \det(\mathbf{M})=\det(\boldsymbol{\varphi}_{II})S_{BB}
    \label{eq:Schur_1node_det}
\end{equation}
where $\det(\boldsymbol{\varphi}_{II})\ne 0$ when using strictly positive definite or conditionally positive definite RBFs $\Phi(r)$ of order 1.

Because of the symmetry of the RBFs, i.e., $\varphi_m(\bx_j)=\varphi_j(\bx_m)$, the column vector $\boldsymbol{\varphi}_{IB}$ can also be interpreted as the values of the RBFs centered at the $m_I$ internal nodes, evaluated at the boundary node $\bx_m$. Therefore,  the column vector $\boldsymbol{\varphi}_{II}^{-1}\boldsymbol{\varphi}_{IB}$ in the first matrix expression in equation \eqref{eq:Schur_1node} can be interpreted as the values of the cardinal functions $\bar{\psi}_j$ of the RBF interpolation based on the internal nodes only, evaluated at $\bx_m$:
\begin{equation}
    \bar{\boldsymbol{\psi}} = \boldsymbol{\varphi}_{II}^{-1}\boldsymbol{\varphi}_{IB}
    \label{eq:OneBoundaryNode_psi_at_boundary}
\end{equation}
where $\bar{\boldsymbol{\psi}}=\{\bar{\psi}_j(\bx_m)\}_{j=1}^{m_I}$.

In the end, by comparing equation \eqref{eq:OneBoundaryNode_DetCofactors} with equation \eqref{eq:Schur_1node_det} and by using the first matrix expression in equation \eqref{eq:Schur_1node}, we can see that, up to the multiplicative constant $-\det(\boldsymbol{\varphi}_{II})$, the cofactors $C_{m,j}$ in equation \eqref{eq:OneBoundaryNode_DetCofactors} correspond to the aforementioned cardinal functions $\bar{\psi}_j$:
\begin{equation}
    C_{m,j}=-\det(\boldsymbol{\varphi}_{II})\bar{\psi}_j(\bx_m) \, , \quad j = 1, \dots, m_I
\end{equation}
Therefore, $\boldsymbol{v}$ can be computed from equation \eqref{eq:OneBoundaryNode_DetCofactors_optimal_v} recalling that  $w_{j,m}=C_{m,j}\Phi'(r_{j,m})$:
\begin{equation}
    \boldsymbol{v}=\sum_{j=1}^{m_I}C_{m,j}\Phi'(r_{j,m})\boldsymbol{e}_{j,m}
    = -\det(\boldsymbol{\varphi}_{II})
    \sum_{j=1}^{m_I}\bar{\psi}_j(\bx_m)\Phi'(r_{j,m})\boldsymbol{e}_{j,m}
    \label{eq:OneBoundaryNode_Schur_optimal_v}
\end{equation}
where $\det(\boldsymbol{\varphi}_{II})$ does not need to be evaluated since it is a multiplicative constant and the values $\bar{\psi}_j(\bx_m)$ can be computed from equation \eqref{eq:OneBoundaryNode_psi_at_boundary}, intended as the solution of the linear system $\boldsymbol{\varphi}_{II}\bar{\boldsymbol{\psi}} = \boldsymbol{\varphi}_{IB}$.

\subsubsection{Qualitative interpretation of \texorpdfstring{$S_{BB}$}{S-BB}}
The second matrix expression in equation \eqref{eq:Schur_1node}, i.e., $S_{BB} =-\boldsymbol{\varphi}_{IB}^T\boldsymbol{\varphi}_{II}^{-1}\partial\boldsymbol{\varphi}_{BI}^T$, provides the following explicit interpretation of $S_{BB}$. 

The entries of $\partial\boldsymbol{\varphi}_{BI}$ are antisymmetric, i.e., $\partial \varphi_j (\bx_m)=-\partial \varphi_m (\bx_j)$, from equation \eqref{eq:normal_phi}, since we always refer to the same normal $\boldsymbol{n}$. Therefore, the vector $\partial\boldsymbol{\varphi}_{BI}$ can be interpreted as the values of the normal derivative of the function $-\varphi_m$ 
evaluated at the $m_I$ internal nodes. The matrix-vector operation $\boldsymbol{\varphi}_{II}^{-1}\partial\boldsymbol{\varphi}_{BI}^T$ hence gives the column vector of the RBF expansion coefficients $\{\alpha_i\}_1^{m_I}$ of the RBF interpolation, based on the $m_I$ internal nodes only, of the function $-\partial\varphi_m$. 
The product $-\boldsymbol{\varphi}_{IB}^T (\boldsymbol{\varphi}_{II}^{-1}\partial\boldsymbol{\varphi}_{BI}^T)$ then gives the extrapolation of the function $\partial\varphi_m$ at $\bx_m$. 

Since $\partial\varphi_m(\bx_m)=0$, $S_{BB}$ can now be interpreted as the difference between the value of the function $\partial\varphi_m$ at $\bx_m$, i.e., 0, and the extrapolation of the same function at that same location $\bx_m$ when using a RBF scheme based on the internal nodes only. 
This interpretation, which will continue to hold in the case of multiple boundary nodes presented in section \ref{sss:RBF_multiple_bnd_nodes_S}, tells us that the interpolation matrix $\mathbf{M}$ in equation \eqref{eq:InterpM_pureRBF_OneBoundaryNode} will be ill-conditioned if the internal nodes are sufficient for a good RBF extrapolation of the function $\partial\varphi_m$ at the location of the boundary node. Indeed, if the RBF interpolation based on the internal nodes is not capable of that, the additional information provided by the Neumann BC, i.e., the last row of $\mathbf{M}$, is not redundant, thus leading to a well-conditioned interpolation matrix.

\subsubsection{Additional remarks}
\begin{figure}[t]
    \def\SpaceBelowText{.25em}
    \def\ImagesWidth{.32\textwidth}
    \def\SpaceBelowFirstRow{.015\textwidth}
    \setlength\fboxsep{0pt}
    \setlength\fboxrule{0.5pt}
    \centering
    \begin{minipage}[b]{.85\textwidth}
        \begin{subfigure}[b]{\ImagesWidth}
            \centering
            $m_I=3$\\[\SpaceBelowText]
            \boxedgraphics{single_node_det_zeros/3_color.pdf}
        \end{subfigure}
        \hfill
        \begin{subfigure}[b]{\ImagesWidth}
            \centering
            $m_I=5$\\[\SpaceBelowText]
            \boxedgraphics{single_node_det_zeros/5_color.pdf}
        \end{subfigure}
        \hfill
        \begin{subfigure}[b]{\ImagesWidth}
            \centering
            $m_I=12$\\[\SpaceBelowText]
            \boxedgraphics{single_node_det_zeros/12_color.pdf}
        \end{subfigure}\\[\SpaceBelowFirstRow]
        \begin{subfigure}[b]{\ImagesWidth}
            \centering
            \boxedgraphics{single_node_det_zeros/3_color_rand.pdf}
        \end{subfigure}
        \hfill
        \begin{subfigure}[b]{\ImagesWidth}
            \centering
            \boxedgraphics{single_node_det_zeros/5_color_rand.pdf}
        \end{subfigure}
        \hfill
        \begin{subfigure}[b]{\ImagesWidth}
            \centering
            \boxedgraphics{single_node_det_zeros/12_color_rand.pdf}
        \end{subfigure}
    \end{minipage}
    \hfill
    \begin{minipage}[b][0.38\textwidth][c]{.12\textwidth}
        \includegraphics[width=\textwidth]{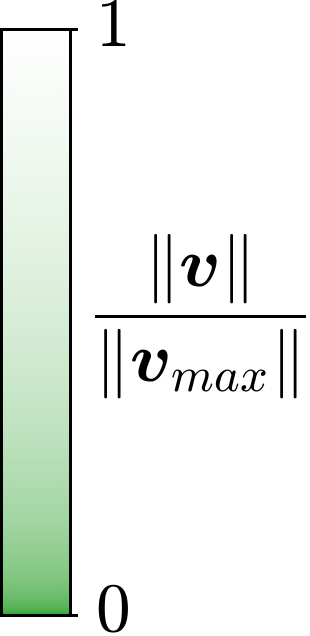}
    \end{minipage}
    \caption{norm of $\boldsymbol{v}$ for different node arrangements as the position of the boundary node changes. The $m_I$ internal nodes are shown as black empty circles and the red dots represent the positions of the boundary node for which $\boldsymbol{v}=\boldsymbol{0}$. Red dots represent therefore forbidden positions for the boundary node. Top row: hexagonal arrangements; bottom row: hexagonal arrangements perturbed by random displacements. $\|\boldsymbol{v}_{max}\|$ represents the maximum norm of $\boldsymbol{v}$ for each subfigure.}
    \label{FIG_single_node_det_zeros}
\end{figure}

A final remark before we conclude the discussion about the case with a single boundary node concerns the following questions. Given a node arrangement for the $m_I$ internal nodes, are there particular locations of the boundary node for which $\det(\mathbf{M})=0$ regardless of the direction of the normal $\boldsymbol{n}$, i.e., for which $\boldsymbol{v}=\boldsymbol{0}$ in equation \eqref{eq:OneBoundaryNode_DetCofactors_v} or, equivalently, $S_{BB}=0$ in equation \eqref{eq:Schur_1node_det}, provided that $\det(\boldsymbol{\varphi}_{II})\ne 0$? 
If so, is it possible that a Neumann boundary condition ends up being defined in one of such singular locations?

Fortunately, this seems not to be the case: Figure \ref{FIG_single_node_det_zeros} shows these singular locations as red dots for different node arrangements in the case of MQ RBF with $\varepsilon s=0.5$, together with a representation of $\|\boldsymbol{v}\|$. From this figure it can be observed that all singular locations fall ``inside" the stencil and dangerous locations, i.e., locations of the boundary node for which $\|\boldsymbol{v}\|$ is close to 0, are also in the neighbourhood of the stencil. In other words, the further the boundary node is from the internal nodes, the further $\mathbf{M}$ is from being singular, if the employed normal is not orthogonal to the optimal direction. Whenever we have a boundary node ``outside" the region delimited by the internal nodes, we do not need therefore to worry about the singular locations.
This principle has been verified numerically for each type of RBF presented in Table \ref{tab:RBFs} which is at least conditionally positive definite of order 1, i.e., RBFs that can be used for RBF interpolation without the need of any polynomial augmentation. For instance, the use of $\Phi(r)=r^2\log r$, which is conditionally positive definite of order 2 and therefore would require an augmentation with a linear polynomial, leads to the existence of singular locations far away from internal nodes.

\subsection{Bare RBF with multiple boundary nodes} 
\subsubsection{Preliminary results}\label{sss:RBF_multiple_bnd_nodes_S}
Consider now a stencil with $m_B$ boundary nodes and, again, an interpolant with no polynomial augmentation. The $m\times m$ interpolation matrix is therefore:
\begin{equation}
    \mathbf{M}=\boldsymbol{\varphi}_{BC}=
    \begin{bmatrix}
        \varphi_1(\bx_1)         \mysep \cdots \mysep \varphi_m(\bx_1)        \\
        \vdots                              \mysep \ddots \mysep \vdots                             \\
        \varphi_1(\bx_{m_I})     \mysep \cdots \mysep \varphi_m(\bx_{m_I})    \\
        \partial\varphi_1(\bx_{m_J}) \mysep \cdots \mysep \partial\varphi_m(\bx_{m_J})\\
        \vdots                              \mysep \ddots \mysep \vdots                             \\
        \partial\varphi_1(\bx_m) \mysep \cdots \mysep \partial\varphi_m(\bx_m)
    \end{bmatrix}
    \label{eq:InterpM_pureRBF_MultipleBoundaryNode}
\end{equation}
where $m_J=m_I+1$. Similarly to what was done in section \ref{sss:practical_v}, consider the following splitting of $\mathbf{M}$:
\begin{equation}
\begin{split}
    \mathbf{M}=&
    \begin{bmatrix}
    \begin{array}{ccc|ccc}
\varphi_1(\bx_1)           \mysep \cdots \mysep \varphi_{m_I}(\bx_1) &
\varphi_{m_J}(\bx_1)     \mysep \cdots \mysep \varphi_{m}(\bx_1)   \\
\vdots                                \mysep \ddots \mysep \vdots &
\vdots                                \mysep \ddots \mysep \vdots \\
\varphi_1(\bx_{m_I})       \mysep \cdots \mysep \varphi_{m_I}(\bx_{m_I}) &
\varphi_{m_J}(\bx_{m_I}) \mysep \cdots \mysep \varphi_{m}(\bx_{m_I})   \\ \hline
\partial\varphi_1(\bx_{m_J})       \mysep \cdots \mysep \partial\varphi_{m_I}(\bx_{m_J}) &
\partial\varphi_{m_J}(\bx_{m_J}) \mysep \cdots \mysep \partial\varphi_{m}(\bx_{m_J})   \\
\vdots                                          \mysep \ddots \mysep \vdots &
\vdots                                          \mysep \ddots \mysep \vdots \\
\partial\varphi_1(\bx_{m})           \mysep \cdots \mysep \partial\varphi_{m_I}(\bx_{m}) &
\partial\varphi_{m_J}(\bx_{m})     \mysep \cdots \mysep \partial\varphi_{m}(\bx_{m})
    \end{array}
    \end{bmatrix} \\
    =& 
    \begin{bmatrix}
    \begin{array}{c|c}
\boldsymbol{\varphi}_{II} & \boldsymbol{\varphi}_{IB} \\ \hline
\partial\boldsymbol{\varphi}_{BI} & \partial\boldsymbol{\varphi}_{BB}
    \end{array}
    \end{bmatrix}
\end{split}
\label{eq:InterpM_pureRBF_MultipleBoundaryNode_splitting}
\end{equation}
where the $m_B$ diagonal entries of $\partial\boldsymbol{\varphi}_{BB}$ are $\partial\varphi_i(\bx_i)=0$, $i=m_J,\ldots,m$, as previously stated. $\boldsymbol{\varphi}_{II}$ is again the symmetric $m_I\times m_I$ RBF interpolation matrix obtained when considering the $m_I$ internal nodes only. $\boldsymbol{\varphi}_{IB}$ is the $m_I\times m_B$ matrix whose columns are the RBFs centered at the $m_B$ boundary nodes and evaluated at the $m_I$ internal nodes, while $\partial\boldsymbol{\varphi}_{BI}$ is the $m_B\times m_I$ matrix whose rows are the normal derivatives of the RBFs centered at the $m_I$ internal nodes and evaluated at the $m_B$ boundary nodes.

By taking the Schur complement $\mathbf{S}_{BB}$ of the $\boldsymbol{\varphi}_{II}$ block of $\mathbf{M}$ in equation \eqref{eq:InterpM_pureRBF_MultipleBoundaryNode_splitting}, we obtain the following $m_B\times m_B$ Schur complement matrix:
\begin{equation}
    \mathbf{S}_{BB} = \partial\boldsymbol{\varphi}_{BB}-\partial\boldsymbol{\varphi}_{BI}\boldsymbol{\varphi}_{II}^{-1}\boldsymbol{\varphi}_{IB}
    \label{eq:Schur_Multiple_node}
\end{equation}

\comment{whose transpose is:
\begin{equation}
    \mathbf{S}_{BB}^T = \partial\boldsymbol{\varphi}_{BB}^T-\boldsymbol{\varphi}_{IB}^T\boldsymbol{\varphi}_{II}^{-1}\partial\boldsymbol{\varphi}_{BI}^T
    \label{eq:Schur_Multiple_node_T}
\end{equation}
since $\boldsymbol{\varphi}_{II}$ is symmetric.}

The following equality also holds:
\begin{equation}
    \det(\mathbf{M})=
    \det(\boldsymbol{\varphi}_{II})\det(\mathbf{S}_{BB})
    \label{eq:Schur_Multiple_node_det}
\end{equation}
where, again, $\det(\boldsymbol{\varphi}_{II})\ne 0$ when using strictly positive definite or conditionally positive definite functions $\Phi(r)$ or order 1.

\subsubsection{Qualitative interpretation of \texorpdfstring{$\mathbf{S}_{BB}$}{SS-BB}}\label{sss:qualitative_multiple}
Equation \eqref{eq:Schur_Multiple_node} provides again an explicit interpretation of $\mathbf{S}_{BB}$ as follows. For RBFs $\varphi_k$ centered at the $m_B$ boundary nodes $\bx_k$ we define:
\begin{equation}
    \partial^*\!\varphi_k (\bx_i) := \frac{\partial\varphi_k}{\partial\boldsymbol{n}_k}(\bx_i)
    \label{eq:normal_der_notation_reversed}
\end{equation}

Because of the symmetry of the RBFs, from equation \eqref{eq:normal_phi} we have:
\begin{equation}
    \partial  \varphi_i (\bx_k) =
    \frac{\partial \varphi_i}{\partial \boldsymbol{n}_k}(\bx_k) =
   -\frac{\partial \varphi_k}{\partial \boldsymbol{n}_k}(\bx_i) =
   -\partial^*\!\varphi_k (\bx_i)
    \label{eq:normal_phi_reversed}
\end{equation}
and then from equation \eqref{eq:Schur_Multiple_node} we have:
\begin{equation}
    \mathbf{S}_{BB} = -\big(
    \partial^*\!\boldsymbol{\varphi}_{BB}-\partial^*\!\boldsymbol{\varphi}_{BI}\boldsymbol{\varphi}_{II}^{-1}\boldsymbol{\varphi}_{IB} \big)
    \label{eq:Schur_Multiple_node_T_star}
\end{equation}
where $\partial^*\!\boldsymbol{\varphi}_{BB}$ and $\partial^*\!\boldsymbol{\varphi}_{BI}$ are obtained from $\partial\boldsymbol{\varphi}_{BB}$ and $\partial\boldsymbol{\varphi}_{BI}$, respectively, by replacing $\partial$ with $\partial^*$ and reversing the indexes of RBF centers with the indexes of evaluation points, as follows:
\begin{equation}
    \partial^*\!\boldsymbol{\varphi}_{BB}=
    \begin{bmatrix}
\partial^*\!\varphi_{m_J}(\bx_{m_J}) \mysep \cdots \mysep \partial^*\!\varphi_{m_J}(\bx_{m})   \\
\vdots                                          \mysep \ddots \mysep \vdots \\
\partial^*\!\varphi_{m}(\bx_{m_J})     \mysep \cdots \mysep \partial^*\!\varphi_{m}(\bx_{m})   \\
    \end{bmatrix}
\label{eq:InterpM_pureRBF_MultipleBoundaryNode_PhiBB_star}
\end{equation}

\begin{equation}
    \partial^*\!\boldsymbol{\varphi}_{BI}=
    \begin{bmatrix}
\partial^*\!\varphi_{m_J}(\bx_1)       \mysep \cdots \mysep \partial^*\!\varphi_{m_J}(\bx_{m_I}) \\
\vdots                                          \mysep \ddots \mysep \vdots  \\
\partial^*\!\varphi_m(\bx_1)           \mysep \cdots \mysep \partial^*\!\varphi_m(\bx_{m_I})
    \end{bmatrix}
\label{eq:InterpM_pureRBF_MultipleBoundaryNode_PhiBI_star}
\end{equation}

By following the same arguments presented in section \ref{sss:practical_v}, we have that column $j$ of $\boldsymbol{\varphi}_{II}^{-1}\boldsymbol{\varphi}_{IB}$ in equation \eqref{eq:Schur_Multiple_node_T_star} contains the values of the $m_I$ cardinal functions of the RBF interpolation, based on the $m_I$ internal nodes only, evaluated at the $j^{th}$ boundary node. Row $i$ of $\partial^*\!\boldsymbol{\varphi}_{BI}$ in equation \eqref{eq:InterpM_pureRBF_MultipleBoundaryNode_PhiBI_star}, on the other hand, contains the values of the function $\partial^*\!\varphi_{m_I+i}$, thus centered at the $i^{th}$ boundary node, and evaluated at the $m_I$ internal nodes. Therefore, entry $(i,j)$ of matrix $\partial^*\!\boldsymbol{\varphi}_{BI}\boldsymbol{\varphi}_{II}^{-1}\boldsymbol{\varphi}_{IB}$ in equation \eqref{eq:Schur_Multiple_node_T_star} is the extrapolation at the $j^{th}$ boundary node of the function $\partial^*\!\varphi_{m_I+i}$ based on its values at the internal nodes only. On the other hand, entry $(i,j)$ of matrix $\partial^*\!\boldsymbol{\varphi}_{BB}$ in equation \eqref{eq:InterpM_pureRBF_MultipleBoundaryNode_PhiBB_star} is the actual value of $\partial^*\!\varphi_{m_I+i}$ at the $j^{th}$ boundary node. Given the definition of $\partial^*$, we conclude that entry $(i,j)$ of $\mathbf{S}_{BB}$ in equation \eqref{eq:Schur_Multiple_node_T_star} can be interpreted as the difference between two values: the value of the normal derivative, intended with respect to the normal at the $i^{th}$ boundary node, of the RBF centered at the same $i^{th}$ boundary node and the extrapolated value of the same normal derivative, obtained through a RBF interpolation based on the $m_I$ internal nodes only, both evaluated at the $j^{th}$ boundary node. Similarly to the case with one boundary node, the previous interpretation tells us that the interpolation matrix $\mathbf{M}$ in equation \eqref{eq:InterpM_pureRBF_MultipleBoundaryNode} will be qualitatively well-conditioned if the additional information provided by the Neumann BCs, i.e., the last $m_B$ rows in $\mathbf{M}$, are not redundant with the ability of the RBF interpolation (based on the internal nodes only) to provide a good extrapolation at the boundary nodes of the normal derivatives of the RBFs centered at the same boundary nodes.

\subsubsection{Dependance upon the normals}
Equation \eqref{eq:Schur_Multiple_node} can be written as follows:
\begin{equation}
    \mathbf{S}_{BB} = \partial\boldsymbol{\varphi}_{BB}-\partial\boldsymbol{\varphi}_{BI}\bar{\boldsymbol{\psi}}
    \label{eq:Schur_Multiple_node_cardinal}
\end{equation}
where $\bar{\boldsymbol{\psi}}=\boldsymbol{\varphi}_{II}^{-1}\boldsymbol{\varphi}_{IB}$, already defined in equation \eqref{eq:OneBoundaryNode_psi_at_boundary}, is now a $m_I\times m_B$ matrix in the case of $m_B$ boundary nodes.

Despite the matrices in equation \eqref{eq:Schur_Multiple_node_cardinal} could be handled by standard linear algebra operations, it is anyhow convenient to use the notation introduced in section \ref{App:s1} for simplicity as follows. By denoting with $\boldsymbol{d}_{i,k}$ the following term appearing in equation \eqref{eq:normal_phi}:
\begin{equation}
    \boldsymbol{d}_{i,k}=\Phi'(r_{i,k})\,\boldsymbol{e}_{i,k}
    \label{eq:dik}
\end{equation}
the normal derivatives in $\partial\boldsymbol{\varphi}_{BB}$ and in $\partial\boldsymbol{\varphi}_{BI}$ of equation \eqref{eq:InterpM_pureRBF_MultipleBoundaryNode_splitting} can be written as simple dot products:
\begin{equation}
    \partial \varphi_i (\bx_k) = 
    \boldsymbol{d}_{i,k} \cdot \boldsymbol{n}_k
    \label{eq:normal_phi_d}
\end{equation}
and therefore $\partial\boldsymbol{\varphi}_{BB}$ and $\partial\boldsymbol{\varphi}_{BI}$ can be written using the operator $H$, as defined in equation \eqref{App1_opH}, as follows:
\begin{equation}
    \begin{split}
    \partial\boldsymbol{\varphi}_{BB}&=H(\mathcal{D}_{BB},\mathcal{N})\\
    \partial\boldsymbol{\varphi}_{BI}&=H(\mathcal{D}_{BI},\mathcal{N})
    \label{eq:Schur_Multiple_node_cardinal_dmat}
\end{split}
\end{equation}

In the previous equations $\mathcal{D}_{BB}\in\Rd{m_B}{m_B}$ and $\mathcal{D}_{BI}\in\Rd{m_B}{m_I}$ are $d$-matrices storing purely geometric information:
\begin{equation}
    \mathcal{D}_{BB}=
    \begin{bmatrix}
\boldsymbol{d}_{m_J,m_J} \mysep \cdots \mysep \boldsymbol{d}_{m,m_J} \\
\vdots                   \mysep \ddots \mysep \vdots \\
\boldsymbol{d}_{m_J,m} \mysep \cdots \mysep \boldsymbol{d}_{m,m} \\
    \end{bmatrix} \, , \quad
    \mathcal{D}_{BI}=
    \begin{bmatrix}
\boldsymbol{d}_{1,m_J} \mysep \cdots \mysep \boldsymbol{d}_{m_I,m_J} \\
\vdots                   \mysep \ddots \mysep \vdots \\
\boldsymbol{d}_{1,m} \mysep \cdots \mysep \boldsymbol{d}_{m_I,m} \\
    \end{bmatrix}
    \label{eq:Schur_Multiple_DBB}
\end{equation}
and $\mathcal{N}=\{\boldsymbol{n}_{m_J},\ldots,\boldsymbol{n}_{m}\}\in\Rdvec{m_B}$ is the $d$-matrix storing the unit normals at the $m_B$ boundary nodes. By considering the vector entries of $\mathcal{D}_{BB}$, we see that $\boldsymbol{d}_{i,k}=-\boldsymbol{d}_{k,i}$ for $k\neq i$ because of equation \eqref{eq:dik}, while the diagonal entries are null vectors because of equation \eqref{eq:PhiCondition}, i.e., $\mathcal{D}_{BB}$ is an antisymmetric $d$-matrix.

By using the previous notations and properties \eqref{App1_prop_sum}-\eqref{App1_prop_prod}, equation \eqref{eq:Schur_Multiple_node_cardinal} becomes:
\begin{equation}
    \mathbf{S}_{BB} = H(\mathcal{D}_{BB},\mathcal{N})-H(\mathcal{D}_{BI},\mathcal{N})\bar{\boldsymbol{\psi}} = H(\underbrace{\mathcal{D}_{BB}-\mathcal{D}_{BI}\bar{\boldsymbol{\psi}}}_{\displaystyle \mathcal{G}_{BB}},\mathcal{N})
    \label{eq:Schur_Multiple_node_cardinal_H}
\end{equation}
where $\mathcal{G}_{BB}=(\boldsymbol{g}_{ij})\in\Rd{m_B}{m_B}$. Given the interpretation of $\mathbf{S}_{BB}$ given in section \ref{sss:qualitative_multiple}, vector $\boldsymbol{g}_{ij}$ of $\mathcal{G}_{BB}$ represents therefore the difference between the gradient of the RBF centered at the $i^{th}$ boundary node and the extrapolated value of the same gradient, obtained through a RBF interpolation based on the $m_I$ internal nodes only, both evaluated at the $j^{th}$ boundary node. When considering increasingly flat RBFs \cite{fornberg2004some}, i.e., $\varepsilon\to 0$ for the RBFs reported in Table \ref{tab:RBFs}, the previous RBF extrapolation gets more and more accurate and therefore $\boldsymbol{g}_{ij}\to\boldsymbol{0}$ if the internal stencil is unisolvent for the lowest-order interpolating polynomial. Hence, in the last case, given the definition of $\mathcal{G}_{BB}$ in equation \eqref{eq:Schur_Multiple_node_cardinal_H}, we have $\mathcal{D}_{BI}\bar{\boldsymbol{\psi}}\to\mathcal{D}_{BB}$ in the limit $\varepsilon\to 0$, i.e., $\mathcal{D}_{BI}\bar{\boldsymbol{\psi}}$ tends to an antisymmetric $d$-matrix as well.

Equation \eqref{eq:Schur_Multiple_node_cardinal_H} thus accounts for the explicit influence of the boundary normals on $\mathbf{S}_{BB}$ and, by equation \eqref{eq:Schur_Multiple_node_det}, $\mathbf{S}_{BB}$ in turn influences the conditioning of the whole interpolation matrix.

From equation \eqref{eq:Schur_Multiple_node_cardinal_H}, the determinant of $\mathbf{S}_{BB}$ can be explicitly written as follows:
\begin{equation}
    \det(\mathbf{S}_{BB}) = \det\!\big(H(\mathcal{G}_{BB},\mathcal{N})\big)=
    \begin{vmatrix}
\boldsymbol{g}_{1,1  }\cdot\bar{\boldsymbol{n}}_{1} & \cdots & \boldsymbol{g}_{1,m_B}\cdot\bar{\boldsymbol{n}}_{1} \\
\vdots & \ddots & \vdots \\
\boldsymbol{g}_{m_B,1  }\cdot\bar{\boldsymbol{n}}_{m_B} & \cdots & \boldsymbol{g}_{m_B,m_B}\cdot\bar{\boldsymbol{n}}_{m_B} 
    \end{vmatrix}
    \label{eq:Schur_Multiple_node_cardinal_det_Gbb}
\end{equation}
where $\bar{\boldsymbol{n}}_{i}=\boldsymbol{n}_{m_I+i}$, i.e., the unit normal at the $i^{th}$ boundary node, is introduced for simplicity of notation. We note that in the case with $m_B=1$ boundary nodes, equation \eqref{eq:Schur_Multiple_node_cardinal_det_Gbb} becomes:
\begin{equation}
    \det(\mathbf{S}_{BB}) = S_{BB} = \boldsymbol{g}_{1,1}\cdot\bar{\boldsymbol{n}}_1
    \label{eq:Schur_Multiple_node_cardinal_det_Gbb_mB_single}
\end{equation}
which corresponds to equation \eqref{eq:OneBoundaryNode_DetCofactors_v} with $\boldsymbol{v}=\det(\boldsymbol{\varphi}_{II})\boldsymbol{g}_{1,1}$ because of equation \eqref{eq:Schur_1node_det}.

Finally we note that even in the case with $m_B>1$ boundary nodes, the determinant of the interpolation matrix can still be expressed in the form of a dot product with respect to a specific normal $\bar{\boldsymbol{n}}_i$ as in equation \eqref{eq:OneBoundaryNode_DetCofactors_v}, i.e., $\det(\mathbf{M})=\boldsymbol{v}_i\cdot\bar{\boldsymbol{n}}_i$ where $\boldsymbol{v}_i$ therefore depends upon the remaining $m_B-1$ normals $\bar{\boldsymbol{n}}_{j\neq i}$. This fact can be derived by directly expanding the determinant in both equations \eqref{eq:InterpM_pureRBF_MultipleBoundaryNode_splitting} and \eqref{eq:Schur_Multiple_node_cardinal_det_Gbb} with respect to the row corresponding to the $i^{th}$ boundary node.

\subsubsection{Optimal directions}
Given the formulation of equation \eqref{eq:Schur_Multiple_node_cardinal_det_Gbb}, one may ask which are the $m_B$ unit vectors $\bar{\boldsymbol{n}}_{i}$ that maximize $\det(\mathbf{S}_{BB})$, i.e., the optimal directions of the $m_B$ normals leading to a well-conditioned interpolation matrix. Such optimal directions can therefore be obtained by solving the following constrained maximization problem:
\begin{equation}
    \begin{array}{lll}
        \text{maximize} && \det\!\big(H(\mathcal{G}_{BB},\mathcal{N})\big)\\[.25em]
        \text{subject to} && \|\bar{\boldsymbol{n}}_i\|_2^2=1,\quad i=1,\ldots,m_B\\
    \end{array}
    \label{eq:constr_opt_mB_nodes}
\end{equation}
where $\mathcal{N}=\{\bar{\boldsymbol{n}}_1,\ldots,\bar{\boldsymbol{n}}_{m_B}\}\in\Rdvec{m_B}$.

\begin{figure}[t]
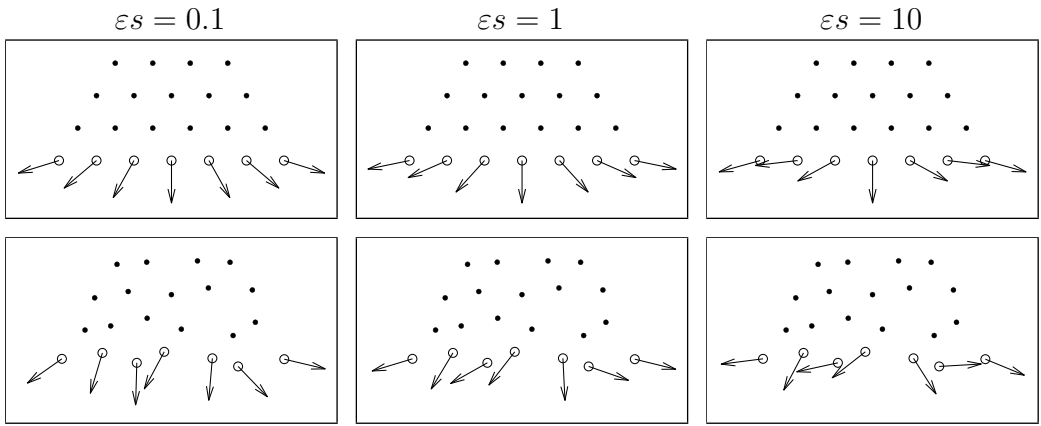

    \def\SpaceBelowText{.25em}
    \def\ImagesWidth{.32\textwidth}
    \def\SpaceBelowFirstRow{0.5em}
    \setlength\fboxsep{0pt}
    \setlength\fboxrule{0.5pt}
    \centering

    \begin{subfigure}[b]{\ImagesWidth}
        \centering
        $\varepsilon s=0.1$\\[\SpaceBelowText]
        \boxedgraphics{optimal_normals/opt_eps01.pdf}
    \end{subfigure}
    \hfill
    \begin{subfigure}[b]{\ImagesWidth}
        \centering
        $\varepsilon s=1$\\[\SpaceBelowText]
        \boxedgraphics{optimal_normals/opt_eps1.pdf}
    \end{subfigure}
    \hfill
    \begin{subfigure}[b]{\ImagesWidth}
        \centering
        $\varepsilon s=10$\\[\SpaceBelowText]
        \boxedgraphics{optimal_normals/opt_eps10.pdf}
    \end{subfigure}\\[\SpaceBelowFirstRow]
    \begin{subfigure}[b]{\ImagesWidth}
        \centering
        \boxedgraphics{optimal_normals/opt_eps01_rand.pdf}
    \end{subfigure}
    \hfill
    \begin{subfigure}[b]{\ImagesWidth}
        \centering
        \boxedgraphics{optimal_normals/opt_eps1_rand.pdf}
    \end{subfigure}
    \hfill
    \begin{subfigure}[b]{\ImagesWidth}
        \centering
        \boxedgraphics{optimal_normals/opt_eps10_rand.pdf}
    \end{subfigure}
    \caption{optimal directions for the reference stencil of Figure \ref{FIG:ref_stencil} (top row) and for the same stencil with random nodal displacements (bottom row).}
    \label{FIG_opt_normals}
\end{figure}

Problem \eqref{eq:constr_opt_mB_nodes} can be solved by the Lagrange multipliers method as presented in \ref{App:s2}. Such optimal directions are shown in Figure \ref{FIG_opt_normals} for the reference stencil of Figure \ref{FIG:ref_stencil} when using MQ RBF with different values of the shape factor $\varepsilon$. In order to highlight the effect of the relative position of the nodes, optimal directions are also shown when the reference stencil is perturbed with random nodal displacements.

In the case of the reference stencil, the optimal directions are approximately pointing towards some reference mid point of the stencil, which however depends upon the shape factor $\varepsilon$, analogously to the case with one boundary node depicted in Figure \ref{FIG_internal_weights_v}. When $\varepsilon s=0.1$, the long range influence of the internal nodes on the optimal normals seems to dominate: as a result, the optimal directions have a radial pattern originating from the geometrical center of the internal nodes. For larger shape factors, i.e., $\varepsilon s=1$ and especially $\varepsilon s=10$, conversely, the short range mutual influence of the boundary nodes seems to dominate, leading to a more flat pattern, again in accordance with the qualitative arguments presented in section \ref{sss:Possible remedies} related to the coefficients $w_{j,m}$ shown in Figure \ref{FIG_internal_weights_v}. In the case of the perturbed stencil, however, the optimal directions are quite unpredictable, as there are boundary nodes close to each other with quite different optimal directions, confirming again that a simple geometric approximation of the optimal directions is not feasible.

\subsection{Influence of the polynomial part}
In order to highlight the effect of the polynomial part of the interpolant \eqref{eq:polyRBF} on the previous results, let us consider again the interpolation matrix $\mathbf{M}$ in equation \eqref{eq:interpSystem}. By using the same matrix splitting employed in equation \eqref{eq:InterpM_pureRBF_MultipleBoundaryNode_splitting}, submatrix $\boldsymbol{\varphi}_{BC}$ initially defined in equation \eqref{eq:phimat} can be expressed as:
\begin{equation}
    \boldsymbol{\varphi}_{BC}=
    \begin{bmatrix}
    \begin{array}{c|c}
\boldsymbol{\varphi}_{II} & \boldsymbol{\varphi}_{IB} \\ \hline
\partial\boldsymbol{\varphi}_{BI} & \partial\boldsymbol{\varphi}_{BB}
    \end{array}
    \end{bmatrix}
\label{eq:InterpPhiBC_wpoly_MultipleBoundaryNode_splitting}
\end{equation}
while matrices $\mathbf{P}_{BC}$ and $\mathbf{P}^T$, defined in equations \eqref{eq:pmat} and \eqref{eq:ptmat}, can be expressed as follows: 
\begin{equation}
    \mathbf{P}_{BC}=
    \begin{bmatrix}
    \begin{array}{ccc}
p_1(\bx_1)               \mysep \cdots \mysep p_q(\bx_1)             \\
\vdots                   \mysep \ddots \mysep \vdots                 \\
p_1(\bx_{m_I})           \mysep \cdots \mysep p_q(\bx_{m_I})         \\ \hline
\partial p_1 (\bx_{m_J}) \mysep \cdots \mysep \partial p_q(\bx_{m_J})\\
\vdots                   \mysep \ddots \mysep \vdots                 \\
\partial p_1(\bx_m)      \mysep \cdots \mysep \partial p_q(\bx_m)
    \end{array}
    \end{bmatrix}
    =
    \begin{bmatrix}
    \begin{array}{c}
\mathbf{P}_I \\ \hline
\partial \mathbf{P}_B
    \end{array}
    \end{bmatrix}
    \label{eq:pmat_splitting}
\end{equation}

\begin{equation}
    \mathbf{P}^T=
    \begin{bmatrix}
    \begin{array}{ccc|ccc}
p_1(\bx_1) \mysep \cdots \mysep p_1(\bx_{m_I}) & p_1(\bx_{m_J}) \mysep \cdots \mysep p_1(\bx_m) \\
\vdots     \mysep \ddots \mysep \vdots         & \vdots         \mysep \ddots \mysep \vdots     \\
p_q(\bx_1) \mysep \cdots \mysep p_q(\bx_{m_I}) & p_q(\bx_{m_J}) \mysep \cdots \mysep p_q(\bx_m)
\end{array}
    \end{bmatrix}
    =
    \begin{bmatrix}
    \begin{array}{c|c}
\mathbf{P}^T_I & \mathbf{P}^T_B
    \end{array}
    \end{bmatrix}
    \label{eq:ptmat_splitting}
\end{equation}

By rearranging the rows and the columns of $\mathbf{M}$ according to internal nodes, polynomial terms and boundary nodes, respectively, the interpolation matrix in block form becomes:
\begin{equation}
    \mathbf{M}=
    \begin{bmatrix}
    \begin{array}{cc|c}
        \boldsymbol{\varphi}_{II} & \mathbf{P}_I & \boldsymbol{\varphi}_{IB} \\[.25em]
        \mathbf{P}^T_I            & \boldsymbol{0} & \mathbf{P}^T_B            \\ \hline
        \partial\boldsymbol{\varphi}_{BI} & \partial\mathbf{P}_B & \partial\boldsymbol{\varphi}_{BB}
    \end{array}
    \end{bmatrix}
    =
    \begin{bmatrix}
    \begin{array}{c|c}
        \mathbf{M}_{II} & \mathbf{Q}_{IB} \\ \hline
        \partial\mathbf{Q}_{BI} & \partial\boldsymbol{\varphi}_{BB}
    \end{array}
    \end{bmatrix}
    \label{eq:M_wpoly_splitting}
\end{equation}

Similarly to the approach employed in section \ref{sss:RBF_multiple_bnd_nodes_S}, we consider the Schur complement $\mathbf{S}_{BB}$ of the $\mathbf{M}_{II}$ block of $\mathbf{M}$:
\begin{equation}
    \mathbf{S}_{BB}=
    \partial\boldsymbol{\varphi}_{BB} -
    \underbrace{\partial\mathbf{Q}_{BI} \mathbf{M}_{II}^{-1} \mathbf{Q}_{IB}}_{\displaystyle\mathbf{W}_{BB}}
    \label{eq:M_wpoly_Schur}
\end{equation}
where $\mathbf{M}_{II}$ is not singular as long as the internal stencil is unisolvent for the chosen degree of the polynomial augmentation.

By comparing equation \eqref{eq:M_wpoly_Schur} with the corresponding equation \eqref{eq:Schur_Multiple_node}, obtained without polynomial augmentation, we see that the only difference is in the double matrix-product term $\mathbf{W}_{BB}$, that in the case of polynomial augmentation can be written as follows:
\begin{equation}
    \mathbf{W}_{BB}  =
    \begin{bmatrix}
        \partial\boldsymbol{\varphi}_{BI} & \partial\mathbf{P}_B
    \end{bmatrix}
    \begin{bmatrix}
        \boldsymbol{\varphi}_{II} & \mathbf{P}_I \\[.25em]
        \mathbf{P}^T_I            & \boldsymbol{0} 
    \end{bmatrix}^{-1}
    \begin{bmatrix}
        \boldsymbol{\varphi}_{IB} \\[.25em]
        \mathbf{P}^T_B  
    \end{bmatrix}
    =
    \partial\boldsymbol{\varphi}_{BI}\mathbf{K}_1 + \partial\mathbf{P}_B\mathbf{K}_2
    \label{eq:WBB_wpoly}
\end{equation}
where $\mathbf{K}_1$ and $\mathbf{K}_2$ can be obtained by simple algebraic operations as follows:
\begin{equation}
\begin{split}
    \mathbf{K}_1 &= \boldsymbol{\varphi}_{II}^{-1}
    \big(
    \boldsymbol{\varphi}_{IB} - \mathbf{P}_I\mathbf{K}_2
    \big) \\
    \mathbf{K}_2 &=
    \underbrace{\big(
    \mathbf{P}_I^T \boldsymbol{\varphi}_{II}^{-1} \mathbf{P}_I
    \big)^{-1}}_{\displaystyle\mathbf{S}}
    \big(
    \boldsymbol{\varphi}_{IB}^T \boldsymbol{\varphi}_{II}^{-1} \mathbf{P}_I - \mathbf{P}_B
    \big)^T
    \end{split}
    \label{eq:K12_wpoly}
\end{equation}

Further algebraic manipulations lead to the following form:
\begin{equation}
    \mathbf{W}_{BB}  =
    \partial\boldsymbol{\varphi}_{BI}\boldsymbol{\varphi}_{II}^{-1}\boldsymbol{\varphi}_{IB} - 
    \underbrace{\big(
    \partial\boldsymbol{\varphi}_{BI} \boldsymbol{\varphi}_{II}^{-1} \mathbf{P}_I
    -\partial\mathbf{P}_B
    \big)
    \mathbf{K}_2}_{\displaystyle\Delta \mathbf{W}}
    \label{eq:WBB_wpoly_end}
\end{equation}
which, compared to the corresponding term $\partial\boldsymbol{\varphi}_{BI}\boldsymbol{\varphi}_{II}^{-1}\boldsymbol{\varphi}_{IB}$ in equation \eqref{eq:Schur_Multiple_node}, reveals that the difference in the Schur complement of the interpolation matrix, with and without polynomial augmentation, is solely due to the term $\Delta \mathbf{W}$:
\begin{equation}
    \Delta \mathbf{W}  =
    \underbrace{\big(
    \partial\boldsymbol{\varphi}_{BI} \boldsymbol{\varphi}_{II}^{-1} \mathbf{P}_I
    -\partial\mathbf{P}_B
    \big)}_{\displaystyle\partial \mathbf{E}}
    \mathbf{S}
    \big(\underbrace{
    \boldsymbol{\varphi}_{IB}^T \boldsymbol{\varphi}_{II}^{-1} \mathbf{P}_I - \mathbf{P}_B
    }_{\displaystyle \mathbf{E}}\big)^T
    \label{eq:DeltaW_wpoly_end}
\end{equation}

Similarly to the qualitative interpretations given in section \ref{sss:qualitative_multiple}, we have that entry $(i,j)$ of matrix $\partial \mathbf{E}$ can be interpreted as the difference between two values: the value of the normal derivative, intended with respect to the normal at the $i^{th}$ boundary node, of the $j^{th}$ function of the polynomial base, i.e., $p_j$, and the extrapolated value of the same normal derivative, obtained through a RBF interpolation based on the $m_I$ internal nodes only, both evaluated at the $i^{th}$ boundary node. On the other hand, entry $(i,j)$ of matrix $\mathbf{E}$ can again be interpreted as the difference between two values: the value of $p_j$ and the extrapolated value of the same function $p_j$, obtained through a RBF interpolation based on the $m_I$ internal nodes only, both evaluated at the $i^{th}$ boundary node.

Therefore, if the internal stencil is unisolvent for the chosen degree of the polynomial augmentation, each entry of matrices $\partial \mathbf{E}$ and $\mathbf{E}$ vanishes as $\varepsilon\to 0$ and the same holds for the entries of matrix $\mathbf{S}$ except for the diagonal entry corresponding to the constant monomial of the polynomial base, which tends to a constant value. In the end, matrix $\Delta \mathbf{W}$ vanishes rapidly as $\varepsilon\to 0$, and therefore in this limit the polynomial augmentation does not affect the formulation presented in the previous sections to study the effect of the normals on the interpolation matrix.

\begin{figure}[t]
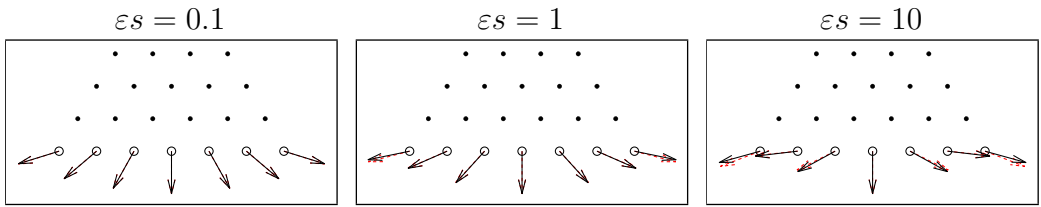

    \def\SpaceBelowText{.25em}
    \def\ImagesWidth{.32\textwidth}
    \setlength\fboxsep{0pt}
    \setlength\fboxrule{0.5pt}
    \centering

    \begin{subfigure}[b]{\ImagesWidth}
        \centering
        $\varepsilon s=0.1$\\[\SpaceBelowText]
        \boxedgraphics{poly_influence/poly_eps01.pdf}
    \end{subfigure}
    \hfill
    \begin{subfigure}[b]{\ImagesWidth}
        \centering
        $\varepsilon s=1$\\[\SpaceBelowText]
        \boxedgraphics{poly_influence/poly_eps1.pdf}
    \end{subfigure}
    \hfill
    \begin{subfigure}[b]{\ImagesWidth}
        \centering
        $\varepsilon s=10$\\[\SpaceBelowText]
        \boxedgraphics{poly_influence/poly_eps10.pdf}
    \end{subfigure}
    \caption{optimal directions for the reference stencil of Figure \ref{FIG:ref_stencil} without polynomial augmentation (black solid vectors) and with polynomial augmentation with degree $P=2$ (red dashed vectors).}
    \label{FIG_poly_influence}
\end{figure}

The influence of the polynomial augmentation, with degree $P=2$, on the optimal directions for the reference stencil of Figure \ref{FIG:ref_stencil} is shown in Figure \ref{FIG_poly_influence} in the case of MQ RBF for different values of the shape parameter $\varepsilon$. In the case $\varepsilon s=0.1$ there is no noticeable difference in the optimal directions with and without polynomial augmentation, confirming the previous theoretical observations that the influence of the polynomial augmentation is negligible for sufficiently small values of $\varepsilon$. For larger values of $\varepsilon$, i.e., in the cases $\varepsilon s=1$ and $\varepsilon s=10$, the differences in the optimal direction start to be visually noticeable, although very small and almost negligible from a practical point of view.

In conclusion we remark that, since in practical cases $\varepsilon s \in [0.2,0.6]$ with MQ RBF \cite{ZAMOLO2020109730} and $m > 2q$ \cite{bayona2017role} in equation \eqref{eq:polyRBF}, the calculation of the optimal directions can be reasonably carried out without taking into account the polynomial terms.
\section{Techniques for an improved interpolation}

\subsection{Approach 1: boundary node selection based on optimal directions}\label{sss:dot_product}
Once the optimal directions $\hat{\boldsymbol{n}}_i$ have been calculated according to equation \eqref{eq:constr_opt_mB_nodes}, it is possible to extend Approach 1 in section \ref{sss:Possible remedies} to the case with multiple boundary nodes, i.e., discarding from the stencil those boundary nodes whose normal vectors are too different from the corresponding optimal normals. This strategy allows to rule out the possibility of a singular interpolation matrix $\mathbf{M}$ by dropping both the rows and columns associated with dangerous nodes.

The boundary node selection can be implemented as an iterative process (see \ref{App_comput_optimal}) that allows to discard the $i^{th}$ boundary node whenever the dot product between the associated normal $\bar{\boldsymbol{n}}_i$ and the corresponding optimal direction $\hat{\boldsymbol{n}}_i$ is less than a certain threshold value $d_{min}$:
\begin{equation}
    \text{discard boundary node $i$ if }|\bar{\boldsymbol{n}}_i \cdot \hat{\boldsymbol{n}}_i| < d_{min}
    \label{eq:dot_product}
\end{equation}
where the absolute value comes from the fact that the reversal of any normal has no effect on the interpolation.

\begin{figure}[t]
    \def\SpaceBelowText{.25em}
    \def\ImagesWidth{.49\textwidth}
    \centering
    \begin{subfigure}[b]{\ImagesWidth}
        \centering
        \includegraphics[width=\textwidth]{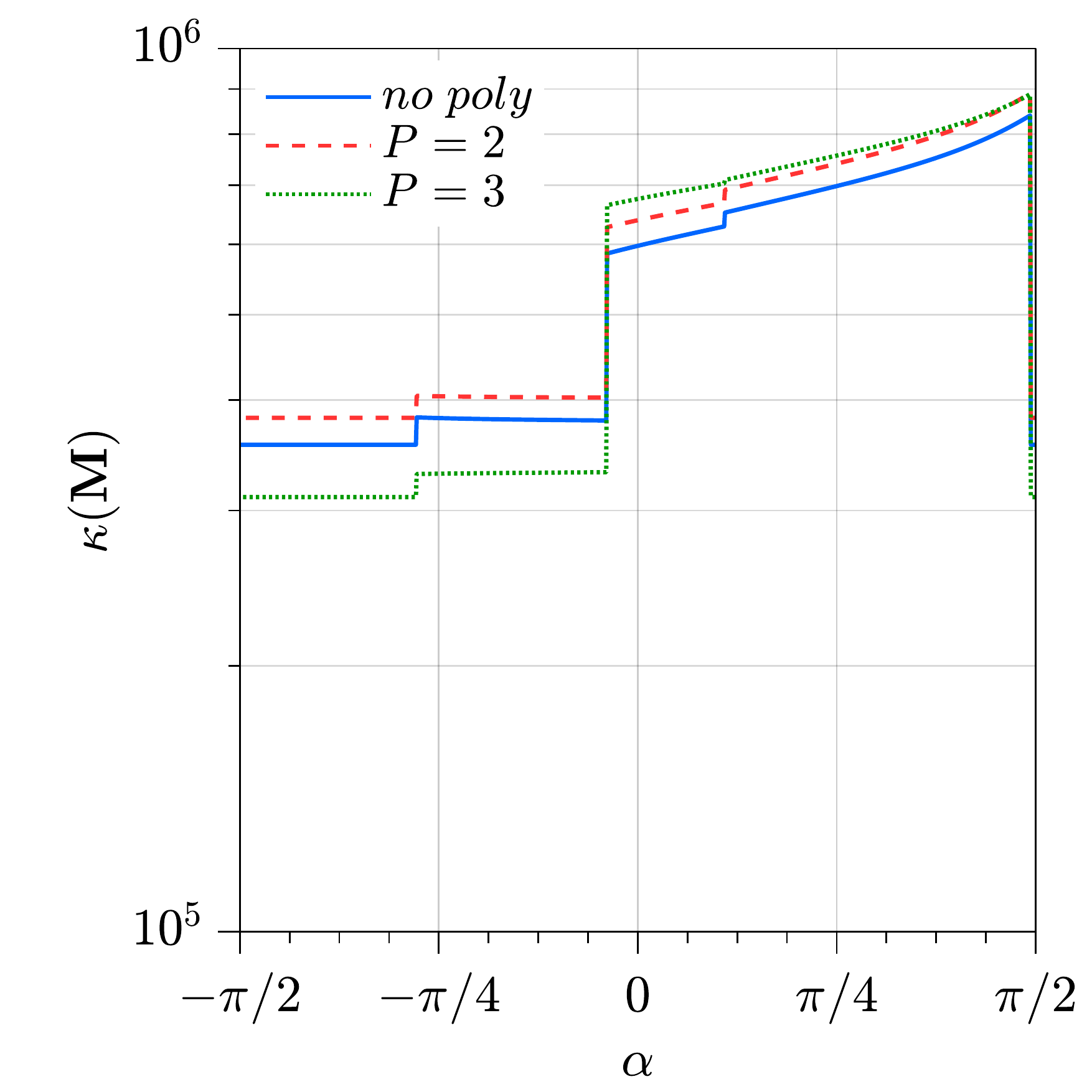}
    \end{subfigure}
    \hfill
    \begin{subfigure}[b]{\ImagesWidth}
        \centering
        \includegraphics[width=\textwidth]{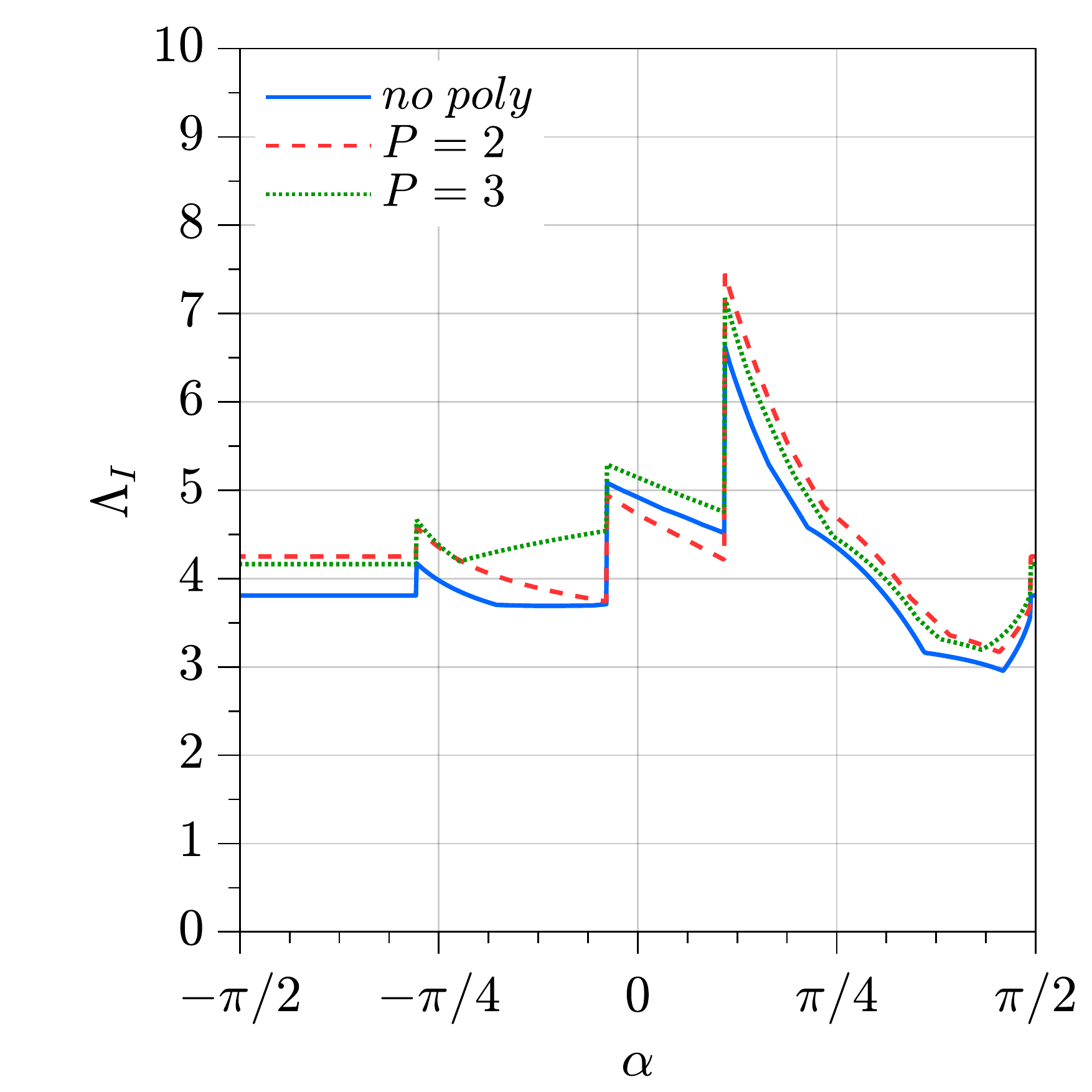}
    \end{subfigure}\\
    \begin{subfigure}[b]{\ImagesWidth}
        \centering
        \includegraphics[width=\textwidth]{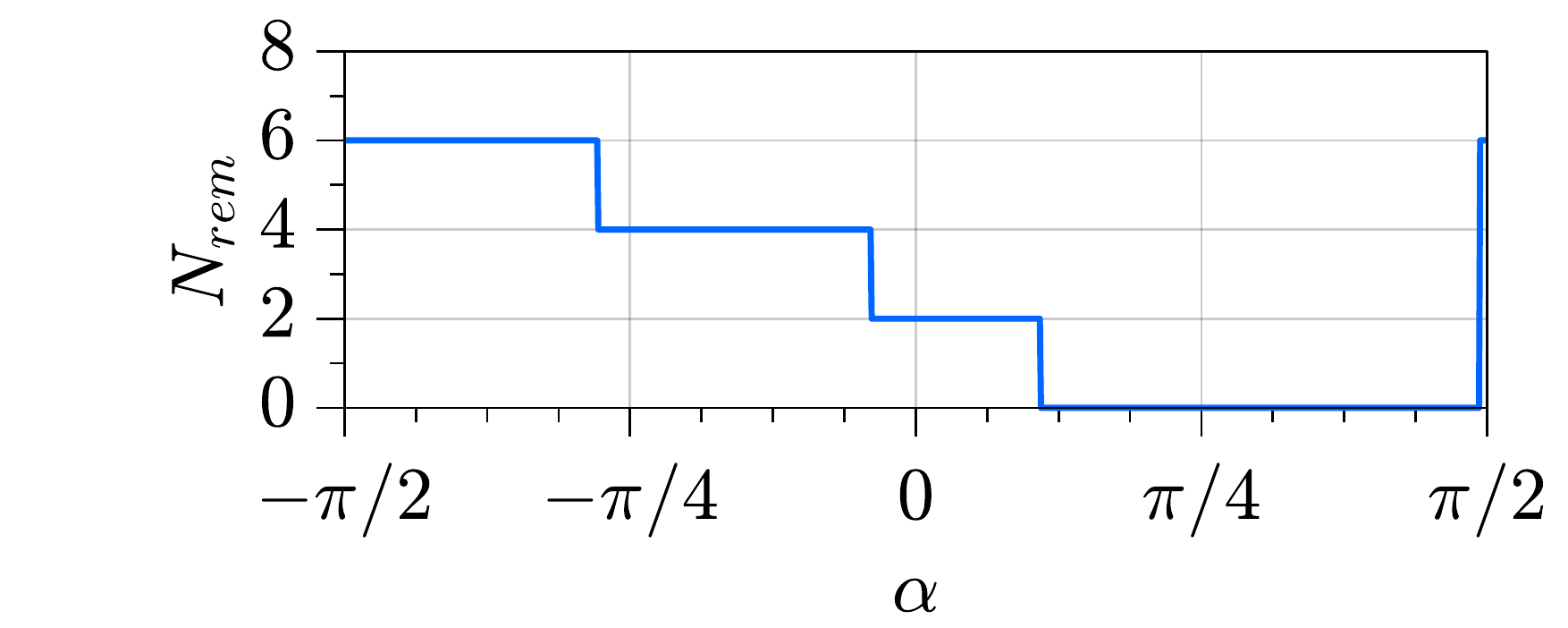}
    \end{subfigure}
    \hfill
    \begin{subfigure}[b]{\ImagesWidth}
        \centering
        \includegraphics[width=\textwidth]{graphics/cardinal_alfa_normals/Removed.pdf}
    \end{subfigure}
    \caption{condition number $\kappa(\mathbf{M})$ of the interpolation matrix (top left), and Lebesgue constant $\Lambda_I$ (top right) for the RBF interpolation on the reference stencil of Figure \ref{FIG:ref_stencil} when the selection of boundary nodes is based on the optimal directions. Bottom row: number of removed nodes $N_{rem}$ (the same plot is reported twice in order to simplify comparisons).}
    \label{FIG_alfa_normals}
\end{figure}

The condition number $\kappa(\mathbf{M})$, the Lebesgue constant $\Lambda_I$ and the number of removed nodes $N_{rem}$ for the reference stencil of Figure \ref{FIG:ref_stencil} are shown in Figure \ref{FIG_alfa_normals} as functions of the angle $\alpha$. From a comparison with Figure \ref{FIG_singularities_alfa} we can see that no more singularities appear for $\alpha < 0$, being now $\kappa(\mathbf{M})$
and $\Lambda_I$ bounded, thus leading to a more stable interpolation. Figure \ref{FIG_alfa_normals} is obtained with MQ RBF ($\varepsilon  s=0.5$) and with a threshold value $d_{min} = 0.6$ in equation \eqref{eq:dot_product}. Lower values are not able to bound the Lebesgue constant $\Lambda_I$ when $\alpha$ approaches the limit values of $\pi/2$ and $-\pi/2$. When 6 nodes are removed, only the central boundary node associated to a vertical normal is included in the stencil. $N_{rem}$ is always even due to the symmetry of the reference stencil. 

\subsection{Approach 2: optimal placement for boundary nodes} \label{ss:optimal_placement}
As hinted in section \ref{sss:Possible remedies}, another strategy for avoiding a singular interpolation matrix consists in moving the stencil nodes in an appropriate way. It follows that any implementation of this strategy involves a modification of the node placement which also depends on the direction of the normals.

In section \ref{App:s3} the optimization of the boundary node position is discussed, it emerges that the interpolation can be improved by projecting on the boundary along the normals those inner nodes which are in the immediate neighborhood.
In in Figure \ref{fig:optimal_bnd} the projected nodes are depicted as asterisks, when the curvature of the boundary is limited they provide a very good approximation for the optimal locations for boundary nodes.
In practice, the adoption of the projected boundary nodes can be implemented with minimal computational cost, here we discuss the attained improvement on the interpolation properties when this is done for the reference stencil of Figure \ref{FIG:ref_stencil}.

\begin{figure}[t]
    \def\ImagesWidth{.49\textwidth}
    \centering
    \begin{subfigure}[b]{\ImagesWidth}
        \centering
        \includegraphics[width=\textwidth]{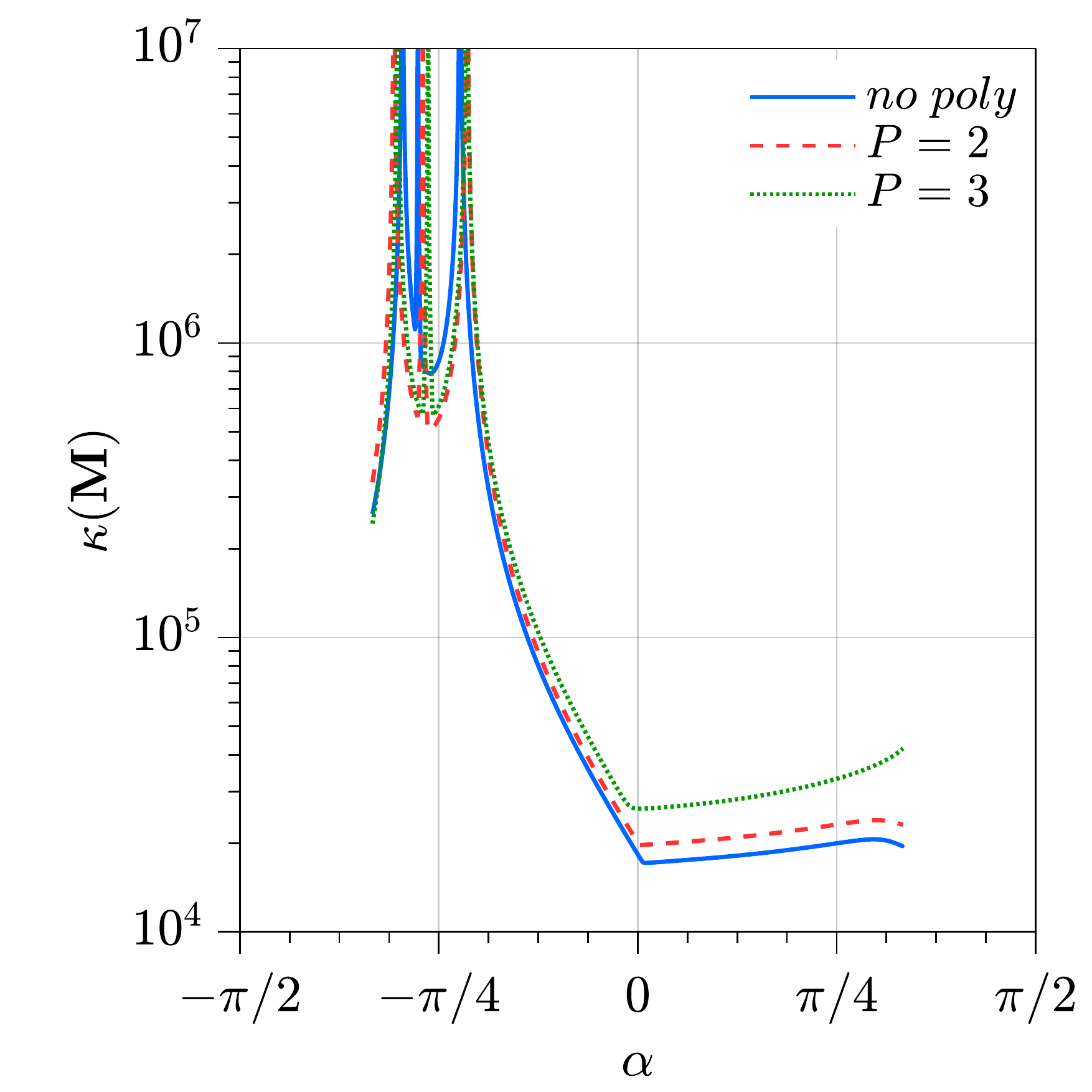}
    \end{subfigure}
    \hfill
    \begin{subfigure}[b]{\ImagesWidth}
        \centering
        \includegraphics[width=\textwidth]{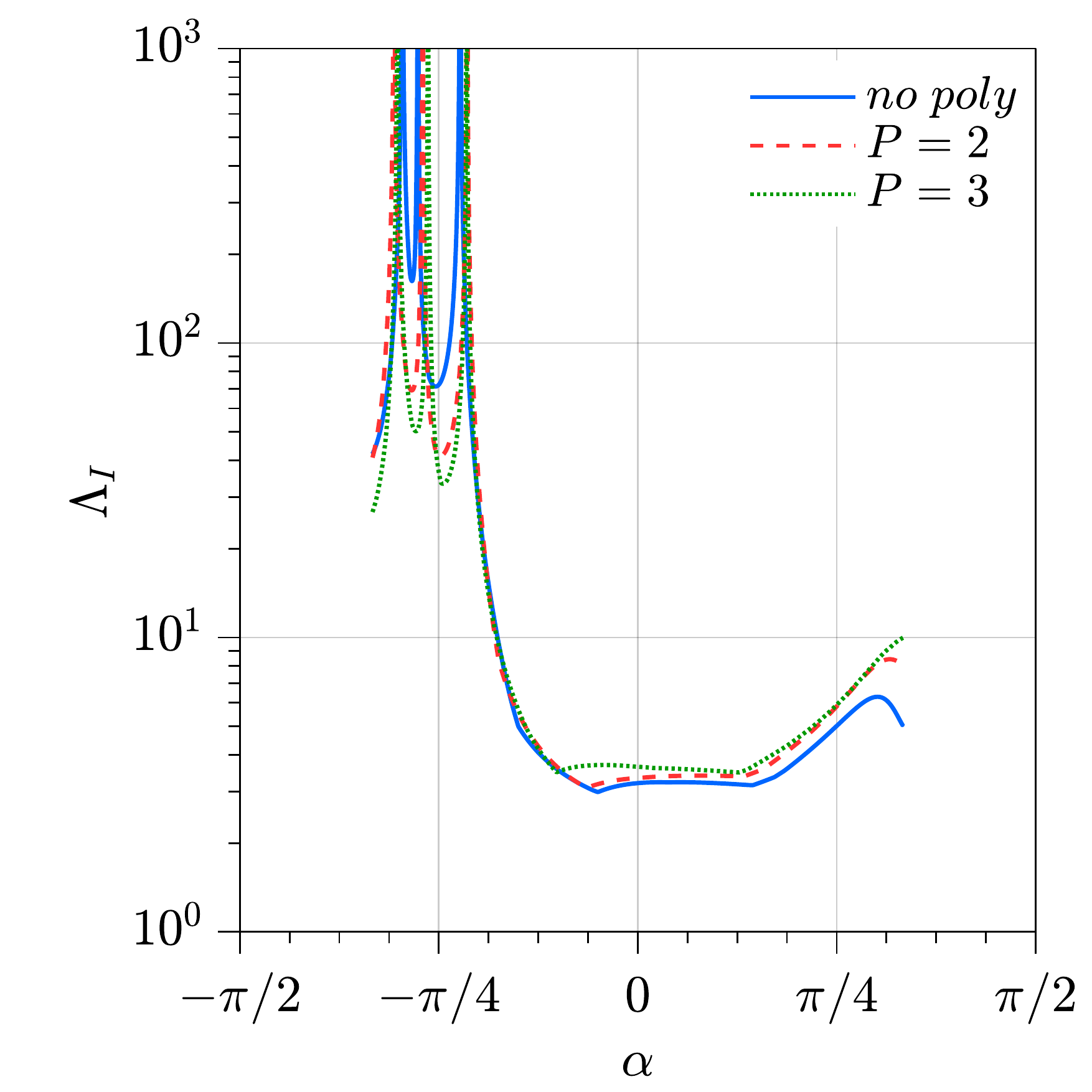}
    \end{subfigure}
    \caption{condition number of the interpolation matrix (left) and Lebesgue constant $\Lambda_I$ (right) for the RBF interpolation on the reference stencil of Figure \ref{FIG:ref_stencil} with projected boundary nodes.}
    \label{fig:Neumann_proj}
\end{figure}

In Figure \ref{fig:Neumann_proj} the condition number $\kappa(\mathbf{M})$ and the Lebesgue constant $\Lambda_I$ are plotted against the angle $\alpha$ when the boundary nodes of the reference stencil are replaced by the projected ones.
The angle $\alpha$ varies in a smaller interval, $\alpha \in  [ -\pi/3 \, , \, \pi/3 ] $ since for large values the projected nodes are placed too close together or too far from one another.
From Figure \ref{fig:Neumann_proj} emerges that, for moderate angles, i.e. approx $\alpha \geq -\pi /8$, the node projection alone is able to avoid the appearance of ill-conditioning problems.
About this, we remark that the reference stencil of Figure \ref{FIG:ref_stencil} can be used to adequately model the case of a curved boundary only for small angles.

Finally, since the use of the projected boundary nodes doesn't impact the solution process it can be applied in conjunction with the selection based on the optimal normals explained above, the two strategies are not mutually exclusive.

\begin{figure}[t!]
    \def\SpaceBelowText{.25em}
    \def\ImagesWidth{.49\textwidth}
    \centering
    \begin{subfigure}[b]{\ImagesWidth}
        \centering
        \includegraphics[width=\textwidth]{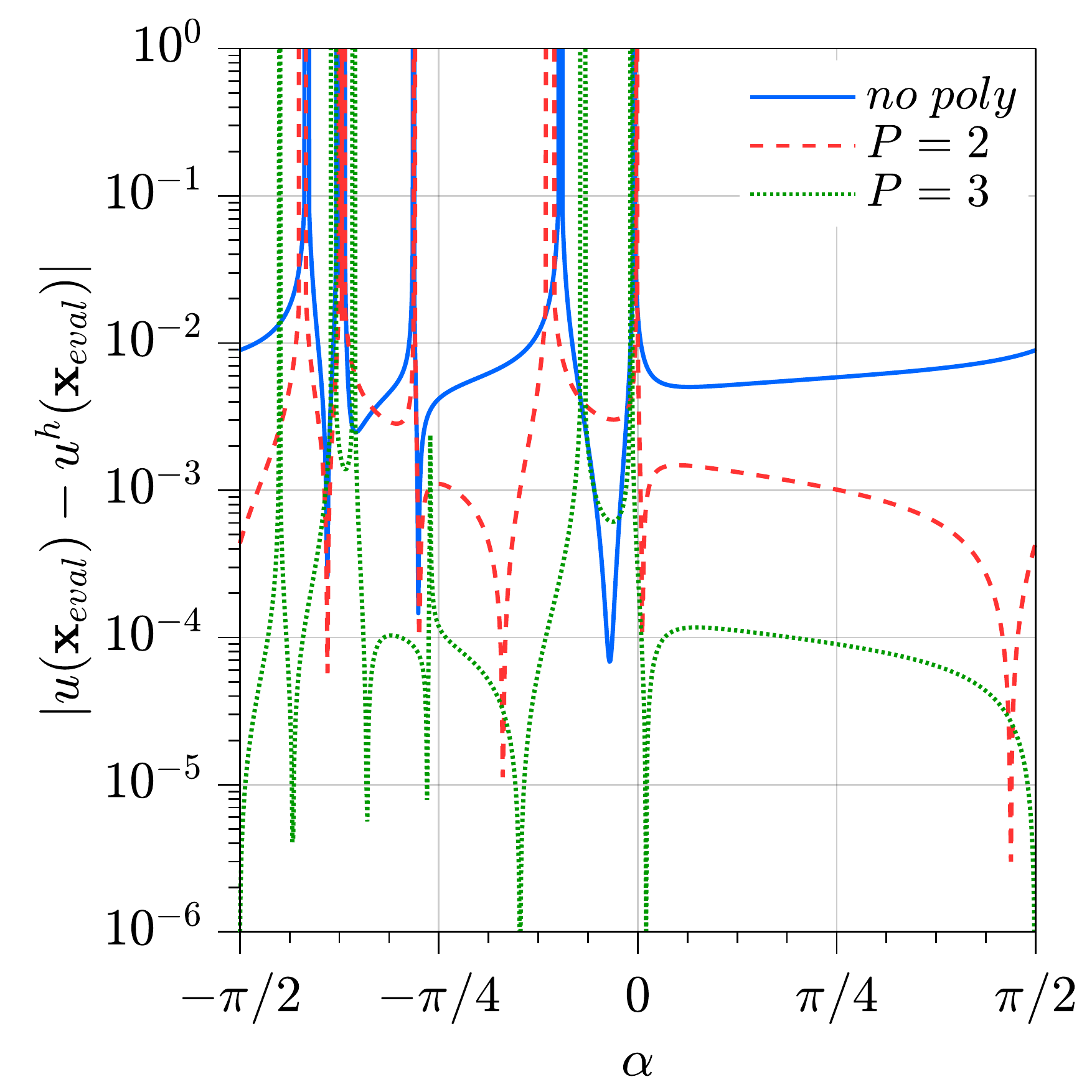}
    \end{subfigure}
    \hfill
    \begin{subfigure}[b]{\ImagesWidth}
        \centering
        \includegraphics[width=\textwidth]{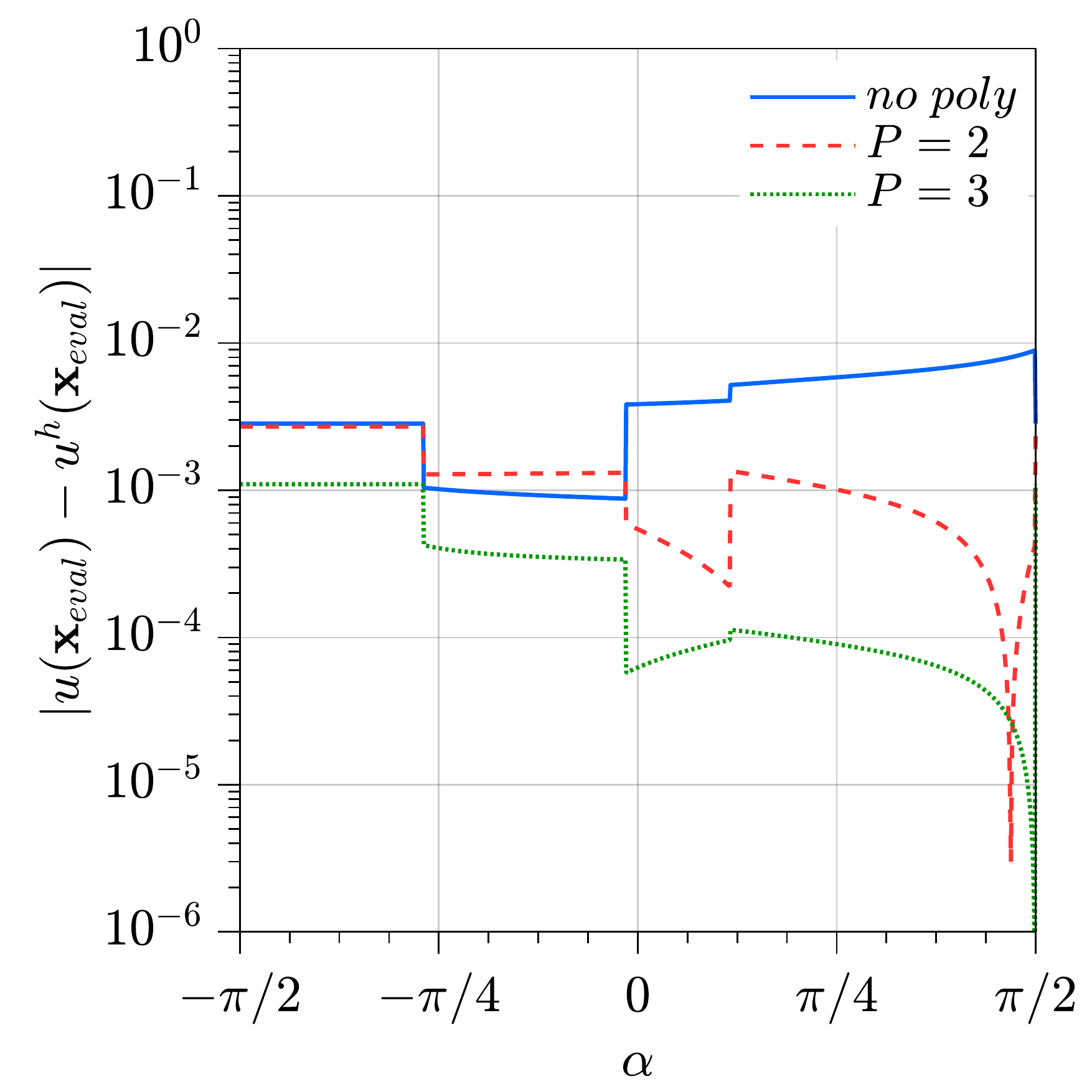}
    \end{subfigure}\\
    \begin{subfigure}[b]{\ImagesWidth}
        \centering
        \includegraphics[width=\textwidth]{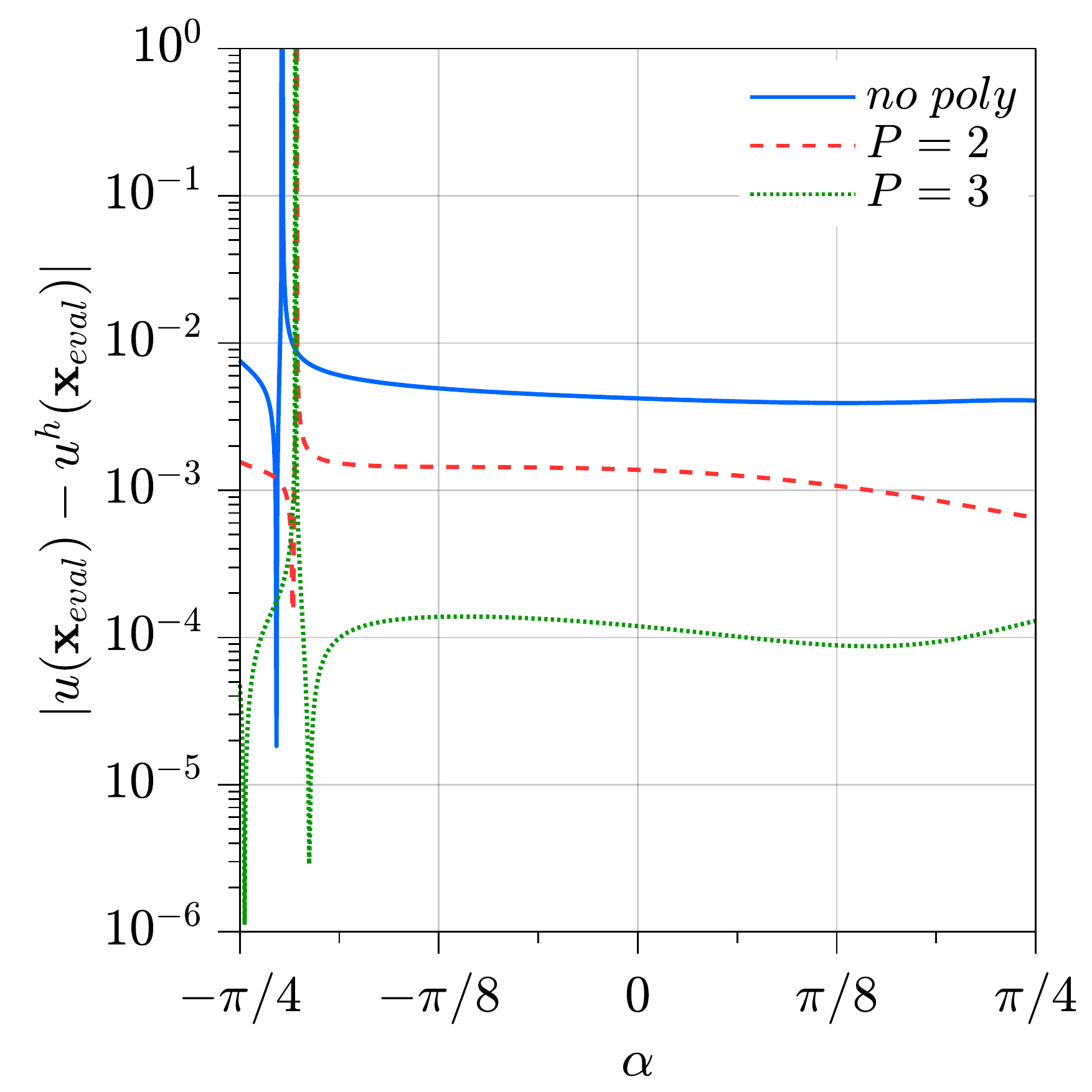}
    \end{subfigure}
    \hfill
    \begin{subfigure}[b]{\ImagesWidth}
        \centering
        \includegraphics[width=\textwidth]{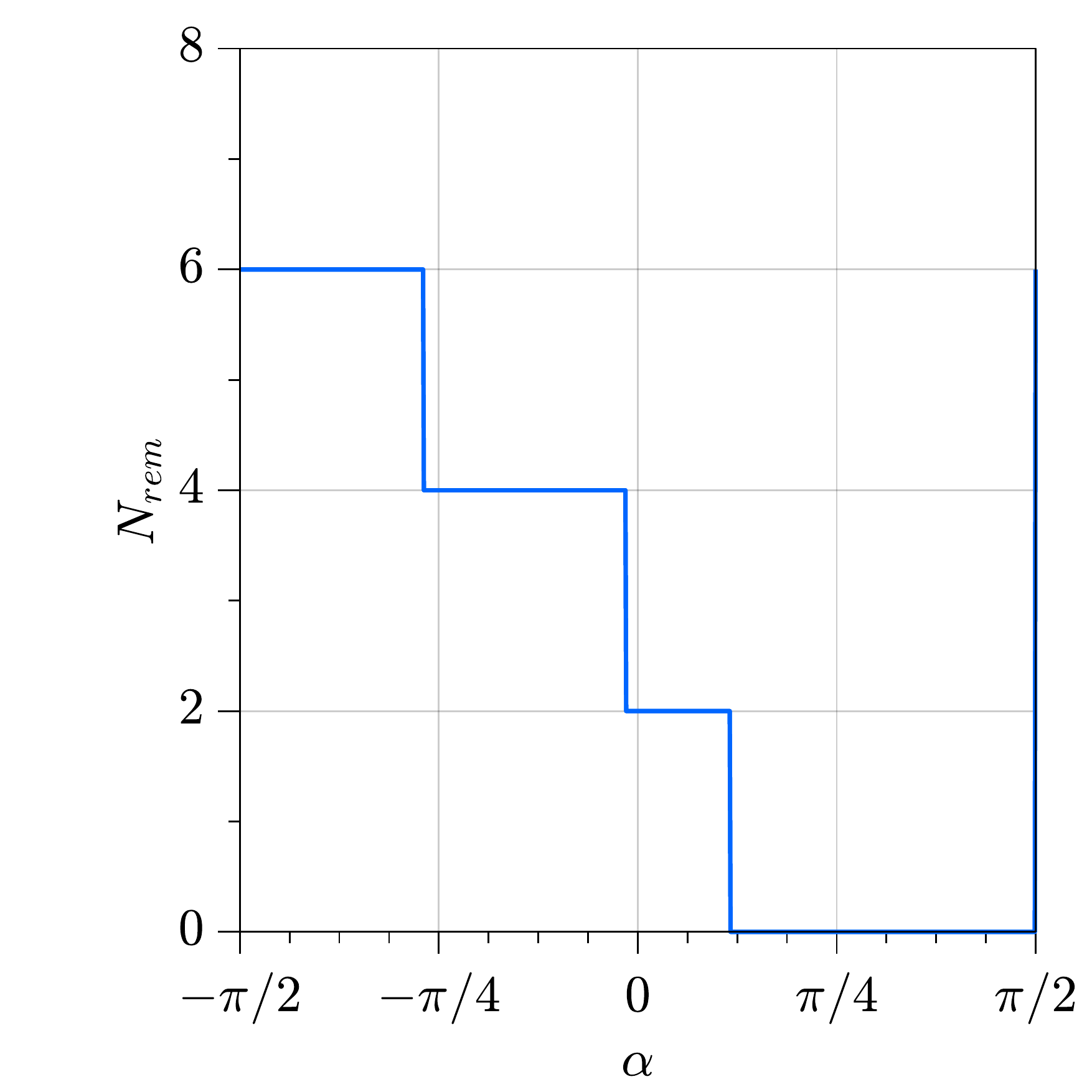}
    \end{subfigure}
    \caption{interpolation error $|u(\boldsymbol{x}_{eval})-u^h(\boldsymbol{x}_{eval})|$ at a given location $\boldsymbol{x}_{eval}=(s/2,s/2)$ for the reference stencil of Figure \ref{FIG:ref_stencil}. No improvement strategy (top left), node selection based on the optimal directions (top right) with the corresponding number $N_{rem}$ of removed boundary nodes (bottom right), boundary node projection (bottom left). The interpolated function is $u=\exp{(x+2y)}$.}
    \label{FIG_alfa_errors}
\end{figure}

In Figure \ref{FIG_alfa_errors} the interpolation error obtained with the different techniques described so far is displayed against the angle $\alpha$. The interpolation error $|u(\boldsymbol{x}_{eval})-u^h(\boldsymbol{x}_{eval})|$ is evaluated at the point $\boldsymbol{x}_{eval}= (s/2,s/2)$, willingly chosen not to be on the symmetry axis, where $s$ is the distance between any two inner nodes of the reference stencil of Figure \ref{FIG:ref_stencil} and the analytic function is defined as $u = \exp{(x+2y)}$. For the selection of boundary nodes based on the optimal directions (right column) the same settings were used as in Figure \ref{FIG_alfa_normals}.

We can see once again that both stabilization techniques succeed in improving the interpolation for small variations of the angle $\alpha$. For the projection of the boundary nodes, $\alpha$ varies within the interval $\alpha \in [-\pi/4 \, , \, \pi/4]$ for the reasons explained above. We can see that the spikes in the errors on the left column, corresponding to singular configurations, are located as in Figures \ref{FIG_singularities_alfa} and \ref{fig:Neumann_proj}, this confirms the reliability of the Lebesgue constant $\Lambda_I$ as a measure for accuracy and stability of the interpolation scheme.
\section{Applications}
\subsection{Stability of the Helmholtz-Hodge Decomposition}\label{ss_Appl_Stability}
In this section the boundary node selection and the node projection strategies, discussed above, are applied to the stabilization of the Helmholtz-Hodge decomposition (HHD) \cite{helmholtz1867lxiii}. The HHD plays a pivotal role in many theoretical and practical applications \cite{bhatia2012helmholtz}, among others, it is at the base of the projection methods for the numerical solution of the incompressible Navier-Stokes equations \cite{CHORIN196712,10.2307_2004575}. According to the HHD, under certain hypothesis, any vector field $\boldsymbol{u}^*$ can be expressed as a sum of the gradient of a scalar potential $\phi$ and a divergence-free vector field $\boldsymbol{u}$:
\begin{equation}
    \boldsymbol{u}^* = \nabla\phi + \boldsymbol{u}
    \label{eq_HHD}
\end{equation}

When solving incompressible Navier-Stokes equations through projection methods, the computation of the unknown velocity field $\boldsymbol{u}$ is carried out iteratively or by marching in time. In any case, each iteration or time step is subdivided into two substeps: first an intermediate velocity $\boldsymbol{u}^*$ is obtained from the momentum equations, then $\boldsymbol{u}$ is obtained from equation \eqref{eq_HHD} as follows:
\begin{equation}
    \boldsymbol{u} = \boldsymbol{u}^* -\nabla\phi 
    \label{eq_HHD_u}
\end{equation}
where the unknown scalar potential $\phi$ is calculated from equation \eqref{eq_HHD_Poisson}, that results from taking the divergence of equation \eqref{eq_HHD}, recalling $\nabla\cdot\boldsymbol{u} = 0$:
\begin{equation}
     \nabla^2\phi = \nabla\cdot\boldsymbol{u}^*
    \label{eq_HHD_Poisson}
\end{equation}

In the case of known velocity at the boundary, e.g., in a cavity, homogeneous Neumann BCs must be enforced in order to solve Poisson equation \eqref{eq_HHD_Poisson}: 
\begin{equation}
    \partial \phi/ \partial\boldsymbol{n} = 0
    \label{eq_HHD_Poisson_BC}
\end{equation}

\begin{figure}[t!]
    \centering
    \includegraphics[width=0.6\textwidth]{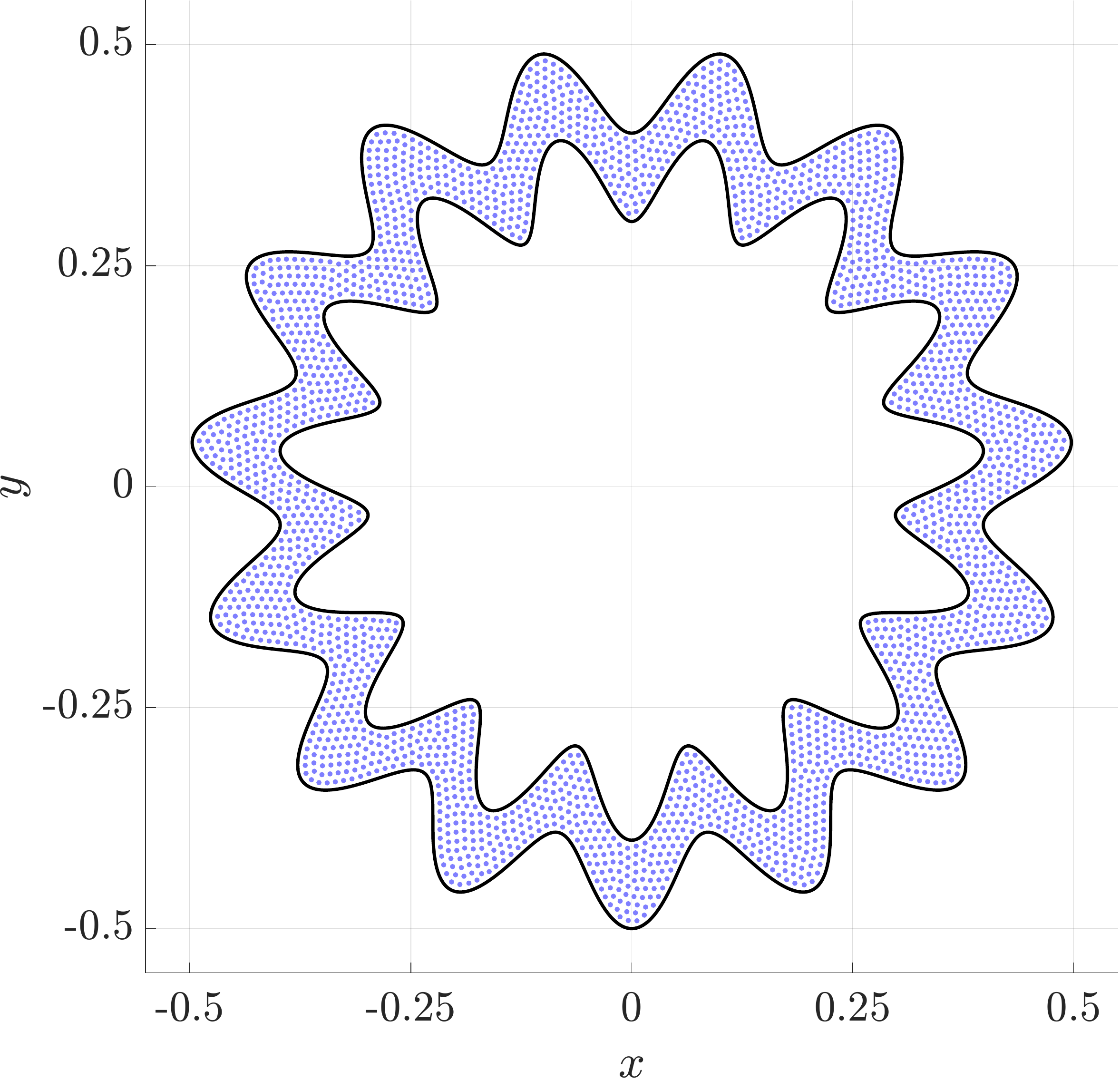}
    \caption{test domain with $N_I\approx 3,300$ internal nodes.}
    \label{fig:geometry_HHD}
\end{figure}

Since the HHD expressed by equations \eqref{eq_HHD_u}-\eqref{eq_HHD_Poisson} is repeated many times, once for every iteration (or time step), it is critical that the discretization of such equations is stable, i.e., it does not introduce any spurious mode which grows indefinitely, thus compromizing the whole simulation. In order to assess the stability properties of the strategies presented in this work, it was decided to apply multiple times the projection scheme described in equations \eqref{eq_HHD_u}-\eqref{eq_HHD_Poisson} starting from a given vector field $\boldsymbol{w}$, defined on the complex-shaped domain depicted in Figure \ref{fig:geometry_HHD}. The boundary is characterized by locally convex and concave features with varying curvature, designed to trigger ill-conditioning issues described in section \ref{s:Neumann}. A Lagrange multiplier is employed to solve the linear systems arising from the RBF-FD discretization of equation \eqref{eq_HHD_Poisson} since Neumann BCs are prescribed on the whole boundary \cite{Zamolo_2022}.

The pseudocode for the stability test is presented in Algorithm \ref{alg_stab}. It is worth mentioning that any conservative discretization scheme, e.g., FVM, implicitly satisfies the previous stability test since the value of the discretized divergence operator vanishes after the first iteration. On the other hand, the basic formulation of the RBF-FD discretization, as presented in section \ref{ss:RBF_FD}, is not conservative, therefore the value of the discretized divergence operator never vanishes at the nodes. Algorithm \ref{alg_stab} might therefore diverge if the Poisson equation is solved with the basic RBF-FD method without a proper treatment of the Neumann BCs.

\begin{figure}[t!]
    \def\SpaceBelowText{.25em}
    \def\ImagesWidth{.32\textwidth}
    \setlength\fboxsep{0pt}
    \setlength\fboxrule{0.5pt}
    \centering

    \begin{subfigure}[b]{\ImagesWidth}
        \centering
        \hspace{1.25em} $P=2$\\[\SpaceBelowText]
        \includegraphics[width=\textwidth]{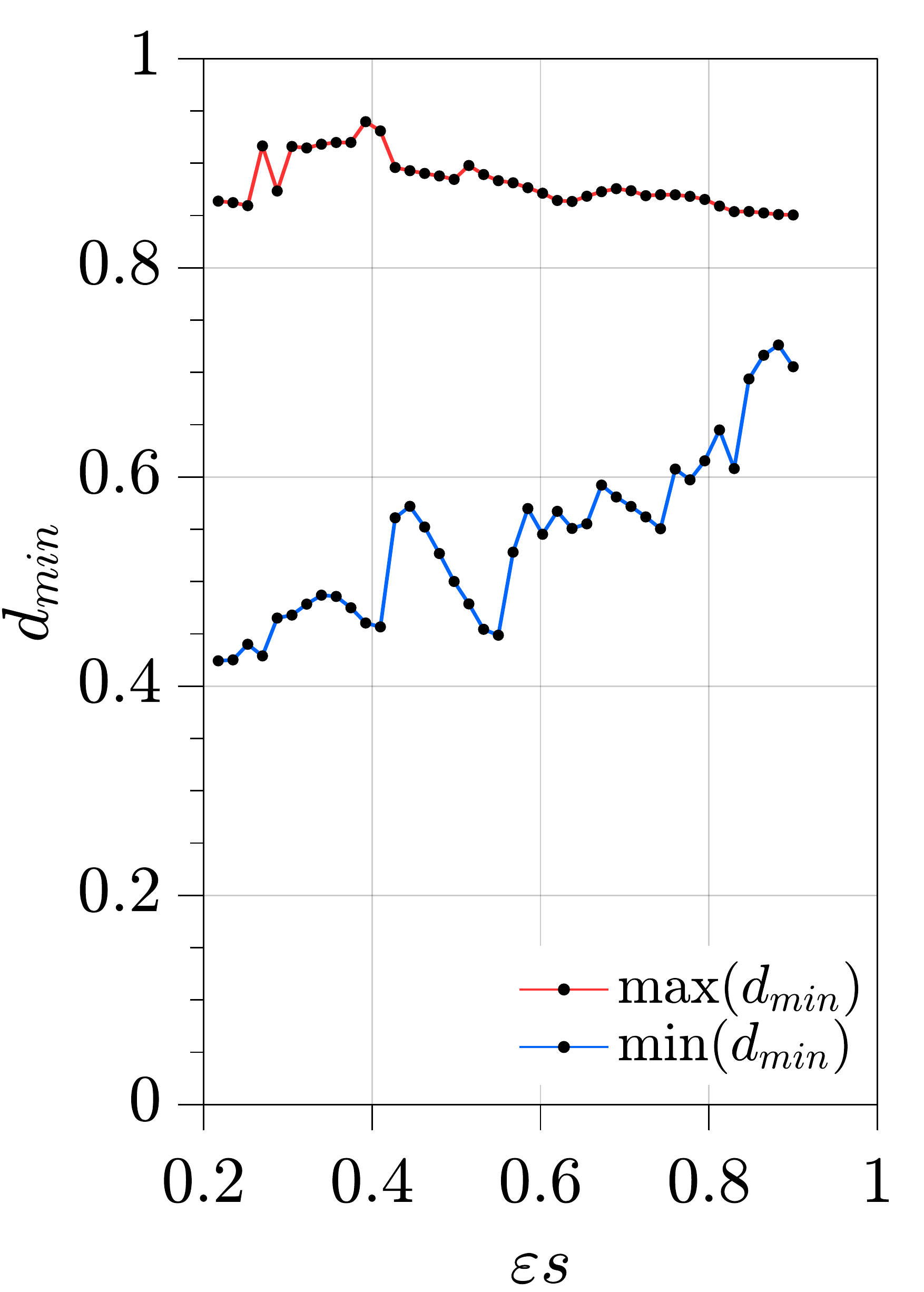}
    \end{subfigure}
    \hfill
    \begin{subfigure}[b]{\ImagesWidth}
        \centering
        \hspace{1.25em} $P=3$\\[\SpaceBelowText]
        \includegraphics[width=\textwidth]{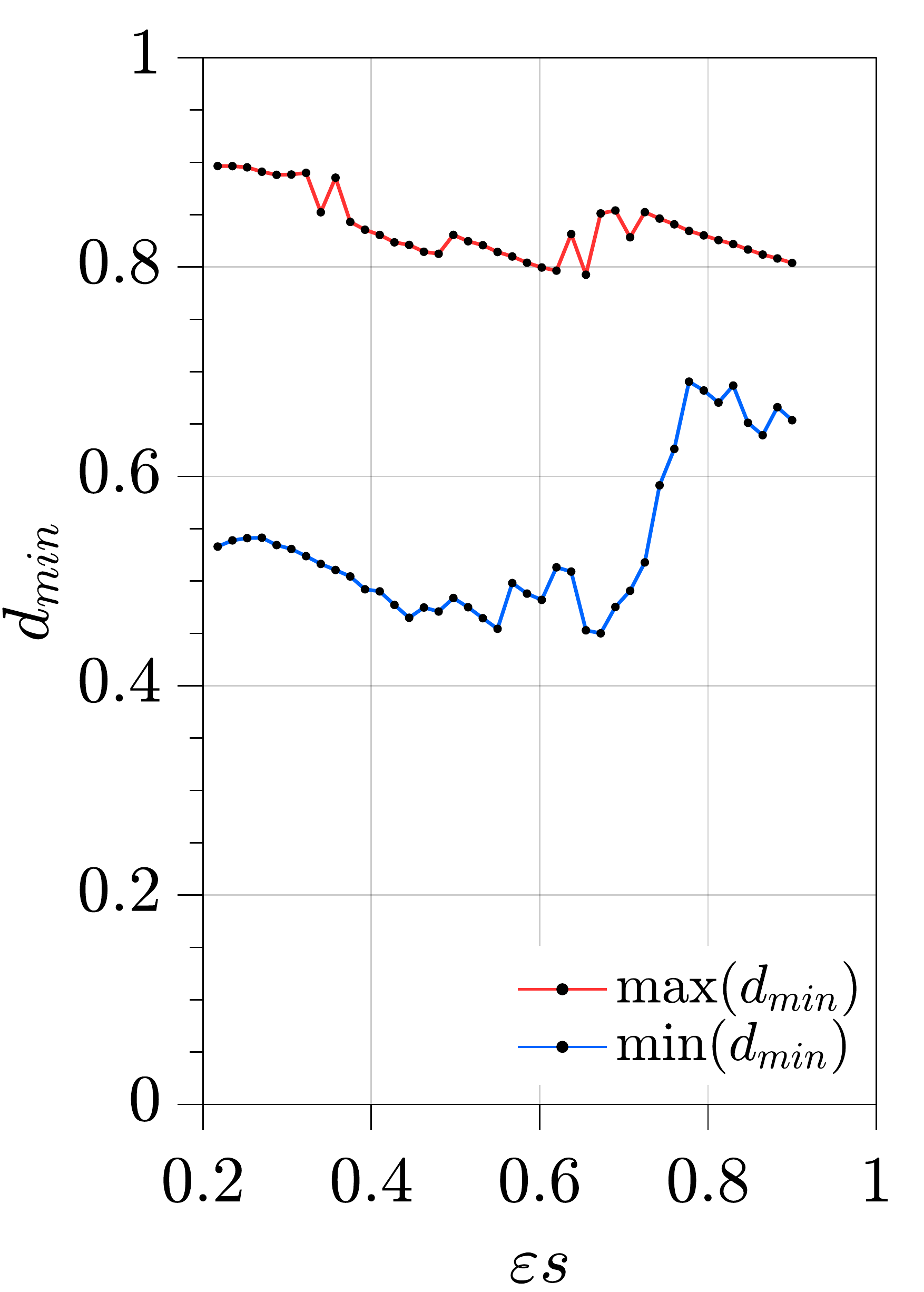}
    \end{subfigure}
    \hfill
    \begin{subfigure}[b]{\ImagesWidth}
        \centering
        \hspace{1.25em} $P=4$\\[\SpaceBelowText]
        \includegraphics[width=\textwidth]{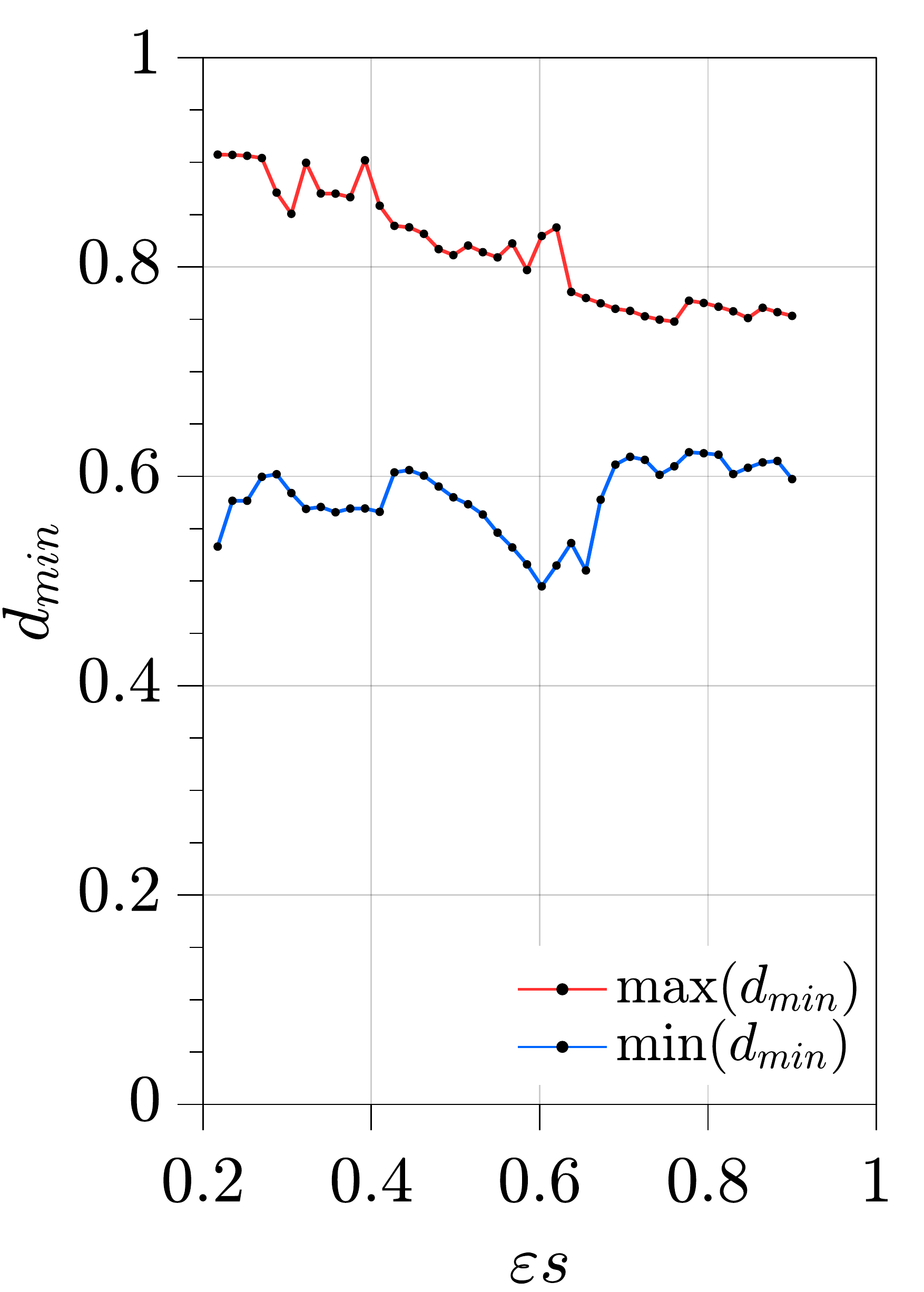}
    \end{subfigure}
    \caption{stability range for the $d_{min}$ parameter in the case of approach 1, i.e., boundary node selection based on optimal directions. Upper and lower bounds for $d_{min}$ values are displayed with red and blue curves, respectively.}
    \label{FIG_norm_sel_stab}
\end{figure}

\begin{algorithm}[b!]
	\caption{Pseudocode for the stability test of repeated HHD}
	\label{alg_stab}
	\begin{algorithmic}[1]
	    \State Initialize: $\boldsymbol{u}^{(1)}\gets\boldsymbol{w}$
		\For {$i=1,2,\ldots$}
				\State Solve $\nabla^2\phi = \nabla\cdot\boldsymbol{u}^{(i)}$ for $\phi$
				with homogeneous Neumann BCs
				\State Update: $\boldsymbol{u}^{(i+1)}\gets\boldsymbol{u}^{(i)}-\nabla\phi$
		\EndFor
	\end{algorithmic} 
\end{algorithm}

Since the equations involved in Algorithm \ref{alg_stab} are linear, the choice of the initial vector field $\boldsymbol{w}$ does not affect the stability of the iterative process. The employed initial vector field $\boldsymbol{w}$ is chosen to be an irrotational vortex for simplicity: the tangential velocity is $w_\theta=r^{-1}$ where $r$ is the distance from the origin.

The results for approach 1 presented in section \ref{sss:dot_product}, i.e., boundary node selection, are shown in Figure \ref{FIG_norm_sel_stab} in the case of MQ RBF with polynomial degrees $P=2,3,4$ and for different values of the shape parameter $\varepsilon$, such that $\varepsilon s\in[0.2,0.9]$. The total number of nodes inside the domain is $N_I\approx 15,000$. The number of internal nodes for each stencil is chosen to follow the rule $m_I = 2q$ \cite{bayona2017role}, see equation \eqref{eq:polyRBF}: $m_I=20$ for $P=3$ and $m_I=30$ for $P=4$. In the case $P=2$ it was decided to use slightly larger stencils with $m_I=15$ internal nodes, since this enhanced stability. In Figure \ref{FIG_norm_sel_stab} the values of $d_{min}$ allowing a stable computation were identified as those lying between the two curves.
As expected, too small values of $d_{min}$, e.g., $d_{min}<0.2$, lead to instability due to stencils with boundary nodes whose actual normal is too different from the corresponding optimal direction, resulting in ill-conditioned local interpolants. On the other hand, too large values of $d_{min}$, e.g., $d_{min}>0.95$, lead to instability due to stencils with too few boundary nodes enforcing the prescribed Neumann BCs \eqref{eq_HHD_Poisson_BC}. From the same figure it is possible to observe that the stability range is slightly reduced when $P$ is increased from $P=2$ to $P=4$, while a general and reasonably good choice is $d_{min}=0.7$.

\begin{figure}[t!]
    \def\SpaceBelowText{.5em}
    \def\ImagesWidth{.32\textwidth}
    \setlength\fboxsep{0pt}
    \setlength\fboxrule{0.5pt}
    \centering

    \begin{subfigure}[b]{\ImagesWidth}
        \centering
        \hspace{2em}No stabilization\\[\SpaceBelowText]
        \includegraphics[width=\textwidth]{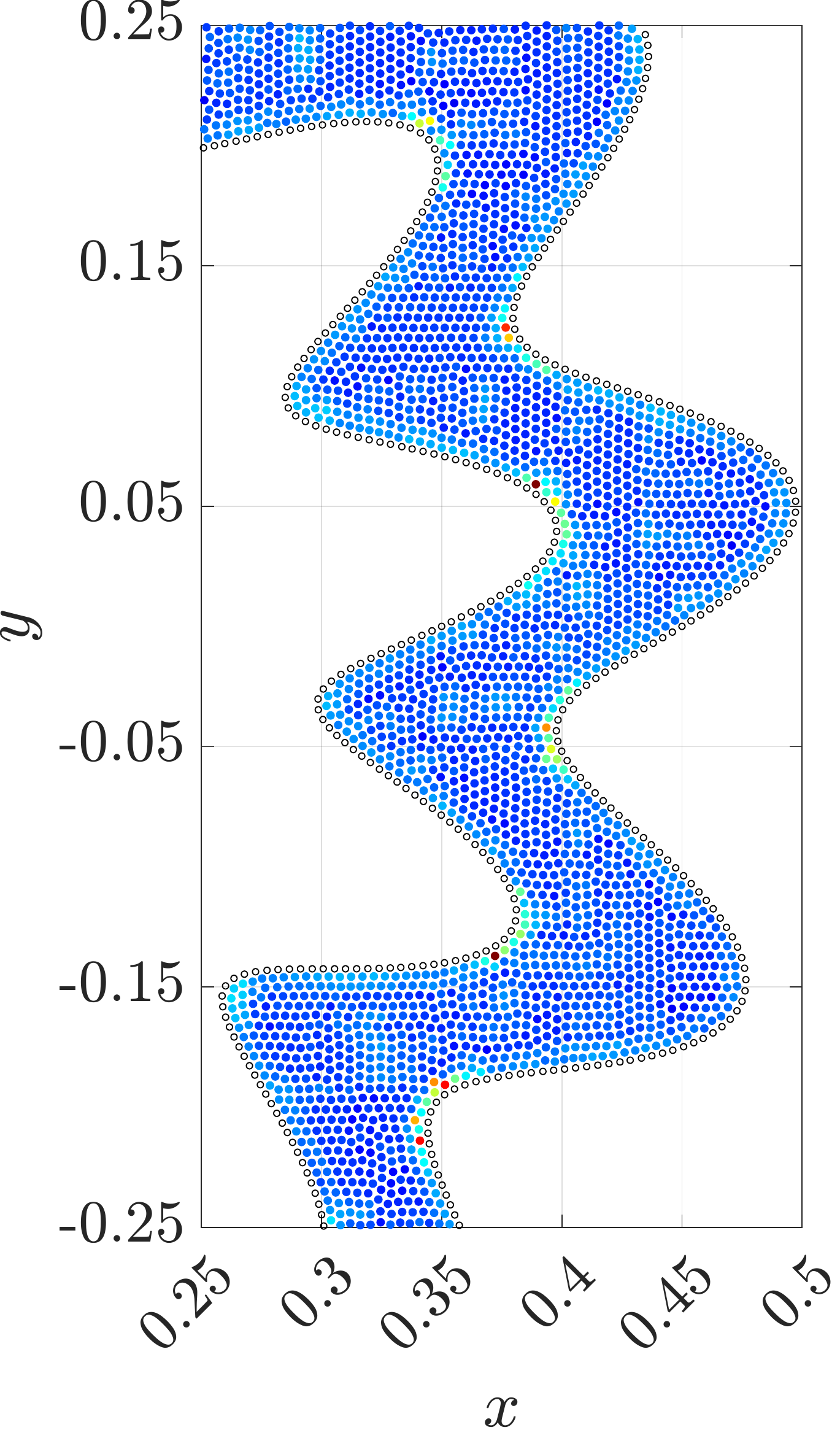}
    \end{subfigure}
    \hfill
    \begin{subfigure}[b]{\ImagesWidth}
        \centering
        \hspace{2em}Approach 1\\[\SpaceBelowText]
        \includegraphics[width=\textwidth]{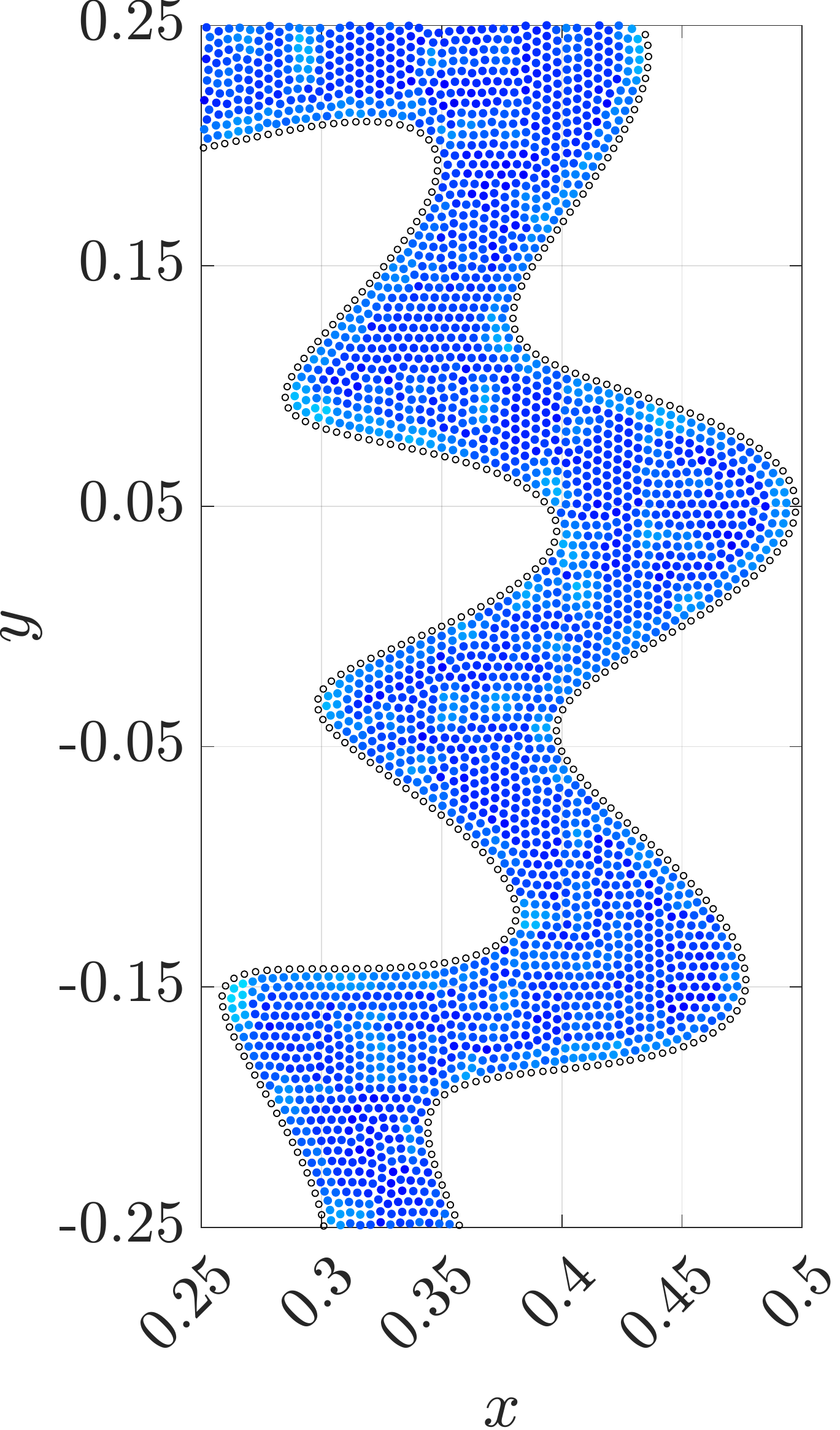}
    \end{subfigure}
    \hfill
    \begin{subfigure}[b]{\ImagesWidth}
        \centering
        \hspace{2em}Approach 2\\[\SpaceBelowText]
        \includegraphics[width=\textwidth]{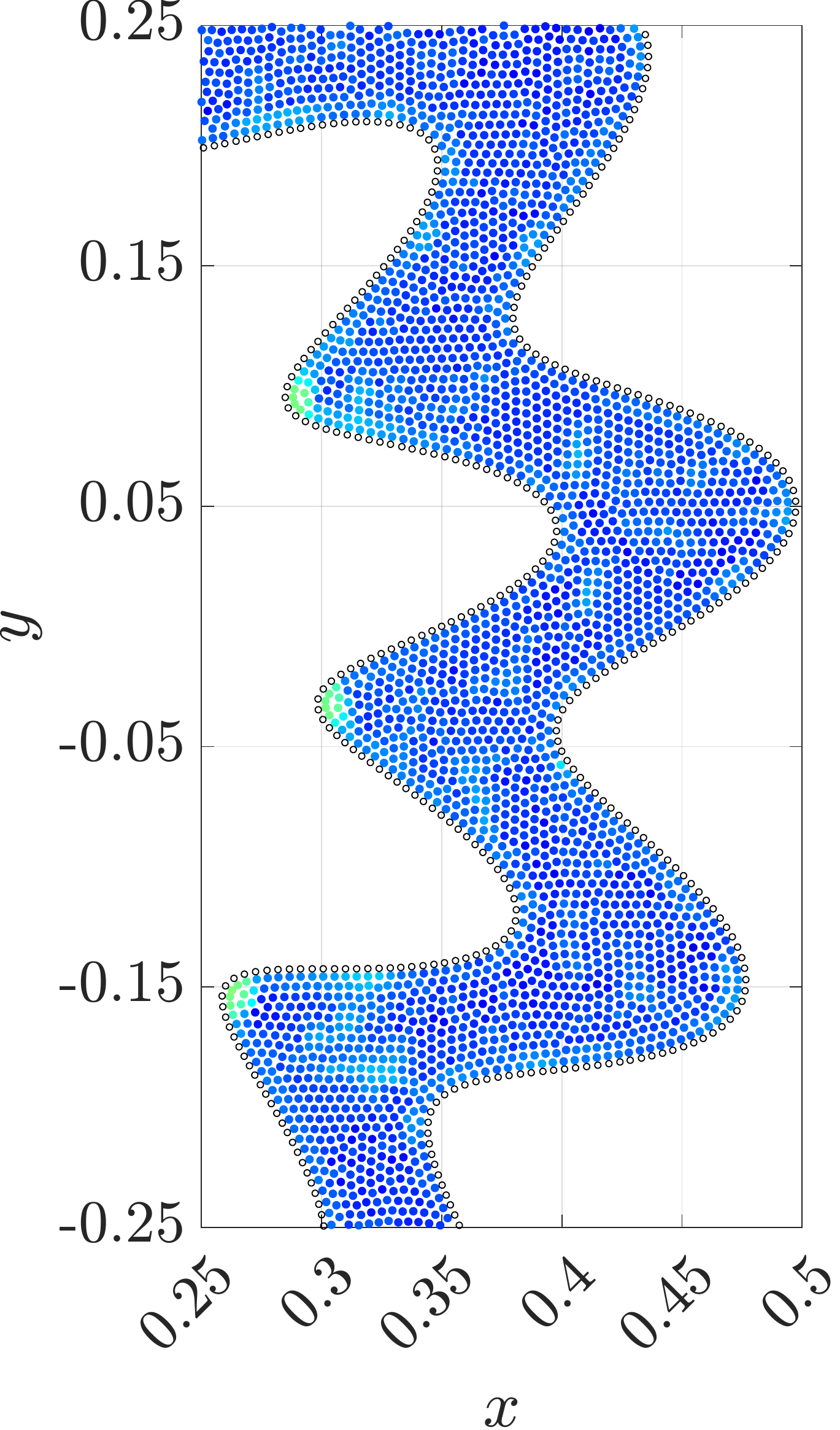}
    \end{subfigure}
    \vspace{1em}
    \begin{subfigure}[b]{\ImagesWidth}
        \centering
        \includegraphics[width=\textwidth]{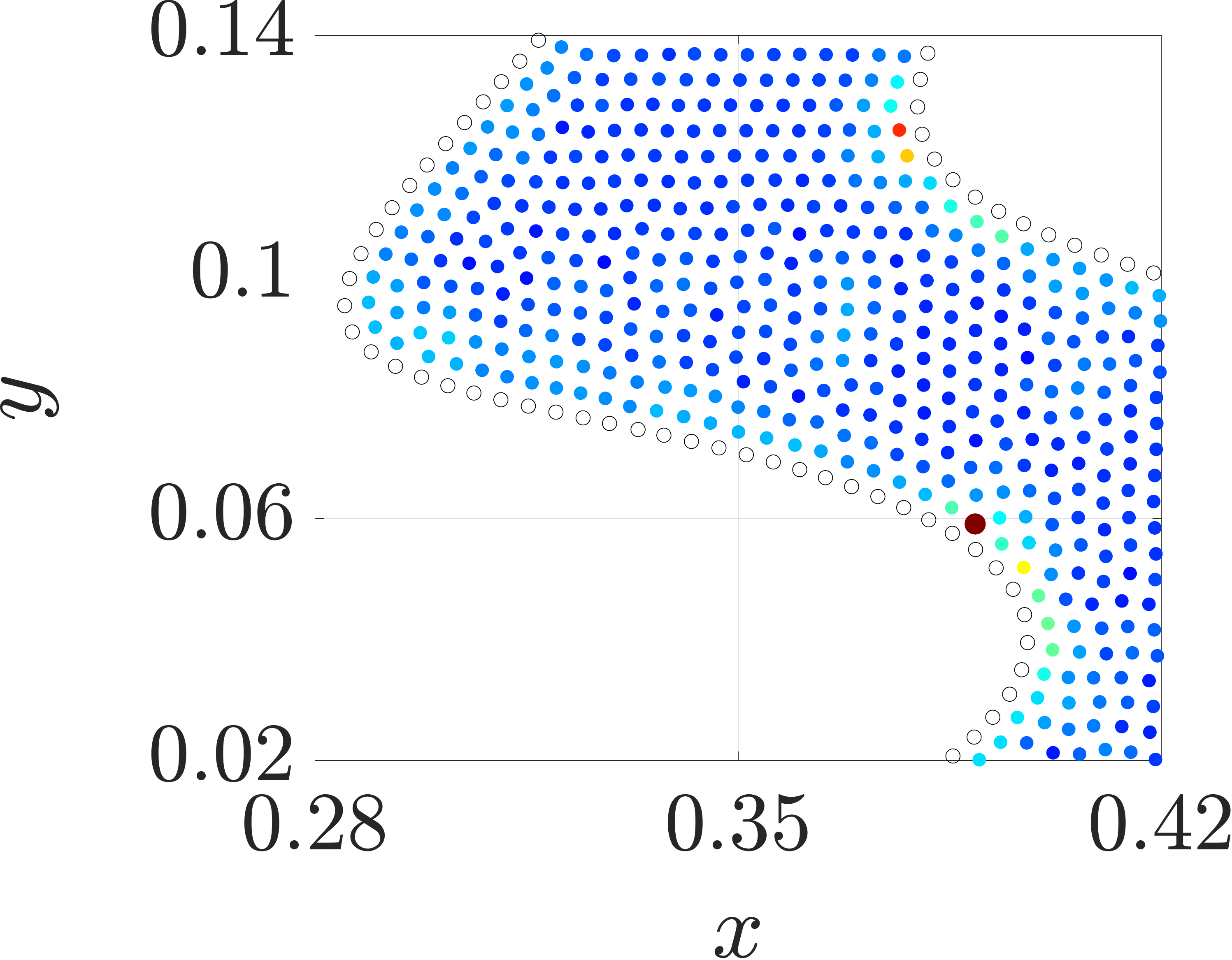}
    \end{subfigure}
    \hfill
    \begin{subfigure}[b]{\ImagesWidth}
        \centering
        \includegraphics[width=\textwidth]{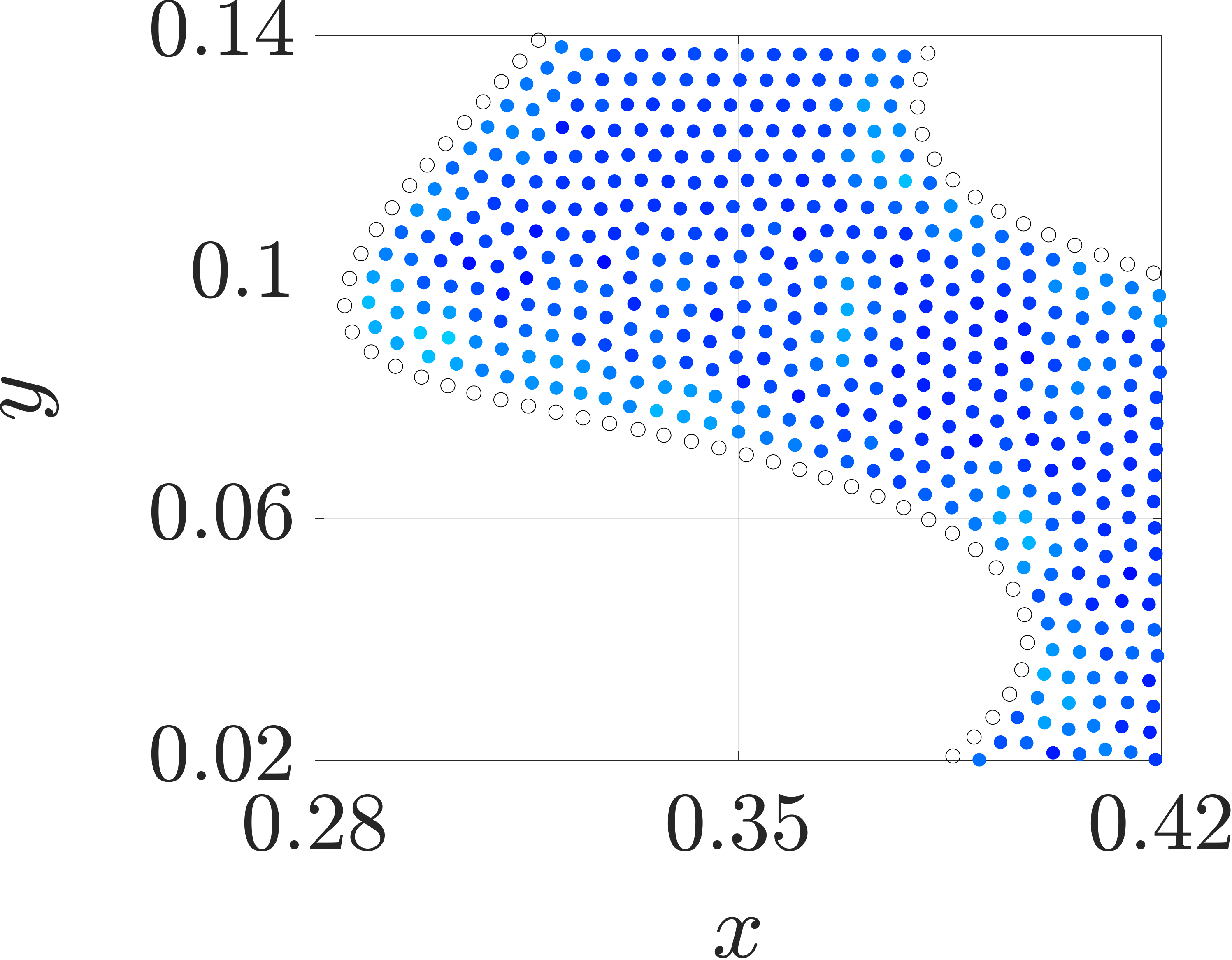}
    \end{subfigure}
    \hfill
    \begin{subfigure}[b]{\ImagesWidth}
        \centering
        \includegraphics[width=\textwidth]{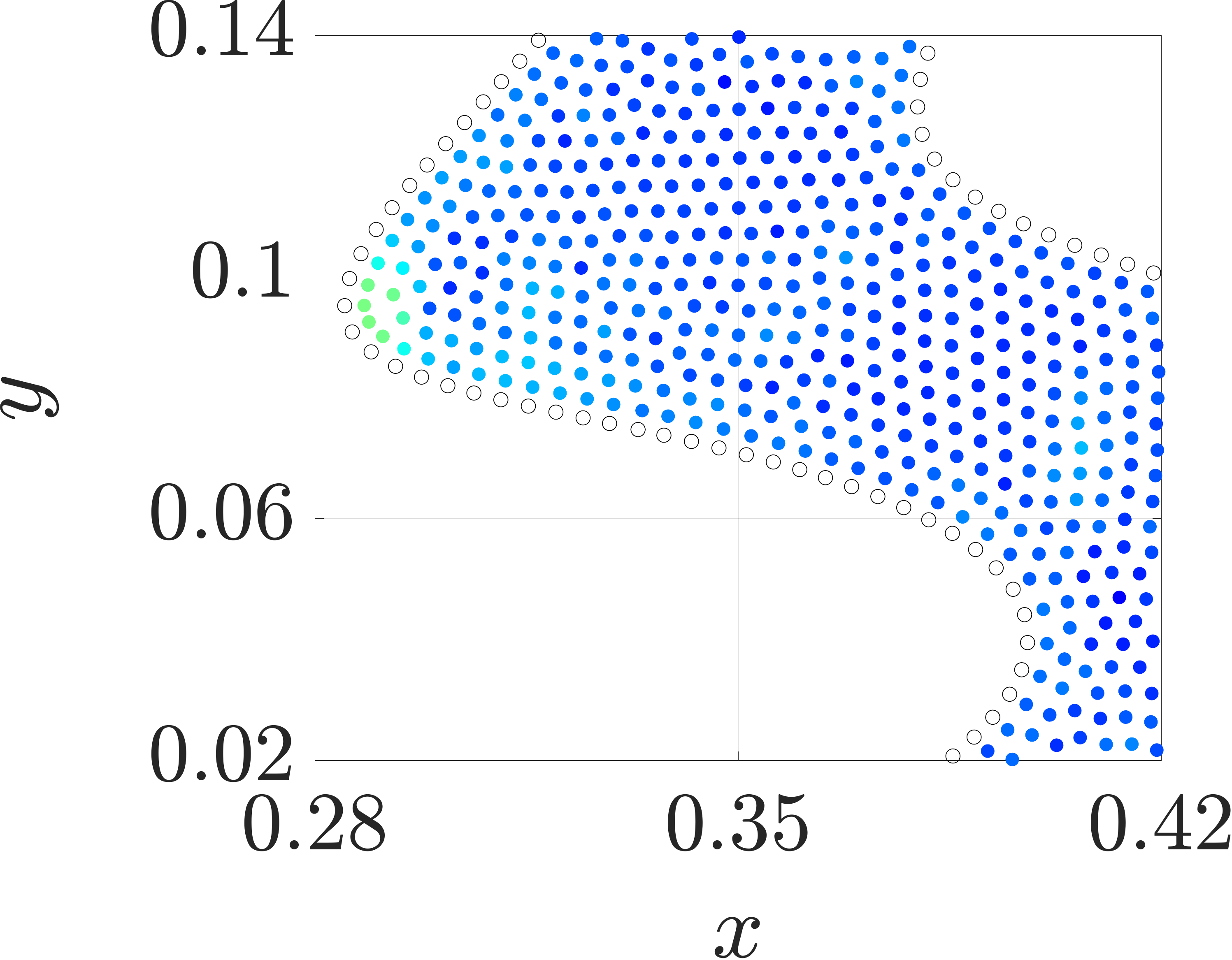}
    \end{subfigure}
    \\[1em] 
    \includegraphics[width=.45\textwidth]{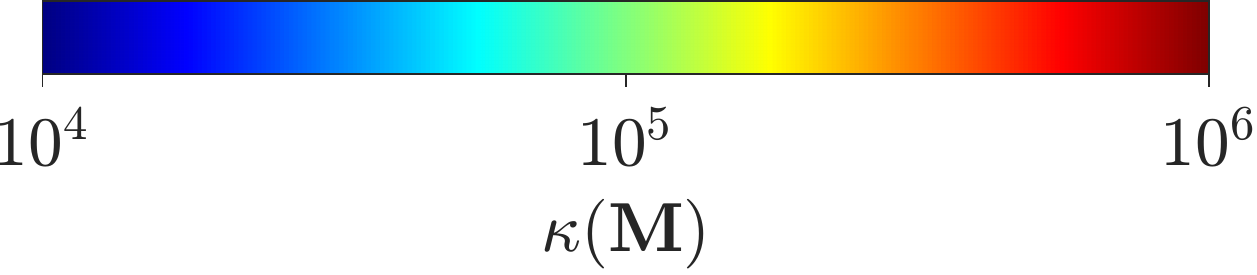}
    \caption{condition number of the local interpolation matrix. Top row, from left to right: no stabilization, boundary selection based on optimal directions, boundary projection. Bottom row: enlarged visualization of a portion of the plots of the top row.}
    \label{FIG_cond_numb}
\end{figure}

In the case of approach 2 presented in section \ref{ss:optimal_placement}, i.e., by employing projected boundary nodes, stability is always attained for $P=2,3,4$ in the range $\varepsilon s\in[0.2,0.9]$ previously considered, highlighting the remarkable stabilization effect of this approach. Nevertheless, stability issues are still encountered for shape factors outside this range.

Figure \ref{FIG_cond_numb} shows the effect of the proposed stabilization approaches on the condition number of the local interpolation matrix, in the case of MQ RBF with $\varepsilon s=0.5$ and $P=3$ (results with different shape parameters $\varepsilon$ and polynomial degrees $P$ are qualitatively similar). As expected from the theoretical observations presented in section \ref{s:Neumann}, stencils close to locally concave boundaries can be very ill-conditioned if no stabilization technique is employed (red nodes in Figure \ref{FIG_cond_numb}, left), leading to large errors and stability issues.

By using the stabilization approach 1 with $d_{min}=0.7$, this issue is properly addressed and ill-conditioned stencil are no longer present (Figure \ref{FIG_cond_numb}, center). The same result is attained by employing the stabilization approach 2 (Figure \ref{FIG_cond_numb}, right). In the latter case we observe an increase in the condition number for the stencils close to locally convex boundaries: this is due to the employed projection technique that is performed from the boundary inward. Therefore, internal nodes projected from locally convex boundaries are very close to each other, leading to an increased condition number. This geometrical issue is common to any simple geometrical projection technique, regardless of the direction and starting locations employed for the projection. With reference to arbitrary geometries, such problem can only be addressed by taking into account the need for projected boundary nodes already at the node generation phase, thus requiring a proper modification of the node generation algorithm.

\subsection{Accuracy}

Besides the stability considerations presented in the previous section, we also provide some insight about the influence of the proposed stabilization approaches on the accuracy of the RBF-FD discretization. In order to do so, we consider the following Poisson equation:
\begin{equation}
    \nabla^2 u = f
    \label{eq_Poisson_complex}
\end{equation}
defined over the previously employed domain depicted in Figure \ref{fig:geometry_HHD}. Neumann BCs are prescribed on the whole boundary:
\begin{equation}
    \frac{\partial u}{\partial \boldsymbol{n}} = g
    \label{eq_Poisson_complex_BC}
\end{equation}
and the analytic solution is chosen to be $u=r^{-1}$, where $r$ is the distance from the origin. The node distributions and the parameters for the RBF-FD discretization are the same as those used in the previous section \ref{ss_Appl_Stability}, and a Lagrange multiplier is again employed to solve the corresponding discretized linear system.

\begin{figure}[t!]
    \def\SpaceBelowText{.25em}
    \def\ImagesWidth{.49\textwidth}
    \centering
    \begin{subfigure}[b]{\ImagesWidth}
        \centering
        \includegraphics[width=\textwidth]{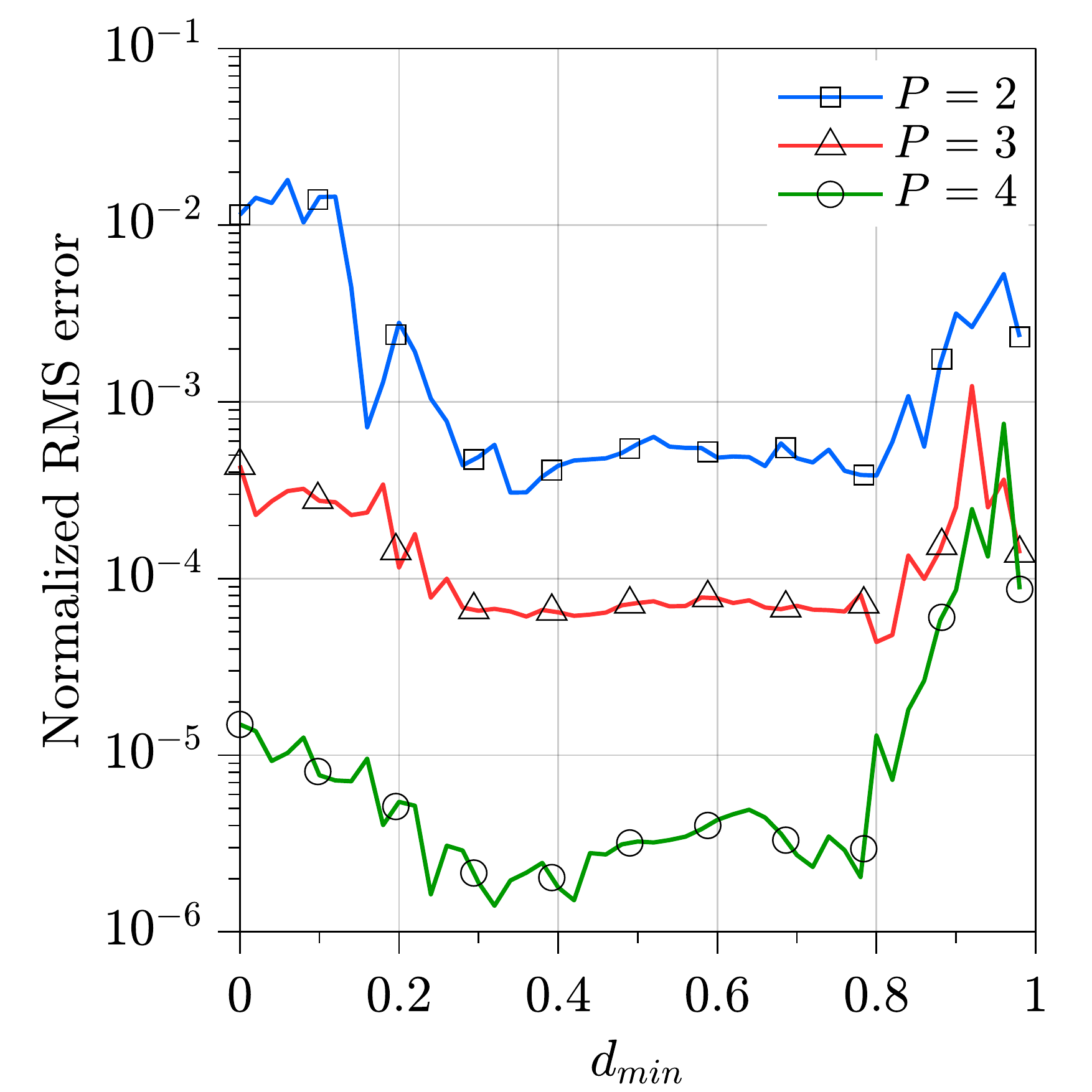}
    \end{subfigure}
    \hfill
    \begin{subfigure}[b]{\ImagesWidth}
        \centering
        \includegraphics[width=\textwidth]{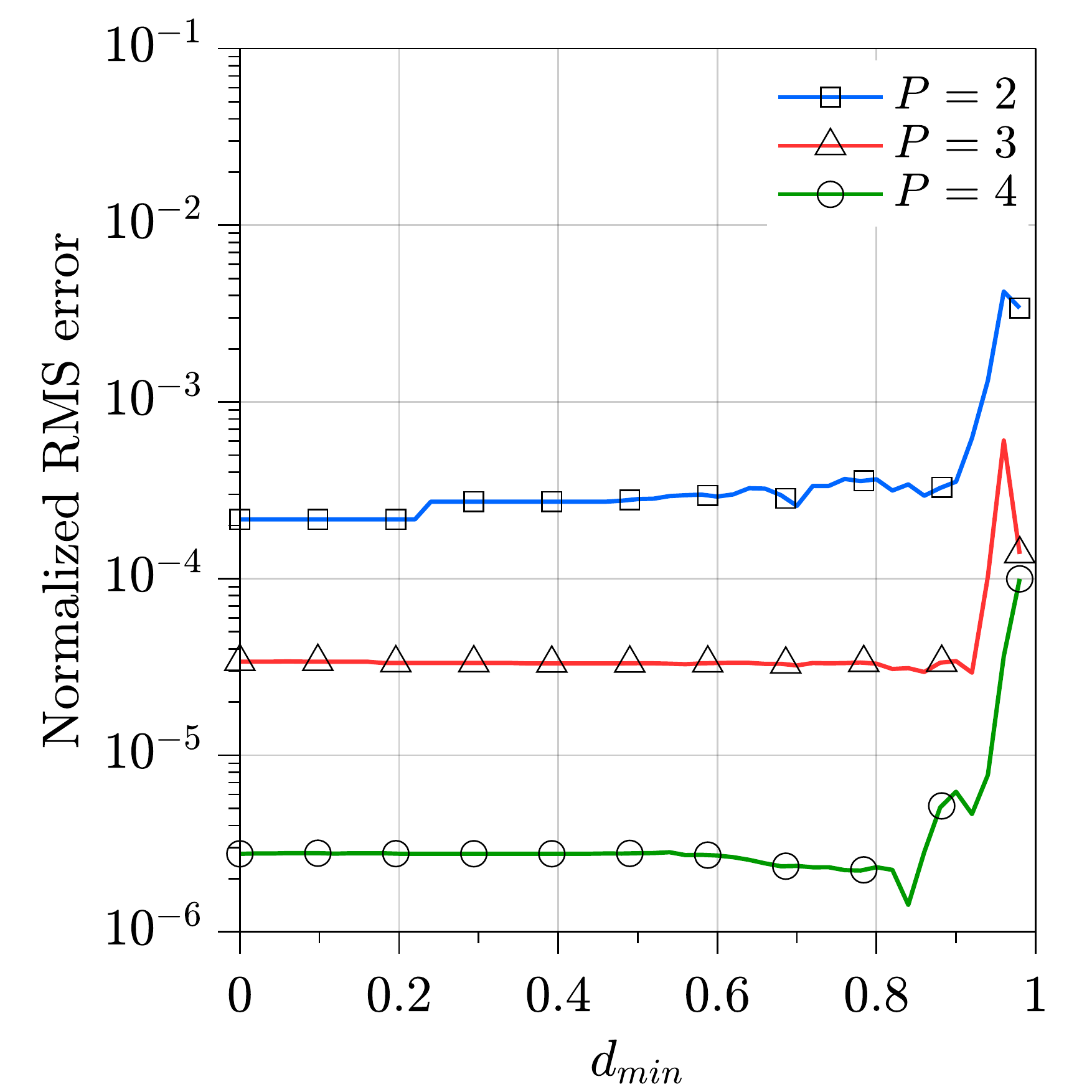}
    \end{subfigure}
    \caption{Errors for the Poisson problem \eqref{eq_Poisson_complex} when employing stabilization approach 1: initial node distribution (left), projected boundary nodes (right).}
    \label{FIG_NRMSE_dmin}
\end{figure}

The normalized RMS errors are shown in Figure \ref{FIG_NRMSE_dmin} when approach 1 is applied to both node distributions without and with projected boundary nodes. The latter case is therefore a combination of both the proposed stabilization techniques. In the former case, and for all polynomial degrees $P=2,3,4$, the error is lower when $d_{min}$ ranges from 0.4 to 0.8. As expected, the error grows when $d_{min}$ is decreased below 0.4 because of the increasing number of ill-conditioned stencils near the boundary. The limit case $d_{min}=0$ corresponds indeed to no stabilization action at all. On the other hand, the error also grows when $d_{min}$ is increased beyond 0.8 because of the increasing number of boundary stencils with an insufficient number of boundary nodes where BCs are enforced.

When both the stabilization approaches are employed, Figure \ref{FIG_NRMSE_dmin} right, the increase in the error for $d_{min}<0.4$ is no longer present since the ill-conditioning issues are already addressed by the projected nodes on the boundary, while the increase for $d_{min}>0.8$ is obviously still present, i.e., when too few boundary nodes are included in the boundary stencils. In this case we also observe that the boundary node selection based on the optimal directions does not affect the accuracy of the discretization at all if $d_{min}<0.8$, as expected.

\section{Conclusions}
In this paper the problem of ill-conditioning and the subsequent instabilities affecting the application of RBF-FD method in presence of Neumann boundary condition have been analyzed both theoretically and numerically. The theoretical motivations for the onset of such issues are derived by highlighting the dependency of the determinant of the local interpolation matrix upon the boundary normals. Qualitative investigations are also carried out numerically by studying a reference stencil and looking for correlations between its geometry and the properties of the associated interpolation matrix.

Once the problem had been described with sufficient detail, two approaches are proposed, namely approach 1, i.e., the boundary node selection based on the optimal normals, and approach 2, i.e., the optimal placement of boundary nodes. Many evidences indicate that both can be applied to the solution of numerical problems of practical interest.

In order to provide a better demonstration of the onset of instabilities due to ill-conditioning of the interpolation matrix, the stabilization of the Helmholtz-Hodge decomposition is discussed in the last section. This demonstration helped to put into context the role of the approaches in the stabilization of the solution process of any boundary value problem involving Neumann boundary conditions.

Summing up, the main pros and cons of the two stabilization strategies are the following. Approach 1 is the most robust and it does not require any control over the node placement. However, the price to pay for its application is a higher computational effort, required to estimate the optimal directions for all the boundary stencils. We remark that this cost can be significantly reduced if the optimal directions are approximated with lower accuracy. Future developments might significantly improve the computational efficiency of approach 1 through the application of machine learning techniques, for example. Furthermore, the mathematical formulation behind this approach leaves great freedom for further developments.
Approach 2, on the other side, does not involve any additional computation and is even more effective at ensuring the stabilization of the RBF-FD method. Unfortunately, it requires the implementation of a more sophisticated node generation algorithm, indeed, more complex geometries might present high curvatures at the boundary thus making it very difficult to develop a generic node projection technique.
At last we remark that both approaches can be applied without any modification to the 3D cases, where all of the considerations made in this paper remain valid.

The application of the stabilization techniques described above does not require the user to intervene during the simulation process and does not require any parameter tuning. Even the $d_{min}$ parameter can be set once for all to be $d_{min} = 0.7$ without any contraindication. This fact is especially important since the major advantage of the RBF-FD method lies in its ability to free the user from all of the considerations involved in the mesh generation, which in turn require a specific expertise. Indeed, for this method to be competitive with traditional ones, it must succeed in making the whole simulation process easier to handle, without sacrificing the reliability of the attained solutions.

Furthermore, the presented approaches can be suitably combined to the appropriate extent to provide even more robust RBF-FD meshless discretizations of boundary value problems.

In conclusion, we believe that the approaches presented in this work are successful not only at solving the stability issues of the RBF-FD method, but also at preserving its potential advantages.

\appendix
\section{\texorpdfstring{$d$}{d}-matrices}\label{App:s1}

\begin{definition}[$d$-matrix]
An $m\times n$ $d$-matrix is defined as the following $m\times n$ matrix $\dmat{A}{a}$ of $d$-dimensional vectors $\boldsymbol{a}_{ij}\in\mathbb{R}^{d}$:
\begin{equation}
    \mathcal{A} =
    \begin{bmatrix}
        \boldsymbol{a}_{11} \mysep \cdots \mysep \boldsymbol{a}_{1n} \\
        \vdots              \mysep \ddots \mysep \vdots \\
        \boldsymbol{a}_{m1} \mysep \cdots \mysep \boldsymbol{a}_{mn} \\
    \end{bmatrix}
    \label{App1_matrix_of_v}
\end{equation}
and the set of all $m\times n$ $d$-matrices will be denoted by $\Rd{m}{n}$. For simplicity $\Rd{m}{1}$ will be called $\Rdvec{m}$.
\end{definition}

The usual operations are defined as follows:
\begin{itemize}
    \setlength\itemsep{-.5em}
    \item product of $\dmat{A}{a}\in\Rd{m}{n}$ with a scalar $\lambda\in\mathbb{R}$:
    \begin{equation}
        \lambda\mathcal{A}=(\lambda\boldsymbol{a}_{ij})\in\Rd{m}{n}
        \label{App1_dmat_prod}
    \end{equation}
    \item sum of $\dmat{A}{a}\in\Rd{m}{n}$ and $\dmat{B}{b}\in\Rd{m}{n}$:
    \begin{equation}
        \mathcal{A}+\mathcal{B}=(\boldsymbol{a}_{ij}+\boldsymbol{b}_{ij})\in\Rd{m}{n}
        \label{App1_dmat_sum}
    \end{equation}
    \item product of $\dmat{A}{a}\in\Rd{m}{n}$ with $\mathbf{Q}=(q_{ij})\in\mathbb{R}^{n\times p}$:
    \begin{equation}
        \mathcal{A}\mathbf{Q}=\left(\sum_{k=1}^{n}\boldsymbol{a}_{ik}q_{kj}\right)\in\Rd{m}{p}
        \label{App1_dmat_matprod}
    \end{equation}
\end{itemize}

For $\dmat{A}{a}\in\Rd{m}{n}$ and $\dvec{V}{v}\in\Rdvec{m}$, by using the dot product $\cdot$ in $\mathbb{R}^d$ we can define the following operator H:
\begin{equation}
    H(\mathcal{A},\mathcal{V})=(\boldsymbol{a}_{ij}\cdot\boldsymbol{v}_{i})\in\mathbb{R}^{m\times n}
    \label{App1_opH}
\end{equation}

Given $\mathcal{A},\mathcal{B}\in\Rd{m}{n}$, $\dvec{V}{v}\in\Rdvec{m}$ and $\mathbf{Q}\in\mathbb{R}^{n\times p}$, the following equalities hold:
\begin{linenomath*}
\begin{align}
    H(\mathcal{A},\mathcal{V})+H(\mathcal{B},\mathcal{V})& =
    H(\mathcal{A}+\mathcal{B},\mathcal{V}) \label{App1_prop_sum}\\
    H(\mathcal{A},\mathcal{V}) \mathbf{Q} & =
    H(\mathcal{A}\mathbf{Q},\mathcal{V}) \label{App1_prop_prod} \\
    \frac{\partial}{\partial v_{i,\eta}}\det\!\big(H(\mathcal{A},\mathcal{V})\big)&=\det\!\big(H(\mathcal{A},\mathcal{V}_{i,\eta})\big) \label{App1_prop_det}
\end{align}
\end{linenomath*}
for $i=1,\ldots,m$ and $\eta=1,\ldots,d$. By considering the canonical basis $\{\boldsymbol{e}_\eta\}_1^d$ of $\mathbb{R}^{d}$ for simplicity, in the last equation $v_{i,\eta}$ is the $\eta^{th}$ component of the $i^{th}$ vector $\boldsymbol{v}_i$ in $\mathcal{V}$, while $\mathcal{V}_{i,\eta}$ is obtained from $\mathcal{V}$ by replacing its $i^{th}$ vector $\boldsymbol{v}_i$ with $\boldsymbol{e}_{\eta}$. Equations \eqref{App1_prop_sum}-\eqref{App1_prop_prod} follow from the simple distributive property of the dot product, while equation \eqref{App1_prop_det} follows from the multilinearity of the determinant and from the definition of the operator $H$ in equation \eqref{App1_opH}. 
\section{Optimal directions}\label{App:s2}
\subsection{Generic properties}
By defining the Lagrangian function $\mathcal{L}$ with multipliers $\boldsymbol{\nu}=\{\nu_1,\ldots,\nu_{m_B}\}$ as follows:
\begin{equation}
    \mathcal{L}(\mathcal{N},\boldsymbol{\nu})=\det\!\big(H(\mathcal{G}_{BB},\mathcal{N})\big)-\sum_{i=1}^{m_B}\nu_i\left(\|\bar{\boldsymbol{n}}_i\|_2^2-1\right)
\end{equation}
the solution of problem \eqref{eq:constr_opt_mB_nodes} is given by the proper solution of the following system of equations:
\begin{equation}
\begin{cases}
\displaystyle
\frac{\partial}{\partial \hat{n}_{i,\eta}} \mathcal{L}(\hat{\mathcal{N}},\boldsymbol{\nu}) =
\frac{\partial}{\partial \hat{n}_{i,\eta}} \det\!\big(H(\mathcal{G}_{BB},\hat{\mathcal{N}})\big)-2\nu_i \hat{n}_{i,\eta}=0\\[1em]
\|\hat{\boldsymbol{n}}_i\|_2^2=1
\end{cases}
    \label{eq:LagrangeSyst_expl}
\end{equation}
for $i=1,\ldots,m_B$ and $\eta=1,\ldots,d$. By looking at the first equation of system \eqref{eq:LagrangeSyst_expl}, we see that the Lagrange multipliers $\nu_i$ act as multiplicative factors allowing the constraints in problem \eqref{eq:constr_opt_mB_nodes} to be fulfilled. Therefore, the solution $\hat{\mathcal{N}}=\{\hat{\boldsymbol{n}}_1,\ldots,\hat{\boldsymbol{n}}_{m_B}\}\in\Rdvec{m_B}$ can be explicitly written through the $m_B$ auxiliary vectors $\boldsymbol{t}_i=\{t_{i,1},\ldots,t_{i,d}\}\in\mathbb{R}^d$ as follows:
\begin{linenomath*}
\begin{align}
    t_{i,\eta}&=\frac{\partial}{\partial \hat{n}_{i,\eta}}\det\!\big(H(\mathcal{G}_{BB},\hat{\mathcal{N}})\big)=\det\!\big(H(\mathcal{G}_{BB},\hat{\mathcal{N}}_{i,\eta})\big)
    \label{eq:ti_opt}\\
    \hat{\boldsymbol{n}}_i&=\pm \frac{\boldsymbol{t}_i}{\|\boldsymbol{t}_i\|_2}
    \label{eq:constr_opt_mB_nodes_solution}
\end{align}
\end{linenomath*}
for $i=1,\ldots,m_B$ and $\eta=1,\ldots,d$. Property \eqref{App1_prop_det} is used in equation \eqref{eq:ti_opt}, where $\hat{\mathcal{N}}_{i,\eta}$ is obtained from $\hat{\mathcal{N}}$ as explained in appendix \ref{App:s1}. We note that the term $\hat{\mathcal{N}}_{i,\eta}$ depends on $\hat{\boldsymbol{n}}_{j\neq i}$, i.e., it does not depend on $\hat{\boldsymbol{n}}_i$ itself. Equation \eqref{eq:constr_opt_mB_nodes_solution} then acts as a normalization of $\boldsymbol{t}_i$ through the Lagrange multiplier $\nu_i$ as previously noted, i.e., $2\nu_i=\pm\|\boldsymbol{t}_i\|_2$. Equations \eqref{eq:ti_opt} and \eqref{eq:constr_opt_mB_nodes_solution} simply state that each solution vector $\hat{\boldsymbol{n}}_i$ is parallel to the gradient of the target function with respect to that normal, in accordance with the principles of the constrained optimization. 

From equation \eqref{eq:Schur_Multiple_node_cardinal_det_Gbb} we see that $\det(\mathbf{S}_{BB})$ is an odd function with respect to each normal $\bar{\boldsymbol{n}}_i$, i.e., the reversal of any $\bar{\boldsymbol{n}}_i$ results in a change of sign of $\det(\mathbf{S}_{BB})$ because of the properties of the determinant. Therefore the solution $\hat{\mathcal{N}}$ of the maximization problem \eqref{eq:constr_opt_mB_nodes} is also a solution when an even number of normals are reversed. Furthermore, analogously to the case with one boundary node, it is also a solution of the corresponding minimization problem when an odd number of normals are reversed, in accordance with the previous interpretation of the Lagrange multipliers. In other words, the target function in equation \eqref{eq:constr_opt_mB_nodes} has $2^{m_B}$ local extrema which are equivalent, hence the symbol $\pm$ in equation \eqref{eq:constr_opt_mB_nodes_solution} for each of the $m_B$ normals. This is obviously expected: the reversal in the directions of the normals in the Neumann BCs \eqref{eq:NeumannBC} has no effect on the properties of the interpolant, indeed we are interested in the normals that maximize the magnitude of the target function $\det(\mathbf{S}_{BB})$ regardless of its sign. When $m_B>1$ the target function has saddle points which therefore are also solutions of equations \eqref{eq:ti_opt}-\eqref{eq:constr_opt_mB_nodes_solution}, and it can also have multiple distinct local extrema, i.e., considering the reversal of any normal as an equivalent solution, as shown as follows.

Let us consider the case $m_B=2$ as an example. The target function $\det(\mathbf{S}_{BB})$ defined by equation \eqref{eq:Schur_Multiple_node_cardinal_det_Gbb} can be expressed as follows:
\begin{equation}
\begin{split}
    \det(\mathbf{S}_{BB}) &=
    \begin{vmatrix}
\boldsymbol{g}_{1,1}\cdot\bar{\boldsymbol{n}}_{1} & \boldsymbol{g}_{1,2}\cdot\bar{\boldsymbol{n}}_{1} \\
\boldsymbol{g}_{2,1}\cdot\bar{\boldsymbol{n}}_{2} & \boldsymbol{g}_{2,2}\cdot\bar{\boldsymbol{n}}_{2}
    \end{vmatrix}=\\
    &=(\boldsymbol{g}_{1,1}\cdot\bar{\boldsymbol{n}}_{1})(\boldsymbol{g}_{2,2}\cdot\bar{\boldsymbol{n}}_{2}) -(\boldsymbol{g}_{1,2}\cdot\bar{\boldsymbol{n}}_{1})(\boldsymbol{g}_{2,1}\cdot\bar{\boldsymbol{n}}_{2})=\bar{\boldsymbol{n}}_{1}\cdot\mathbf{G}\bar{\boldsymbol{n}}_{2}
    \end{split}
    \label{eq:det_optimal_mB_2}
\end{equation}
where:
\begin{equation}
\mathbf{G}=\boldsymbol{g}_{1,1}\otimes \boldsymbol{g}_{2,2}-
           \boldsymbol{g}_{1,2}\otimes \boldsymbol{g}_{2,1}
    \label{eq:det_optimal_mB_2_tens}
\end{equation}
is thus a $d\times d$ matrix having exactly $d-2$ null eigenvalues in the general case where $\boldsymbol{g}_{1,1}$ and $\boldsymbol{g}_{2,2}$ are not parallel to $\boldsymbol{g}_{1,2}$ and $\boldsymbol{g}_{2,1}$, respectively, and all $\boldsymbol{g}_{i,j}$ are non-zero vectors. The explicit solution expressed by equations \eqref{eq:ti_opt}-\eqref{eq:constr_opt_mB_nodes_solution} can therefore be simplified as follows:
\begin{equation}
\begin{split}
    \hat{\boldsymbol{n}}_1 &=\pm\frac{\mathbf{G}\hat{\boldsymbol{n}}_2}{\|\mathbf{G}\hat{\boldsymbol{n}}_2\|_2}\\
    \hat{\boldsymbol{n}}_2 &=\pm\frac{\mathbf{G}^T\hat{\boldsymbol{n}}_1}{\|\mathbf{G}^T\hat{\boldsymbol{n}}_1\|_2}
    \label{eq:explicit_sol_mB_2}
\end{split}
\end{equation}
which can be combined to yield the following decoupled equations:
\begin{equation}
\begin{split}
    \hat{\boldsymbol{n}}_1 &=\frac{\mathbf{Q}_1\hat{\boldsymbol{n}}_1}{\|\mathbf{Q}_1\hat{\boldsymbol{n}}_1\|_2}\\
    \hat{\boldsymbol{n}}_2 &=\frac{\mathbf{Q}_2\hat{\boldsymbol{n}}_2}{\|\mathbf{Q}_2\hat{\boldsymbol{n}}_2\|_2}
    \label{eq:explicit_sol_mB_2_dec}
\end{split}
\end{equation}
where $\mathbf{Q}_1=\mathbf{G}\mathbf{G}^T$ and $\mathbf{Q}_2=\mathbf{G}^T\mathbf{G}$. Equations \eqref{eq:explicit_sol_mB_2_dec} represent therefore two classic eigenvalue problems with matrix $\mathbf{Q}_1$ for $\hat{\boldsymbol{n}}_1$ and matrix $\mathbf{Q}_2$ for $\hat{\boldsymbol{n}}_2$. $\mathbf{Q}_1$ and $\mathbf{Q}_2$ have the same real non-negative eigenvalues $\{\lambda_i\}_1^d$ since these matrices are symmetric and because their definitions imply $\mathbf{G}^T\mathbf{Q}_1=\mathbf{Q}_2\mathbf{G}^T$ and $\mathbf{G}\mathbf{Q}_2=\mathbf{Q}_1\mathbf{G}$. 

In $d=2$ dimensions the eigenvalues are $0<\lambda_1<\lambda_2$ in the general case, therefore the possible distinct solutions of equations \eqref{eq:explicit_sol_mB_2} are four: two for each of the two normals. It can be shown that the solutions corresponding to complementary eigenvalues, i.e., $(\lambda_1,\lambda_2)$ or $(\lambda_2,\lambda_1)$ for $\hat{\boldsymbol{n}}_1$ and $\hat{\boldsymbol{n}}_2$, although solutions of equations \eqref{eq:explicit_sol_mB_2_dec}, are not solutions of equations \eqref{eq:explicit_sol_mB_2}. Therefore the two remaining solutions are those associated to the same eigenvalue for both normals, i.e., $(\lambda_1,\lambda_1)$ or $(\lambda_2,\lambda_2)$. By combining the last expression of equation \eqref{eq:det_optimal_mB_2} with the first equation \eqref{eq:explicit_sol_mB_2} with positive sign, we obtain:
\begin{equation}
\det(\mathbf{S}_{BB})\Big|_{\substack{
    \bar{\boldsymbol{n}}_1=\hat{\boldsymbol{n}}_1\\
    \bar{\boldsymbol{n}}_2=\hat{\boldsymbol{n}}_2}
} = \|\mathbf{G}\hat{\boldsymbol{n}}_2\|_2 = \sqrt{\lambda_i}
    \label{eq:det_optimal_mB_eigen}
\end{equation}
where the last equality comes from the explicit form of the second eigenvalue problem \eqref{eq:explicit_sol_mB_2_dec}, i.e., $\mathbf{G}^T\mathbf{G}\hat{\boldsymbol{n}}_2=\lambda_i\hat{\boldsymbol{n}}_2$. From equation \eqref{eq:det_optimal_mB_eigen} we can see that the solution associated to the largest eigenvalue $\lambda_2$ is the sought solution of problem \eqref{eq:constr_opt_mB_nodes}, while the solution associated to the smallest eigenvalue $\lambda_1$ is therefore a saddle point.

In $d=3$ dimensions there is a null eigenvalue since $\mathbf{G}$ also has one null eigenvalue. This implies that in the 3D case with $m_B=2$ boundary nodes, for each boundary node there is a normal which leads to a singular interpolation matrix regardless of the other normal. Nonetheless, the normals associated to null eigenvalues can not be solutions of equations \eqref{eq:explicit_sol_mB_2}, therefore even in the 3D case there are only two solutions, namely a local extremum and a saddle point. 

Unfortunately, when $m_B=3$ in $d=2$ dimensions, there can be two distinct local extrema in addition to saddle points. Indeed, if we consider the following particular symmetric choice for the vectors $\boldsymbol{g}_{i,j}$:
\begin{equation}
\begin{split}
    \boldsymbol{g}_{1,1}&=\boldsymbol{g}_{2,3}=\boldsymbol{g}_{3,2}=\{1,0\}^T\\
    \boldsymbol{g}_{2,2}&=\boldsymbol{g}_{1,3}=\boldsymbol{g}_{3,1}=\{0,0\}^T\\
    \boldsymbol{g}_{3,3}&=\boldsymbol{g}_{1,2}=\boldsymbol{g}_{2,1}=\{0,1\}^T
    \label{eq:particular_gij}
\end{split}
\end{equation}
then the determinant in equation \eqref{eq:Schur_Multiple_node_cardinal_det_Gbb} becomes:
\begin{equation}
\begin{split}
    \det(\mathbf{S}_{BB}) &=
\begin{vmatrix}
\boldsymbol{g}_{1,1}\cdot\bar{\boldsymbol{n}}_{1} & \boldsymbol{g}_{1,2}\cdot\bar{\boldsymbol{n}}_{1} &
\boldsymbol{g}_{1,3}\cdot\bar{\boldsymbol{n}}_{1} \\
\boldsymbol{g}_{2,1}\cdot\bar{\boldsymbol{n}}_{2} & \boldsymbol{g}_{2,2}\cdot\bar{\boldsymbol{n}}_{2} &
\boldsymbol{g}_{2,3}\cdot\bar{\boldsymbol{n}}_{2} \\
\boldsymbol{g}_{3,1}\cdot\bar{\boldsymbol{n}}_{3} & \boldsymbol{g}_{3,2}\cdot\bar{\boldsymbol{n}}_{3} &
\boldsymbol{g}_{3,3}\cdot\bar{\boldsymbol{n}}_{3}
\end{vmatrix}
=\\
&=\begin{vmatrix}
\cos\alpha_1 & \sin\alpha_1 & 0 \\
\sin\alpha_2 & 0              & \cos\alpha_2 \\
0 & \cos\alpha_3 & \sin\alpha_3
\end{vmatrix}
=\\   
& = -\big(\cos\alpha_1\cos\alpha_2\cos\alpha_3+\sin\alpha_1\sin\alpha_2\sin\alpha_3\big)
    \end{split}
    \label{eq:det_optimal_mB_3}
\end{equation}
where the 2D normals are parametrized as usual as  $\bar{\boldsymbol{n}}_{i}=\{\cos\alpha_i,\sin\alpha_i\}^T$.

From the last trigonometric expression in equation \eqref{eq:det_optimal_mB_3} it can be seen that the extremum $\alpha_1=\alpha_2=\alpha_3=0$, including the equivalent solutions obtained by reversing any of the normals, i.e., $\alpha_i=\pi$, is complemented by another distinct extremum $\alpha_1=\alpha_2=\alpha_3=\pi/2$, including the corresponding equivalent solutions.

It is therefore to be expected that multiple distinct extrema can also occur with $m_B>3$ boundary nodes, in both 2D and 3D cases. Nonetheless, we note that $\boldsymbol{g}_{i,j}$ are not arbitrary vectors but rather they are defined in equation \eqref{eq:Schur_Multiple_node_cardinal_H} according to the RBF interpolant, i.e., according to the relative positions of the stencil nodes. From a practical point of view, we think that a sufficiently regular distribution of nodes in the stencil can guarantee a unique extremum, as confirmed by numerical experiments.

\subsection{Computation of the optimal directions}\label{App_comput_optimal}
Equations \eqref{eq:ti_opt}-\eqref{eq:constr_opt_mB_nodes_solution} represent a nonlinear system of coupled algebraic equations in the $m_B$ unknown vectors $\hat{\boldsymbol{n}}_i$, which therefore requires an iterative solution. Such iterative process can be directly obtained from equation \eqref{eq:ti_opt} as follows:
\begin{equation}
    t_{i,\eta}^{(k+1)}=\det\!\big(H(\mathcal{G}_{BB},\hat{\mathcal{N}}_{i,\eta}^{(k)})\big)
    \label{eq:constr_opt_mB_nodes_solution_iter}
\end{equation}
where the superscript $(k)$ refers to iteration $k$. The values $t_{i,\eta}^{(k+1)}$ for $i=1,\ldots,m_B$ and $\eta=1,\ldots,d$ are therefore computed explicitly by using equation \eqref{eq:constr_opt_mB_nodes_solution_iter} where $\hat{\mathcal{N}}_{i,\eta}^{(k)}$ contains the normals at iteration $k$. Then the actual normals $\hat{\boldsymbol{n}}_i^{(k+1)}$ at iteration $k+1$ are obtained once again through the normalization expressed by equation \eqref{eq:constr_opt_mB_nodes_solution}.

It is convenient to compute the determinant in equation \eqref{eq:constr_opt_mB_nodes_solution_iter} as follows. Let us consider the explicit form of equation \eqref{eq:ti_opt} where both row $i$ and column $i$ of matrix $H(\mathcal{G}_{BB},\hat{\mathcal{N}}_{i,\eta})$ are moved to the last row and to the last column, respectively:
\begin{equation}
\begin{split}
    t_{i,\eta} &= \det\!\big(H(\mathcal{G}_{BB},\hat{\mathcal{N}}_{i,\eta})\big)=\\
    &=\begin{vmatrix}
    \begin{array}{ccc|c}
\boldsymbol{g}_{1,1  }\cdot\hat{\boldsymbol{n}}_{1} & \cdots & \boldsymbol{g}_{1,m_B}\cdot\hat{\boldsymbol{n}}_{1} & \boldsymbol{g}_{1,i  }\cdot\hat{\boldsymbol{n}}_{1} \\
\vdots & \ddots & \vdots & \vdots\\
\boldsymbol{g}_{m_B,1  }\cdot\hat{\boldsymbol{n}}_{m_B} & \cdots & \boldsymbol{g}_{m_B,m_B}\cdot\hat{\boldsymbol{n}}_{m_B} & \boldsymbol{g}_{m_B,i  }\cdot\hat{\boldsymbol{n}}_{m_B}\\[.25em] \hline
\boldsymbol{g}_{i,1  }\cdot\boldsymbol{e}_\eta & \cdots & \boldsymbol{g}_{i,m_B}\cdot\boldsymbol{e}_\eta & \boldsymbol{g}_{i,i  }\cdot\boldsymbol{e}_\eta
\end{array}
    \end{vmatrix}
    =
    \begin{vmatrix}
    \begin{array}{c|c}
\mathbf{S}_{i} & \boldsymbol{c}_{i} \\ \hline
\boldsymbol{r}_{i,\eta} & S_{ii,\eta}
    \end{array}
    \end{vmatrix}
    \label{eq:det_i_row_col}
\end{split}
\end{equation}

In other words: starting from matrix $H(\mathcal{G}_{BB},\hat{\mathcal{N}}_{i,\eta})$, matrix $\mathbf{S}_{i}$ is obtained by removing row $i$ and column $i$, column vector $\boldsymbol{c}_{i}$ is obtained by removing entry $i$ of column $i$, and row vector $\boldsymbol{r}_{i,\eta}$ is obtained by removing entry $i$ of row $i$.

By the properties of the Schur complement of the block $\mathbf{S}_{i}$ of the last matrix in equation \eqref{eq:det_i_row_col}, we have:
\begin{equation}
    t_{i,\eta} =
    \begin{vmatrix}
    \begin{array}{c|c}
\mathbf{S}_{i} & \boldsymbol{c}_{i} \\ \hline
\boldsymbol{r}_{i,\eta} & S_{ii,\eta}
    \end{array}
    \end{vmatrix}=\det(\mathbf{S}_{i})\big(S_{ii,\eta}-\boldsymbol{r}_{i,\eta}\mathbf{S}_{i}^{-1}\boldsymbol{c}_{i}\big)
    \label{eq:det_i_row_col_Schur}
\end{equation}

Since equation \eqref{eq:det_i_row_col_Schur} holds for each of the $d$ components of the vector $\boldsymbol{t}_i=\{t_{i,\eta}\}_{\eta=1}^d$, it can be written explicitly as follows:
\begin{equation}
    \boldsymbol{t}_i = \det(\mathbf{S}_{i})
    \Big(\boldsymbol{g}_{i,i}-\sum_{j=1}^{m_B-1}w_j\boldsymbol{g}_{i,j'}\Big)
    \label{eq:ti_expl}
\end{equation}
where $w_j$ are the $m_B-1$ entries of the column vector $\boldsymbol{w}=\mathbf{S}_{i}^{-1}\boldsymbol{c}_{i}$ and $j'$ is used to index vectors $\boldsymbol{g}_{i,j}$ excluding the case $j=i$, i.e., $j'=j$ if $j<i$, $j'=j+1$ otherwise. We finally note that $\det(\mathbf{S}_{i})$ in equation \eqref{eq:ti_expl} does not need to be computed because of the subsequent normalization of equation \eqref{eq:constr_opt_mB_nodes_solution}.

Moving back to the required iterative scheme formerly expressed by equation \eqref{eq:constr_opt_mB_nodes_solution_iter}, it can finally be expressed by equation \eqref{eq:ti_expl} where the dependence of the right side upon the iteration index $k$ arises only in the terms $w_j=w_j^{(k)}$. We observe that the vector $\boldsymbol{w}$ carrying the values $w_j$ does not depend upon any specific component $\eta$, and therefore it can be computed only once at each iteration for each of the $m_B$ vectors $\boldsymbol{t}_i$.

Since the system of equations \eqref{eq:ti_opt}-\eqref{eq:constr_opt_mB_nodes_solution} is nonlinear, it is important to choose a proper set of initial vectors $\hat{\boldsymbol{n}}_i^{(0)}$ for the iterative scheme \eqref{eq:constr_opt_mB_nodes_solution_iter} to converge to the correct solution. Numerical experiments showed that a good starting point can be obtained from the diagonal vectors of $d$-matrix $\mathcal{G}_{BB}$ defined in equation \eqref{eq:Schur_Multiple_node_cardinal_H}:
\begin{equation}
    \hat{\boldsymbol{n}}_i^{(0)}=
    \frac{\boldsymbol{g}_{ii}}{\|\boldsymbol{g}_{ii}\|_2}
    \label{eq:diag_vectors}
\end{equation}
for $i=1,\ldots,m_B$. We note that this choice corresponds exactly to equation \eqref{eq:ti_expl} with $w_j=0$. Equation \eqref{eq:ti_expl} can therefore be employed to provide an improved starting point if any suitable approximation of $w_j$ is available. Nonetheless, since in general the norm of $\boldsymbol{g}_{ii}$ and $\boldsymbol{g}_{ij'}$ are in the same order of magnitude, naive approximations of $w_j\ne 0$ can lead to large errors.

Another good starting point can be obtained from the following geometric heuristic:
\begin{equation}
    \hat{\boldsymbol{n}}_i^{(0)}=
    \frac{\boldsymbol{x}_{m_I+i}-\bar{\boldsymbol{x}}}{\|\boldsymbol{x}_{m_I+i}-\bar{\boldsymbol{x}}\|_2}
    \label{eq:geom_vectors}
\end{equation}
where $\bar{\boldsymbol{x}}$ is some geometrical reference point for the stencil, e.g., the centroid of the internal nodes.

\subsection{Computational costs}\label{App_comput_cost}
The computation of the optimal directions consists of two phases: the computation of the $d$-matrix $\mathcal{G}_{BB}$, equation \eqref{eq:Schur_Multiple_node_cardinal_H}, and the iterative solution of the associated system, equation \eqref{eq:constr_opt_mB_nodes_solution_iter}. The computation of $\mathcal{G}_{BB}$ requires the computation of the $m_I\times m_B$ matrix $\bar{\boldsymbol{\psi}}$ from equation \eqref{eq:OneBoundaryNode_psi_at_boundary}, whose cost is therefore $\mathcal{O}(m_I^3)+\mathcal{O}(m_B\cdot m_I^2)$, i.e., factorization of $\boldsymbol{\varphi}_{II}$ + solution of $m_B$ associated linear systems. Then, each iteration updates the $m_B$ optimal directions through equation \eqref{eq:ti_expl}, each of which requires the $(m_B-1)\times (m_B-1)$ linear system $\mathbf{S}_{i}\boldsymbol{w}=\boldsymbol{c}_{i}$ to be solved for the weights $\boldsymbol{w}$, thus leading to a cost $\mathcal{O}(N_{it}\cdot m_B^4)$, where $N_{it}$ is the number of iterations which depends upon the desired level of convergence. We note that the optimal directions do not need to be computed with high accuracy since they are only used in equation \eqref{eq:dot_product}.

From a practical point of view, our numerical experiments showed that most of the computational time is due to the iterative phase, even if the number of boundary nodes $m_B$ is usually much smaller than the number of internal nodes $m_I$ of a single stencil. Nonetheless, we trivially note that the computation of the optimal directions is required only for those stencils having boundary nodes, and the number of these boundary stencils is usually much smaller than the total number of stencils because of obvious geometrical reasons.

We finally note that the iterative strategy presented as approach 1 in section \ref{sss:dot_product} requires a single calculation of the $d$-matrix $\mathcal{G}_{BB}$, since the successive elimination of boundary nodes corresponds to the removal of the corresponding rows and columns of $\mathcal{G}_{BB}$.

\section{Optimal position for boundary nodes}\label{App:s3}

\subsection{Optimization process}
We are interested in improving the interpolation properties of the reference stencil of Figure \ref{FIG:ref_stencil}. First of all, let's consider that moving the inner nodes would probably affect the interpolation properties of the neighboring stencils. Indeed, if we decided to move the inner nodes of one stencil, in order to preserve the local uniformity of the node distribution, which ensures that the difference between fill distance and separation distance remains negligible, we would need to adjust the position of many other nodes accordingly.
On top of that, if we take into account the matrix $\mathbf{M}$ in the general case of equation \eqref{eq:InterpM_pureRBF_MultipleBoundaryNode_splitting} and consider that $\det(\boldsymbol{\varphi}_{II}) \neq 0$ in the cases of interest, we are once again prompted to focus our attention on the boundary nodes.

Since the the interpolation becomes unstable as we approach those configurations where $\det(\mathbf{S}_{BB}) = 0$, we might consider the problem of maximizing $\det(\mathbf{M})$ or $\det(\mathbf{S}_{BB})$ by controlling the locations of boundary nodes.
Unfortunately, we would find that, because of the radial basis functions, the further apart we spread the boundary nodes, the larger the determinant becomes.

In section \ref{ss:card_func} the Lebesgue constants were pointed as an upper bound for the interpolation error, it was also shown that an ill conditioned interpolation matrix $\mathbf{M}$ would lead to very high Lebesgue constants.
Starting from these remarks it seems natural to minimize the Lebesgue functions (as defined in equation \eqref{eq:LebFunc}) by controlling the position of boundary nodes.
Here follows a brief discussion of the method adopted for this purpose.

When no polynomial augmentation is used, cardinal functions $\psi_i(\boldsymbol{x})$ can be found as the solution of the following linear system:
\begin{equation}
    \mathbf{M}^T \{\psi_i(\bx)\}_1^m = \{\varphi_i(\boldsymbol{x})\}_1^m
    \label{eq:card_no_poly}
\end{equation}
where the same notation of equation \eqref{eq:adjoint_cardinal} was used.

In some cases $\lambda_I$ and $\lambda_B$ exhibit large values at the boundary when Neumann boundary conditions are employed, see for example Figure \ref{fig:cardinal_alpha}.

Willing to find the configuration of boundary nodes that minimizes the Lebesgue constants $\Lambda_I$ and $\Lambda_B$, as defined in equation \eqref{eq:LebConst_IB}, the cost function is defined as follows:
\begin{equation}
    \mathcal{F}(\boldsymbol{x}_{m_J}, \dots, \boldsymbol{x}_m) =
    \sum_{k = m_J}^m \lambda_I(\boldsymbol{x}_k) =
    \sum_{k = m_J}^m \sum_{i = 1}^{m_I} |\psi_i(\boldsymbol{x}_k)|
    \label{eq:cost_function}
\end{equation}
where the stencil nodes are enumerated as in section \ref{ss:RBF_FD}. In equation \eqref{eq:cost_function} we considered only internal Lebesgue functions $\lambda_I$, however, it was observed that optimizing with respect to the boundary Lebesgue functions $\lambda_B$ leads to very similar optimal configurations.

By using a more compact vector notation, the sensitivities of the cost function with respect to the position of a boundary node $\boldsymbol{x}_j=\{x_{j,1}\, , \, x_{j,2}\}$ in 2D are:
\begin{equation}
    \frac{\partial \mathcal{F}}{\partial x_{j,\eta}}
    =
    \sum_{k=m_J}^m \sgn\big(\boldsymbol{\psi}_I(\boldsymbol{x}_k)\big)^T 
    \frac{\partial \boldsymbol{\psi}(\boldsymbol{x}_k)}{\partial x_{j,\eta}}\\
    \label{eq:sensitivities1}
\end{equation}
where $\boldsymbol{\psi}(\boldsymbol{x}_k)$ is the column vector $\{\psi_i(\boldsymbol{x}_k)\}_{i=1}^{m}$ and $\sgn(\boldsymbol{\psi}_I(\boldsymbol{x}_k))$ is the column vector having the signs of $\boldsymbol{\psi}(\boldsymbol{x}_k)$ in the first $m_I$ entries, and zeros for the last $m_B$ indices corresponding to the boundary nodes.

The term $\frac{\partial \boldsymbol{\psi}(\boldsymbol{x}_k)}{\partial x_{j,\eta}}$ in equation \eqref{eq:sensitivities1} can be obtained by deriving equation \eqref{eq:card_no_poly} as follows:
\begin{equation}
    \frac{\partial \boldsymbol{\psi}(\boldsymbol{x}_k)}{\partial x_{j,\eta}} = (\mathbf{M}^T)^{-1} \left(\frac{\partial \boldsymbol{\varphi}(\boldsymbol{x}_k)}{\partial x_{j,\eta}} -
     \frac{\partial \mathbf{M}^T}{\partial x_{j,\eta}} \, \boldsymbol{\psi}(\boldsymbol{x}_k)
    \right)
    \label{eq:partial_derivatives}
\end{equation}
where $\frac{\partial \mathbf{M}^T}{\partial x_{j,\eta}}$ is obtained by taking the derivative of all the matrix entries.
If we substitute equation \eqref{eq:partial_derivatives} into equation \eqref{eq:sensitivities1}, the following term appears:
\begin{equation}
    \sgn\big(\boldsymbol{\psi}_I(\boldsymbol{x}_k)\big)^T (\mathbf{M}^T)^{-1} = \Big(\mathbf{M}^{-1} \sgn\big(\boldsymbol{\psi}_I(\boldsymbol{x}_k)\big)\Big)^T = \boldsymbol{c}^T
\end{equation}
where $\boldsymbol{c}$ is therefore the column vector given as the solution of the following adjoint equation:
\begin{equation}
    \mathbf{M} \, \boldsymbol{c} = \sgn\big(\boldsymbol{\psi}_I(\boldsymbol{x}_k)\big)
    \label{eq:adjoint_sensitivities}
\end{equation}

The computational cost for the calculation of the sensitivities of \eqref{eq:sensitivities1} then becomes independent of the dimension of the problem and requires the solution of $m_B$ linear systems like equation \eqref{eq:adjoint_sensitivities}, one for each boundary node.

The sensitivities can then be used to perform a gradient based optimization of the boundary nodes positions, this is done by taking the components orthogonal to the normals.

Whenever two nodes ended up overlapping during the optimization, one of them was removed, this is why in Figure \ref{fig:optimal_bnd} only $6$ boundary nodes appear.

\subsection{Results}
\begin{figure}[t]
    \def\SpaceBelowText{.27em}
    \def\ImagesWidth{.32\textwidth}
    \setlength\fboxsep{0pt}
    \setlength\fboxrule{0.5pt}
    \centering
    \begin{minipage}[b]{\textwidth}
        \begin{subfigure}[b]{\ImagesWidth}
            \centering
            $\alpha=-\pi/12$\\
            \includegraphics[width = \textwidth]{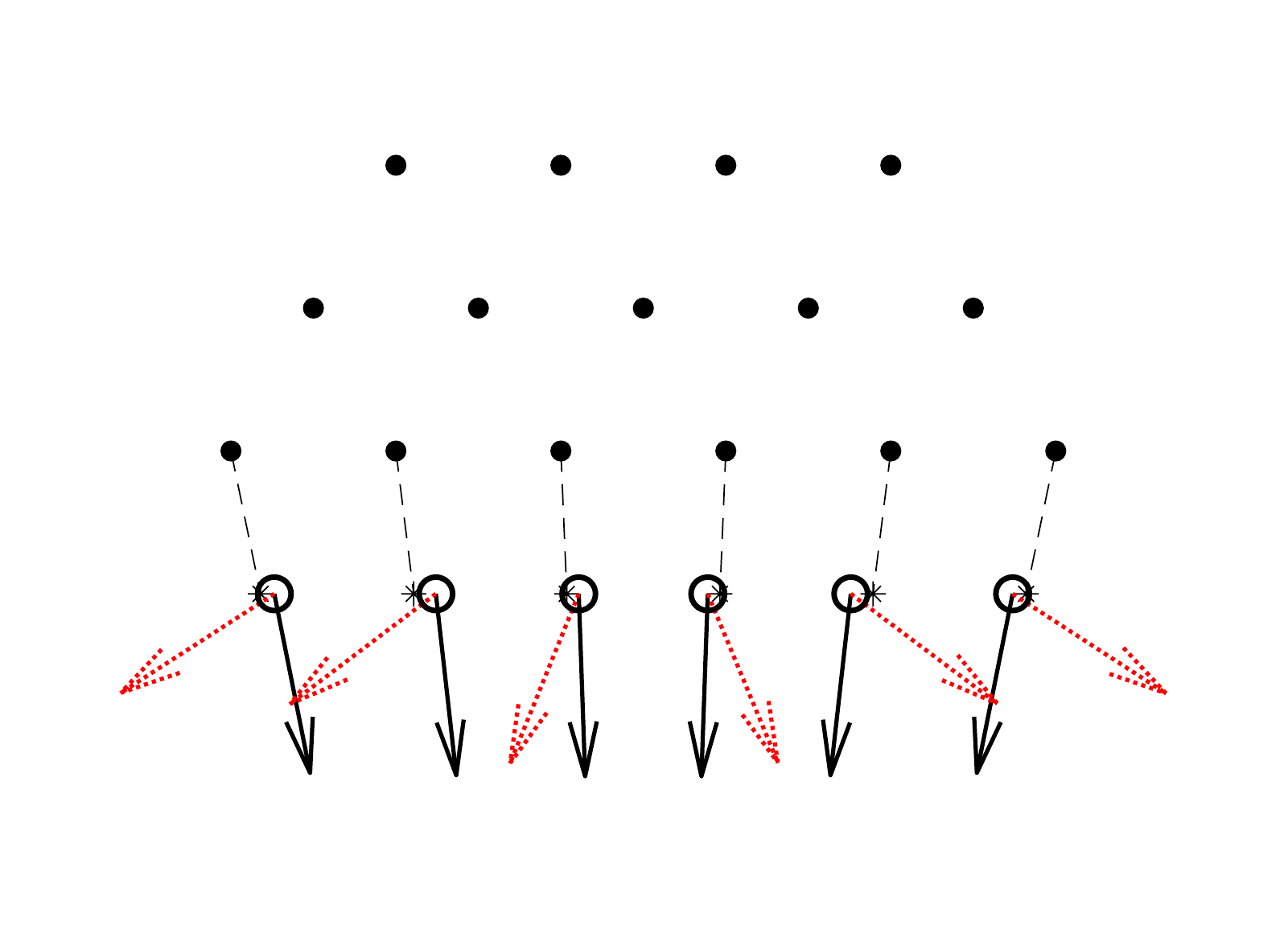}
        \end{subfigure}
        \begin{subfigure}[b]{\ImagesWidth}
            \centering
            $\alpha=0$\\[\SpaceBelowText]
            \includegraphics[width = \textwidth]{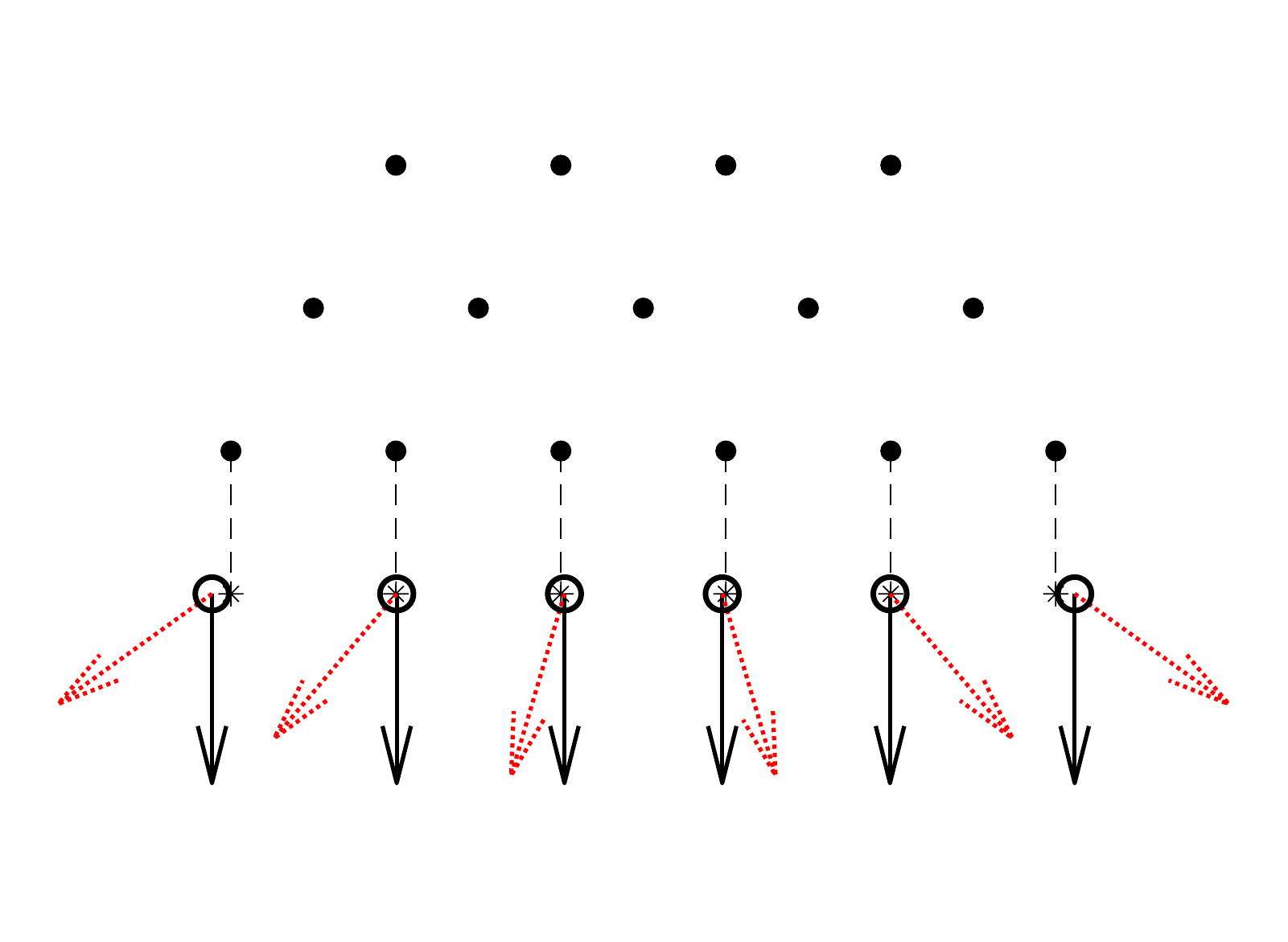}
        \end{subfigure}
        \begin{subfigure}[b]{\ImagesWidth}
            \centering
            $\alpha=\pi/12$\\
            \includegraphics[width = \textwidth]{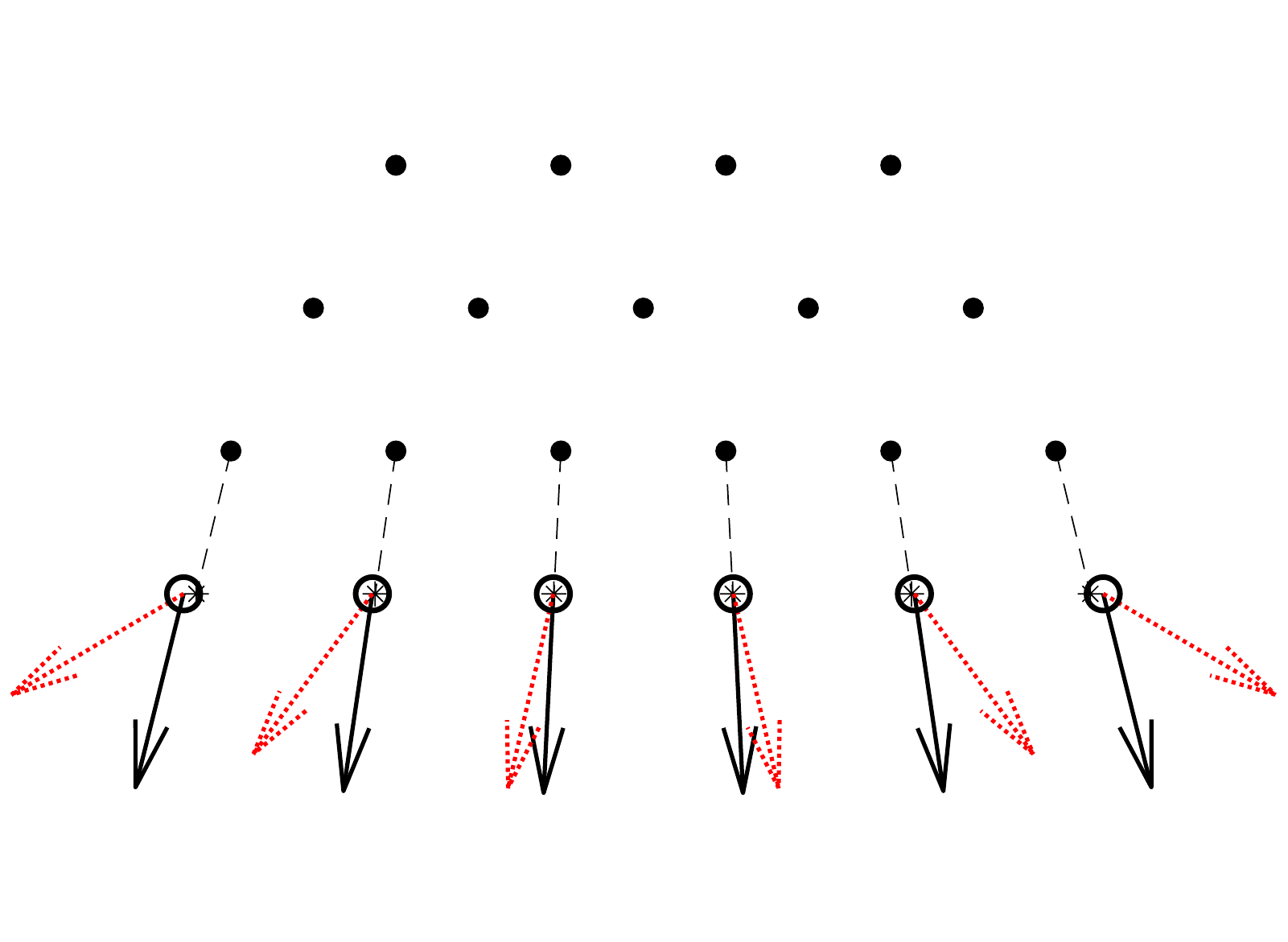}
        \end{subfigure}
    \end{minipage}
    \caption{optimal placement of boundary nodes obtained starting from the reference stencil of Figure \ref{FIG:ref_stencil}. Black solid lines are the actual normals, red dotted lines are the optimal directions, asterisks are the position of boundary nodes obtained by projecting the first layer of inner nodes along the normals.}
    \label{fig:optimal_bnd}
\end{figure}

We have performed the optimization process described above by starting from the reference stencil of Figure
\ref{FIG:ref_stencil} and using MQ RBF with shape parameter $\varepsilon$ satisfying
$\varepsilon s = 0.5$, where $s$ is the distance between the nodes. Once the convergence is reached, we are left with the situation depicted in Figure \ref{fig:optimal_bnd}. We can see that different values of the angle $\alpha$ correspond to different optimal placements of the boundary nodes.

We remark that, even if the procedure described so far gives the best placement for the boundary nodes, this comes at a relatively high computational cost. Indeed, it is required to solve $m_B$ linear systems like \eqref{eq:card_no_poly} in order to calculate the cardinal functions, in addition to $m_B$ linear systems like \eqref{eq:adjoint_sensitivities} in order to estimate the sensitivities at each step of the optimization. On top of that, a control on the normals, like the one explained in section \ref{sss:dot_product}, is still needed when the boundary has high curvature, see the angle between actual and optimal normals for the case $\alpha = -\pi/12$ in Figure \ref{fig:optimal_bnd}.

The results obtained here suggest that the interpolation can be improved by placing the boundary nodes on the locations identified by projecting along the normals the inner nodes of the first layer. Replacing standard boundary nodes with the projected ones can be done after the node generation phase and requires minimal additional computational cost on smooth and regular domains.

\bibliographystyle{elsarticle-num} 
\bibliography{references}

\end{document}